\documentclass[journal,11pt,draftcls,onecolumn]{IEEEtran}
\newcommand{\DRAFT}[1]{#1}
\newcommand{\FINAL}[1]{}

\usepackage{url}
\usepackage{here}
\usepackage[noadjust]{cite}      
\usepackage{graphicx}  
\usepackage{subfigure} 
\usepackage{stfloats}  
\usepackage{amsmath}   
\usepackage{xspace,bbm,epic,eepic,amssymb}
\usepackage[latin1]{inputenc}
\usepackage{theorem}

\usepackage[lined,boxed,linesnumbered]{algorithm2e}

\def\XS{\xspace}
\DeclareMathAlphabet{\mathb}{OML}{cmm}{b}{it}
\def\sbm#1{\ensuremath{\mathb{#1}}}                
\def\scu#1{\ensuremath{\mathcal{#1\XS}}}      
\def\sbl#1{\ensuremath{\mathbbm{#1}}}         
\def\sbmm#1{\ensuremath{\boldsymbol{#1}}}  
\def\Ec{{\scu{E}}\XS}
\def\Jc{{\scu{J}}\XS}   
\def\Sc{{\scu{S}}\XS}   

\def\Rbb{{\sbl{R}}\XS}

\def\Ab{{\sbm{A}}\XS}  
\def\ab{{\sbm{a}}\XS}
\def\hb{{\sbm{h}}\XS}
\def\nb{{\sbm{n}}\XS}
\def\xb{{\sbm{x}}\XS}
\def\yb{{\sbm{y}}\XS}

\def\lambdab {{\sbmm{\lambda}}\XS}
\RequirePackage{ifthen}
\RequirePackage{calc}

\newcommand{\taille}[1][\scad]{%
\ifthenelse{#1 = -5}{}{}%
\ifthenelse{#1 = -4}{\tiny}{}%
\ifthenelse{#1 = -3}{\scriptsize}{}%
\ifthenelse{#1 = -2}{\footnotesize}{}%
\ifthenelse{#1 = -1}{\small}{}%
\ifthenelse{#1 = 0}{\normalsize}{}%
\ifthenelse{#1 = 1}{\large}{}%
\ifthenelse{#1 = 2}{\Large}{}%
\ifthenelse{#1 = 3}{\LARGE}{}%
\ifthenelse{#1 = 4}{\huge}{}%
\ifthenelse{#1 = 5}{\Huge}{}}

\newboolean{@serif}
\setboolean{@serif}{true}
\def\scad{-5} 
\newcounter{taille}
\newcommand{\sca}[2][\scad]{\setcounter{taille}{#1}%
  \ifthenelse{\boolean{@serif}}
  {{\taille[\thetaille]\textsc{#2}}}
  {\setcounter{taille}{\value{taille}-1}{\uppercase{\taille[\thetaille]#2}}}}

\def\apost{\textit{a posteriori}\XS}
\def\eg{\textit{e.g.,}\XS}
\def\ie{\textit{i.e.,}\XS}
\def\iid{i.i.d.\XS}
\def\vs{\textit{vs}\XS}

\newcommand{\figc}[2][]
   {\setlength{\tabcolsep}{0pt}\begin{tabular}{c}\includegraphics[#1]{#2}\end{tabular}}

\def\stdpth#1{(#1)}
\def\stdacc#1{\{#1\}}
\def\stdcro#1{[#1]}
\def\stdbars#1{|#1|}
\def\stdscal#1{\langle#1\rangle}
\def\bigacc#1{\bigl\{#1\bigr\}}            
\def\Bigacc#1{\Bigl\{#1\Bigr\}}   
\def\argmax{\mathop{\mathrm{arg\,max}}} 
\def\argmin{\mathop{\mathrm{arg\,min}}} 
\def\froc#1#2{{#1/#2}}                  
\newtheorem{defi}{Definition}
\newtheorem{lemma}{Lemma}
\newtheorem{theorem}{Theorem}

\graphicspath{{./}}
\setboolean{@serif}{false}

\begin{document}
\title{Homotopy based algorithms for $\ell_0$-regularized least-squares}
\author{
  Charles~Soussen$^\star$,~%
  J\'er\^ome~Idier,~
  Junbo~Duan,~
  and~David~Brie
 \thanks{This work was carried out in part while C. Soussen was visiting 
 IRCCyN during the academic year 2010-2011 with the financial support 
 of CNRS.}
  \thanks{C.~Soussen and
  D.~Brie are with the Universit\'e de Lorraine and CNRS 
  at the Centre de Recherche en Automatique de Nancy (UMR 7039)%
  \FINAL{, F-54506 Vand{\oe}uvre-l\`es-Nancy, France}.%
  \DRAFT{Campus Sciences, B.P. 70239, F-54506 Vand{\oe}uvre-l\`es-Nancy, 
France. Tel: (+33)-3 83 59 56 43, Fax: (+33)-3 83 68 44
    62.} E-mail:
  charles.soussen@univ-lorraine.fr, david.brie@univ-lorraine.fr.}
\thanks{J.~Idier is with L'UNAM Universit\'e, Ecole Centrale Nantes 
    and CNRS at the Institut de Recherche en Communications
    et Cybern\'etique de Nantes (UMR 6597), 
    \DRAFT{1 rue de la No\"e, BP 92101,} F-44321 Nantes\DRAFT{ Cedex~3}, France. \DRAFT{Tel: (+33)-2 40 37 69 09, Fax: (+33)-2 40
    37 69 30.} E-mail: jerome.idier@irccyn.ec-nantes.fr.}
\thanks{J.~Duan was with CRAN. He is now with the Department
  of Biomedical Engineering, Xi'an Jiaotong University%
  \FINAL{, Xi'An 710049, China.}%
  \DRAFT{. No. 28, Xianning West Road, Xi'an 710049, Shaanxi Province, China.
Tel: (+86)-29-82 66 86 68, Fax: (+86)-29 82 66 76 67.}
 E-mail: junbo.duan@mail.xjtu.edu.cn. }
}
\markboth{Technical Report%
\FINAL{, \today}}{Soussen, Idier, Duan, Brie: Using the Document Class IEEEtran.cls}
\maketitle

\begin{abstract}
  Sparse signal restoration is usually formulated as the minimization
  of a quadratic cost function $\|\yb-\Ab\xb\|_2^2$, where \Ab is a
  dictionary and \xb is an unknown sparse vector. It is well-known
  that imposing an $\ell_0$ constraint leads to an NP-hard
  minimization problem. The convex relaxation approach has received
  considerable attention, where the $\ell_0$-norm is replaced by the
  $\ell_1$-norm. Among the many efficient $\ell_1$ solvers, the
  homotopy algorithm minimizes $\|\yb-\Ab\xb\|_2^2+\lambda\|\xb\|_1$
  with respect to \xb for a continuum of $\lambda$'s. It is inspired
  by the piecewise regularity of the $\ell_1$-regularization path,
  also referred to as the homotopy path. In this paper, we address the
  minimization problem $\|\yb-\Ab\xb\|_2^2+\lambda\|\xb\|_0$ for a
  continuum of $\lambda$'s and propose two heuristic search algorithms
  for $\ell_0$-homotopy. Continuation Single Best Replacement is a
  forward-backward greedy strategy extending the Single Best
  Replacement algorithm, previously proposed for $\ell_0$-minimization
  at a given $\lambda$. The adaptive search of the $\lambda$-values is
  inspired by $\ell_1$-homotopy. $\ell_0$ Regularization Path Descent
  is a more complex algorithm exploiting the structural properties of
  the $\ell_0$-regularization path, which is piecewise constant with
  respect to $\lambda$. Both algorithms are empirically evaluated for
  difficult inverse problems involving ill-conditioned
  dictionaries. Finally, we show that they can be easily coupled with
  usual methods of model order selection.
\end{abstract}

\section{Introduction}
\label{sec:intro}
Sparse approximation from noisy data is traditionally addressed as
the constrained least-square problems
\begin{align}
\min_\xb\|\yb-\Ab\xb\|_2^2~~\textrm{subject to}~~\|\xb\|_0\leq
k
\label{eq:cstr_LS}
\end{align}
or
\begin{align}
\min_\xb\|\xb\|_0~~\textrm{subject to}~~\|\yb-\Ab\xb\|_2^2\leq
\varepsilon
\label{eq:cstr_LS2}
\end{align}
where $\|\xb\|_0$ is the $\ell_0$-``norm'' counting the number of
nonzero entries in \xb, and the quadratic fidelity-to-data term
$\|\yb-\Ab\xb\|_2^2$ measures the quality of approximation.
Formulation~\eqref{eq:cstr_LS} is well adapted when one has a
knowledge of the maximum number $k$ of atoms to be selected in the
dictionary \Ab. On the contrary, the choice of~\eqref{eq:cstr_LS2} is
more appropriate when $k$ is unknown but one has a knowledge of the
variance of the observation noise.  The value of $\varepsilon$ may
then be chosen relative to the noise variance. Since
both~\eqref{eq:cstr_LS} and~\eqref{eq:cstr_LS2} are subset selection
problems, they are discrete optimization problems. They are known to
be NP-hard except for specific cases~\cite{Natarajan95}.

When no knowledge is available on either $k$ or $\varepsilon$, the
unconstrained formulation
\begin{align}
\min_\xb\{\Jc(\xb;\lambda)=\|\yb-\Ab\xb\|_2^2+\lambda\|\xb\|_0\}
\label{eq:pen_LS}
\end{align}
is worth being considered, where $\lambda$ expresses the trade-off
between the quality of approximation and the sparsity
level~\cite{Nikolova13}. In a Bayesian viewpoint, \eqref{eq:pen_LS}
can be seen as a (limit) maximum \apost formulation where
$\|\yb-\Ab\xb\|_2^2$ and the penalty $\|\xb\|_0$ are respectively
related to a Gaussian noise distribution and a prior distribution for
sparse signals (a limit Bernoulli-Gaussian distribution
with infinite Gaussian variance)~\cite{Soussen11c}.

\subsection{Classification of methods}
\subsubsection{$\ell_0$-constrained least-squares}
The discrete algorithms dedicated to
problems~\eqref{eq:cstr_LS}-\eqref{eq:cstr_LS2} can be categorized
into two classes. First, the forward greedy algorithms explore subsets
of increasing cardinalities starting from the empty set. At each
iteration, a new atom is appended to the current subset, therefore
gradually improving the quality of
approximation~\cite{Tropp10}. Greedy algorithms include, by increasing
order of complexity: Matching Pursuit (MP)~\cite{Mallat93}, Orthogonal
Matching Pursuit (OMP)~\cite{Pati93}, and Orthogonal Least Squares
(OLS)~\cite{Chen89}, also referred to as forward selection in
statistical regression~\cite{Miller02} and known as Order Recursive
Matching Pursuit (ORMP)~\cite{Cotter99} and Optimized Orthogonal
Matching Pursuit (OOMP)~\cite{RebolloNeira02}. The second category are
thresholding algorithms, where each iteration delivers a subset of
same cardinality $k$. Popular thresholding algorithms include
Iterative Hard Thresholding~\cite{Blumensath08a}, Subspace
Pursuit~\cite{Dai09} and CoSaMP~\cite{Needell09}.

Among these two categories, greedy algorithms are well-adapted to the
resolution of~\eqref{eq:cstr_LS} and~\eqref{eq:cstr_LS2} for
\emph{variable} sparsity levels. Indeed, they yield a series of
subsets for consecutive $k$ (\ie for decreasing approximation errors
$\varepsilon$) since at each iteration, the current subset is
increased by one element.

\subsubsection{$\ell_0$-penalized least-squares}
In~\cite{Soussen11c}, we evidenced that the minimization of
$\Jc(\xb;\lambda)$ using a descent algorithm leads to bidirectional
extensions of forward (orthogonal) greedy algorithms. To be more
specific, consider a candidate subset $S$ corresponding to the support
of \xb. Including a new element into $S$ yields a decrease of the
square error, defined as the minimum of $\|\yb-\Ab\xb\|_2^2$ for \xb
supported by $S$. On the other hand, the penalty term
$\lambda\|\xb\|_0$ is increased by $\lambda$. Overall, the cost
function $\Jc(\xb;\lambda)$ decreases as soon as the square error
variation exceeds $\lambda$. Similarly, a decrease of
$\Jc(\xb;\lambda)$ occurs when an element is removed from $S$ provided
that the squared error increment is lower than $\lambda$. Because both
inclusion and removal operations can induce a decrease of \Jc, the
formulation~\eqref{eq:pen_LS} allows one to design descent schemes
allowing a ``forward-backward'' search strategy, where each iteration
either selects a new atom (forward selection) or de-selects an atom
that was previously selected (backward elimination). The Bayesian
OMP~\cite{Herzet14} and Single Best Replacement
(SBR)~\cite{Soussen11c} algorithms have been proposed in this
spirit. They are extensions of OMP and OLS, respectively. Their
advantage over forward greedy algorithms is that an early wrong atom
selection may be later cancelled. Forward-backward algorithms include
the so-called stepwise regression algorithms which are OLS
extensions~\cite{Efroymson60,Berk80,Miller02}, and OMP based
algorithms of lower complexity~\cite{Zhang11,Herzet14}.

\subsubsection{Connection with the continuous relaxation of the $\ell_0$ norm}
The algorithms described so far are discrete search strategies
dedicated to $\ell_0$-regularized least-squares. A classical
alternative consists in relaxing the $\ell_0$-norm by a continuous
function that is nondifferentiable at 0, and optimizing the resulting
cost function. See, \eg \cite{Figueiredo07,Zibulevsky10} and
\cite{Candes08,Gasso09,Mourad10,Wipf10,Gholami11,Ramirez12,HoaiAn13,Selesnick14}
for convex ($\ell_1$) and nonconvex relaxation, respectively. The
convex problem $\min_\xb\|\yb-\Ab\xb\|_2^2$ s.t. $\|\xb\|_1\leq t$ is
referred to as both Basis Pursuit Denoising (BPDN) and the LASSO.  It
is noticeable that BPDN leads to stepwise
algorithms~\cite{Figueiredo07,Donoho08} including the popular
$\ell_1$-homotopy~\cite{Osborne00,Efron04,Donoho08}, a
forward-backward greedy search whose complexity is close to that of
OMP. $\ell_1$-homotopy is closely connected to the Least Angle
Regression (LARS), a simpler forward strategy allowing only atom
selections. It is referred to as ``LARS with the LASSO modification''
in~\cite{Efron04}. Importantly, $\ell_1$-homotopy solves the BPDN for
a continuum of values of $t$.

\subsection{Main idea}
\label{sec:main_ideas}
Our approach is dedicated to $\ell_0$-penalized least-squares.
It is based on the following geometrical interpretation.

\begin{figure}[t]
\centering
{
\begin{tabular}{c}
\setlength{\unitlength}{0.00043745in}
\begingroup\makeatletter\ifx\SetFigFont\undefined%
\gdef\SetFigFont#1#2#3#4#5{%
  \reset@font\fontsize{#1}{#2pt}%
  \fontfamily{#3}\fontseries{#4}\fontshape{#5}%
  \selectfont}%
\fi\endgroup%
{\renewcommand{\dashlinestretch}{30}
\begin{picture}(7000,4882)(-300,300)
\thicklines
\put(3920,3525){\blacken\ellipse{90}{90}}
\put(3920,3525){\ellipse{90}{90}}
\put(1770,2250){\blacken\ellipse{90}{90}}
\put(1770,2250){\ellipse{90}{90}}
\blacken\thicklines
\path(1095.000,4845.000)(1050.000,4995.000)(1005.000,4845.000)(1095.000,4845.000)
\path(1050,4995)(1050,775)(6700,775)
\blacken\path(6550.000,730.000)(6700.000,775.000)(6550.000,820.000)(6550.000,730.000)
\thinlines
\dottedline{60.000}(1050,3020)(6495,4460)
\dashline{60.000}(1770,2250)(1770,775)
\dashline{60.000}(3920,3525)(3920,775)
\dottedline{60.000}(1050,975)(3275,4917)
\dottedline{60.000}(1050,1822)(6277,4917) 
\thicklines
\path(1050,975)(1770,2250)(3920,3525)(6446,3525)
\thinlines
\dottedline{60.000}(1050,3525)(6446,3525)
\dottedline{60.000}(1050,2683)(4821,4915) 
\put(6420,415){\makebox(0,0)[lb]{\smash{{\SetFigFont{8}{8}{\rmdefault}{\mddefault}{\updefault}$\lambda$}}}}
\put(2600,4700){\makebox(0,0)[lb]{\smash{{\SetFigFont{8}{8}{\rmdefault}{\mddefault}{\updefault}$S^{\star}_2$}}}}
\put(250,3500){\makebox(0,0)[lb]{\smash{{\SetFigFont{8}{8}{\rmdefault}{\mddefault}{\updefault}$\Ec(S^{\star}_0)$}}}}
\put(250,1690){\makebox(0,0)[lb]{\smash{{\SetFigFont{8}{8}{\rmdefault}{\mddefault}{\updefault}$\Ec(S^{\star}_1)$}}}}
\put(-450,4770){\makebox(0,0)[lb]{\smash{{\SetFigFont{8}{8}{\rmdefault}{\mddefault}{\updefault}$\Ec(S)+\lambda\stdbars{S}$}}}}
\put(1680,415){\makebox(0,0)[lb]{\smash{{\SetFigFont{8}{8}{\rmdefault}{\mddefault}{\updefault}$\lambda^\star_2$}}}}
\put(3810,415){\makebox(0,0)[lb]{\smash{{\SetFigFont{8}{8}{\rmdefault}{\mddefault}{\updefault}$\lambda^\star_1$}}}}
\put(250,930){\makebox(0,0)[lb]{\smash{{\SetFigFont{8}{8}{\rmdefault}{\mddefault}{\updefault}$\Ec(S^{\star}_2)$}}}}
\put(5330,4700){\makebox(0,0)[lb]{\smash{{\SetFigFont{8}{8}{\rmdefault}{\mddefault}{\updefault}$S^{\star}_1$}}}}
\put(5550,3150){\makebox(0,0)[lb]{\smash{{\SetFigFont{8}{8}{\rmdefault}{\mddefault}{\updefault}$S^{\star}_0=\emptyset$}}}}
\put(525,415){\makebox(0,0)[lb]{\smash{{\SetFigFont{8}{8}{\rmdefault}{\mddefault}{\updefault}$\lambda^\star_3=0$}}}}
\end{picture}
}
\end{tabular}
}
\caption{Representation of lines $\lambda\mapsto
  \Ec(S)+\lambda\stdbars{S}$ for various subsets $S$. The
  $\ell_0$-curve, in plain line, is the minimal curve
  $\lambda\mapsto\min_S\stdacc{\Ec(S)+\lambda\stdbars{S}}$. It is
  continuous, concave, and piecewise affine with a finite number of
  pieces. The $\ell_0$-penalized regularization path is composed of
  the supports (here, $S_0^{\star}$, $S_1^{\star}$, $S_2^{\star}$)
  that are optimal for some $\lambda$-values. For instance,
  $S_1^{\star}$ is optimal for
  $\lambda\in[\lambda_2^\star,\lambda_1^\star]$. These supports
  $S^\star$ induce global minimizers of $\Jc(\xb;\lambda)$, defined as
  the least-square solutions $\xb_{S^{\star}}$. For instance,
  $\xb_{S_1^{\star}}$ is a global minimizer of $\Jc(\xb;\lambda)$ with
  respect to \xb whenever
  $\lambda\in[\lambda_2^\star,\lambda_1^\star]$.  }
  \label{fig:path}
\end{figure}
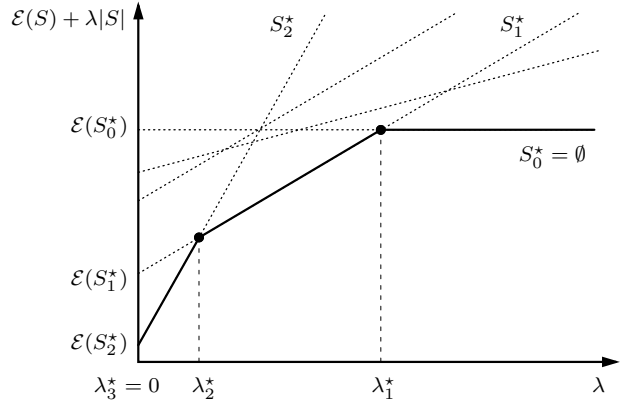

First, for any subset $S$, we can define a linear function
$\lambda\mapsto\Ec(S)+\lambda\stdbars{S}$, where
$\Ec(S)=\|\yb-\Ab\xb\|_2^2$ is the corresponding least-square error
and $\stdbars{S}$ stands for the cardinality of $S$. For each subset
$S$, this function yields a line in the 2D domain $(\lambda,\Jc)$, as
shown on Fig.~\ref{fig:path}.

Second, the set of solutions to~\eqref{eq:pen_LS} is piecewise
constant with respect to $\lambda$ (see Appendix~\ref{app:theory} for
a proof). Geometrically, this result can be easily understood by
noticing that the minimum of $\Jc(\xb;\lambda)$ with respect to \xb is
obtained for all $\lambda$-values by considering the concave envelope
of the set of lines $\lambda\mapsto\Ec(S)+\lambda\stdbars{S}$ for all
subsets $S$. The resulting piecewise affine curve is referred to as
the $\ell_0$-curve (see Fig.~\ref{fig:path}). Its edges are related to
the supports of the sparse solutions for all $\lambda$, and its
vertices yield the breakpoints $\lambda_i^\star$ around which the set
of optimal solutions $\argmin_\xb\Jc(\xb;\lambda)$ is changing.

We take advantage of this interpretation to propose two suboptimal
greedy algorithms that address~\eqref{eq:pen_LS} for a continuum of
$\lambda$-values. Continuation Single Best Replacement (CSBR)
repeatedly minimizes $\Jc(\xb;\lambda)$ with respect to \xb for
decreasing $\lambda$'s. $\ell_0$ Regularization Path Descent
($\ell_0$-PD) is a more complex algorithm maintaining \emph{a list of}
subsets so as to improve (decrease) the current approximation of the
$\ell_0$ curve.

\subsection{Related works}
\subsubsection{Bi-objective optimization}
The formulations~\eqref{eq:cstr_LS},~\eqref{eq:cstr_LS2}
and~\eqref{eq:pen_LS} can be interpreted as the same bi-objective
problem because they all intend to minimize both the approximation
error $\|\yb-\Ab\xb\|_2^2$ and the sparsity measure
$\|\xb\|_0$. Although \xb is continuous, the bi-objective optimization
problem should rather be considered as a discrete one where both
objectives reread $\Ec(S)$ and $\stdbars{S}$. Indeed, the continuous
solutions deduce from the discrete solutions, \xb reading as a
least-square minimizer among all vectors supported by $S$.

\begin{figure}[t]
\begin{center}
{
\setlength{\unitlength}{0.00040in}
\begingroup\makeatletter\ifx\SetFigFont\undefined%
\gdef\SetFigFont#1#2#3#4#5{%
  \reset@font\fontsize{#1}{#2pt}%
  \fontfamily{#3}\fontseries{#4}\fontshape{#5}%
  \selectfont}%
\fi\endgroup%
{\renewcommand{\dashlinestretch}{30}
\begin{picture}(5551,4700)(0,610)
\thicklines
\path(1095,4224)(1095,4044)
\path(1185,4134)(1005,4134)
\path(1095,4674)(1095,4494)
\path(1185,4584)(1005,4584)
\path(2220,4594)(2220,4414)
\path(2310,4504)(2130,4504)
\path(4470,1499)(4470,1319)
\path(4560,1409)(4380,1409)
\path(4470,2099)(4470,1919)
\path(4560,2009)(4380,2009)
\path(4470,2799)(4470,2619)
\path(4560,2709)(4380,2709)
\path(3345,2724)(3345,2544)
\path(3435,2634)(3255,2634)
\path(3345,3734)(3345,3554)
\path(3435,3644)(3255,3644)
\path(3345,3359)(3345,3179)
\path(3435,3269)(3255,3269)
\path(2220,3749)(2220,3569)
\path(2310,3659)(2130,3659)
\path(2220,4054)(2220,3874)
\path(2310,3964)(2130,3964)
\put(4470,1059){\blacken\ellipse{90}{90}}
\put(4470,1059){\ellipse{90}{90}}
\put(3345,1464){\blacken\ellipse{90}{90}}
\put(3345,1464){\ellipse{90}{90}}
\put(1095,3309){\blacken\ellipse{90}{90}}
\put(1095,3309){\ellipse{90}{90}}
\blacken\path(465.000,5009.000)(420.000,5159.000)(375.000,5009.000)(465.000,5009.000)
\path(420,5159)(420,884)(5670,884)
\blacken\path(5520.000,839.000)(5670.000,884.000)(5520.000,929.000)(5520.000,839.000)
\thinlines
\path(2220,3039)(2220,884)
\path(3345,1464)(3345,884)
\path(1095,3309)(1095,884)
\path(3345,1464)(420,1464)
\thicklines
\path(733,4900)(1095,3309)(2220,3039)
	(3345,1464)(4470,1065)(5345,990)
\thinlines
\path(2220,3039)(420,3039)
\dashline{60.000}(3345,1464)(1095,3309)
\path(1095,3309)(420,3309)
\thicklines
\whiten\path(2157.5,3111)(2282.5,3111)(2282.5,2967)
	(2157.5,2967)(2157.5,3111)
\path(2157.5,3111)(2282.5,3111)(2282.5,2967)
	(2157.5,2967)(2157.5,3111)
\put(1800,550){\makebox(0,0)[lb]{\smash{{\SetFigFont{9}{9}{\rmdefault}{\mddefault}{\updefault}$k+1$}}}}
\put(3030,550){\makebox(0,0)[lb]{\smash{{\SetFigFont{9}{9}{\rmdefault}{\mddefault}{\updefault}$k+2$}}}}
\put(3400,3850){\makebox(0,0)[lb]{\smash{{\SetFigFont{9}{9}{\rmdefault}{\mddefault}{\updefault}$S'$}}}}
\put(-600,3300){\makebox(0,0)[lb]{\smash{{\SetFigFont{9}{9}{\rmdefault}{\mddefault}{\updefault}$\Ec(S^{\star a}$)}}}}
\put(-600,1374){\makebox(0,0)[lb]{\smash{{\SetFigFont{9}{9}{\rmdefault}{\mddefault}{\updefault}$\Ec(S^{\star c})$}}}}
\put(5200,470){\makebox(0,0)[lb]{\smash{{\SetFigFont{9}{9}{\rmdefault}{\mddefault}{\updefault}$\stdbars{S}$}}}}
\put(-400,4879){\makebox(0,0)[lb]{\smash{{\SetFigFont{9}{9}{\rmdefault}{\mddefault}{\updefault}$\Ec(S)$}}}}
\put(1000,550){\makebox(0,0)[lb]{\smash{{\SetFigFont{9}{9}{\rmdefault}{\mddefault}{\updefault}$k$}}}}
\put(1150,3450){\makebox(0,0)[lb]{\smash{{\SetFigFont{9}{9}{\rmdefault}{\mddefault}{\updefault}$S^{\star a}$}}}}
\put(2350,3100){\makebox(0,0)[lb]{\smash{{\SetFigFont{9}{9}{\rmdefault}{\mddefault}{\updefault}$S^{\star b}$}}}}
\put(3300,1600){\makebox(0,0)[lb]{\smash{{\SetFigFont{9}{9}{\rmdefault}{\mddefault}{\updefault}$S^{\star c}$}}}}
\put(-600,2900){\makebox(0,0)[lb]{\smash{{\SetFigFont{9}{9}{\rmdefault}{\mddefault}{\updefault}$\Ec(S^{\star b})$}}}}
\end{picture}
}
}
\end{center}
\caption{Sparse approximation seen as a bi-objective optimization
  problem. The Pareto frontier gathers the non-dominated points: no
  other point can strictly decrease both $\stdbars{S}$ and $\Ec(S)$.
  Bullets and squares are all Pareto solutions. A supported solution
  is a minimizer of $\Ec(S)+\lambda\stdbars{S}$ with respect to $S$
  for some $\lambda$. $S^{\star a}$ and $S^{\star c}$ are supported,
  contrary to $S^{\star b}$. }
  \label{fig:pareto}
\end{figure}
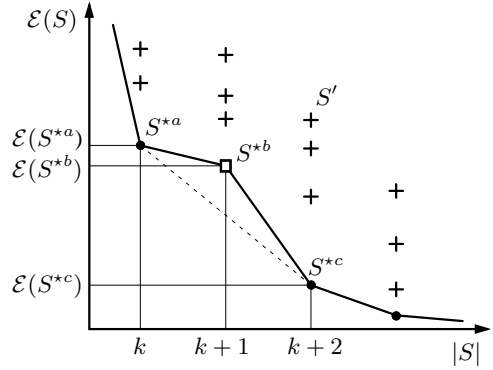

Fig.~\ref{fig:pareto} is a classical bi-objective representation where
each axis is related to a single objective~\cite{Das97}, namely
$\stdbars{S}$ and $\Ec(S)$. In bi-objective optimization, a point $S$
is called Pareto optimal when no other point $S'$ can decrease both
objectives~\cite{Marler04}. In the present context, $\stdbars{S}$
takes integer values, thus the Pareto solutions are the minimizers of
$\Ec(S)$ subject to $\stdbars{S}\leq k$ for consecutive values of
$k$. Equivalently, they minimize $\stdbars{S}$ subject to
$\Ec(S)\leq\varepsilon$ for some $\varepsilon$. They are usually
classified as supported or non-supported. The former lay on the convex
envelope of the Pareto frontier (the bullet points in
Fig.~\ref{fig:pareto}) whereas the latter lay in the nonconvex areas
(the square point). It is well known that a supported solution can be
reached when minimizing the weighted sum of both objectives, \ie when
minimizing $\Ec(S)+\lambda\stdbars{S}$ with respect to $S$ for some
weight $\lambda$. On the contrary, the non-supported solutions
cannot~\cite{Marler04}. Choosing between the weighting sum method and
a more complex method is a nontrivial question. The answer depends on
the problem at-hand and specifically, on the size of the nonconvex
areas in the Pareto frontier.

\subsubsection{$\ell_1$ and $\ell_0$-homotopy seen as a weighted sum method}
It is important to notice that for convex objectives, the Pareto
solutions are all supported. Consider the BPDN; because
$\|\yb-\Ab\xb\|_2^2$ and $\|\xb\|_1$ are convex functions of \xb, the
set of minimizers of $\|\yb-\Ab\xb\|_2^2+\lambda\|\xb\|_1$ for all
$\lambda$ coincides with the set of minimizers of $\|\yb-\Ab\xb\|_2^2$
s.t. $\|\xb\|_1\leq t$ for all $t$~\cite{vdBerg08}. Both sets are
referred to as the (unique) ``$\ell_1$-regularization path''. The
situation is different with $\ell_0$-regularization. Now, the weighted
sum formulation~\eqref{eq:pen_LS} may not yield the same solutions as
the constrained formulations~\eqref{eq:cstr_LS}
and~\eqref{eq:cstr_LS2} because the $\ell_0$-norm is
nonconvex~\cite{Nikolova13}. This will lead us to define two
$\ell_0$-regularization paths, namely the ``$\ell_0$-penalized path''
and the ``$\ell_0$-constrained path'' (Section~\ref{sec:notations}).

On the algorithmic side, the $\ell_0$ problems are acknowledged to be
difficult. Many authors actually discourage the direct optimization of
\Jc because there are a very large number of local
minimizers~\cite{Candes08,Wipf10}. In~\cite{Soussen11c}, however, we
showed that forward-backward extensions of OLS are able to escape from
some local minimizers of $\Jc(\xb;\lambda)$ for a given
$\lambda$. This motivates us to propose efficient OLS-based strategies
for minimizing \Jc for variable $\lambda$-values.

\subsubsection{Positioning with respect to other stepwise algorithms}
In statistical regression, the word ``stepwise'' originally refers to
Efroymson's algorithm~\cite{Efroymson60}, proposed in 1960 as an
empirical extension of forward selection (\ie OLS). Other stepwise
algorithms were proposed in the 1980's~\cite[Chapter~3]{Miller02}
among which Berk's and Broersen's
algorithms~\cite{Berk80,Broersen86}. All these algorithms perform a
single replacement per iteration, \ie a forward selection or a
backward elimination. They were originally applied to over-determined
problems in which the number of columns of \Ab is lower than the
number of rows. Recent stepwise algorithms were designed as either
OMP~\cite{Zhang11,Herzet14} or OLS
extensions~\cite{Haugland07,Chatterjee12}. They all aim to find
subsets of cardinality $k$ yielding a low approximation error $\Ec(S)$
for all $k$. Although our algorithms share the same objective, they
are inspired by \emph{(i)} the $\ell_1$-homotopy algorithm; and
\emph{(ii)} the structural properties of the $\ell_0$-regularization
paths. To the best of our knowledge, the idea of reconstructing an
$\ell_0$-regularization path using $\ell_0$-homotopy procedures is
novel.

CSBR and $\ell_0$-PD both read as descent algorithms in different
senses: CSBR, first sketched in~\cite{Duan09a}, repeatedly minimizes
$\Jc(\xb;\lambda)$ for decreasing $\lambda$'s. On the contrary,
$\ell_0$-PD minimizes $\Jc(\xb;\lambda)$ for any $\lambda$-value
\emph{simultaneously} by maintaining a list of candidate subsets. The
idea of maintaining a list of support candidates was recently
developed within the framework of forward
selection~\cite{Kwon14,Maymon15}. Our approach is different, because a
family of optimization problems are being addressed together. In
contrast, the supports in the list are all candidate solutions to
solve the same problem in~\cite{Kwon14,Maymon15}.

\subsubsection{Positioning with respect to continuation algorithms}
The principle of continuation is to handle a difficult problem by
solving a sequence of simpler problems with warm start initialization,
and gradually tuning some continuous
hyperparameter~\cite{Wasserstrom73}. In sparse approximation, the word
continuation is used in two opposite contexts.

First, the BDPN problem involving the $\ell_1$-norm. BPDN is solved
for decreasing hyperparameter values using the solution for each value
as a warm starting point for the next
value~\cite{Tropp10}. $\ell_1$-homotopy~\cite{Efron04,Malioutov05,Donoho08}
exploits that the $\ell_1$ regularization path is piecewise affine and
tracks the breakpoints between consecutive affine pieces. CSBR is
designed in a similar spirit and can be interpreted as an
``$\ell_0$-homotopy'' procedure (although the $\ell_0$ minimization
steps are solved in a sub-optimal way) working for decreasing
$\lambda$-values.

Second, the continuous approximation of the (discrete) $\ell_0$
pseudo-norm~\cite{Trzasko09} using a Graduated Non Convexity (GNC)
approach~\cite{Mohimani09}: a series of continuous concave metrics is
considered leading to the resolution of continuous optimization
problems with warm start initialization. Although the full
reconstruction of the $\ell_0$-regularization path has been rarely
addressed, it is noticeable that a GNC-like approach, called
SparseNet, aims to gradually update some estimation of the
regularization path induced by increasingly non-convex sparsity
measures~\cite{Mazumder11}. This strategy relies on the choice of a
grid of $\lambda$-values. Because the influence of the grid is
critical~\cite{vdBerg08}, it is useful to adapt the grid while the
nonconvex measure is modified~\cite{Mazumder11}. On the contrary, our
approach does not rely on a grid definition. The $\lambda$-values are
rather adaptively computed similar to the $\ell_1$-homotopy
principle~\cite{Efron04,Donoho08}.

The paper is organized as follows. In Section~\ref{sec:notations}, we
define the $\ell_0$-regularization paths and establish their main
properties. The CSBR and $\ell_0$-PD algorithms are respectively
proposed in Sections~\ref{sec:csbr} and~\ref{sec:track_l0}. In
Section~\ref{sec:simuls}, both algorithms are analyzed and compared
with the state-of-art algorithms based on nonconvex penalties for
difficult inverse problems. Additionally, we investigate the automatic
choice of the cardinality $k$ using classical order selection rules.

\section{$\ell_0$-regularization paths}
\label{sec:notations}

\subsection{Definitions, terminology and working assumptions}
Let $m\times n$ denote the size of the dictionary \Ab (usually, $m\leq
n$ in sparse approximation). The observation signal \yb and the weight
vector \xb are of size $m\times 1$ and $n\times 1$, respectively. We
assume that any $\min(m,n)$ columns of \Ab are linearly independent so
that for any subset $S\subset\stdacc{1,\ldots,n}$, the submatrix of
\Ab gathering the columns indexed by $S$ is full rank, and the
least-square error $\Ec(S)$ can be numerically computed. This
assumption is however not necessary for the theoretical results
provided hereafter.

We denote by $\stdbars{S}$ the cardinality of a subset $S$. We use the
alternative notations ``$S+\{i\}$'' and ``$S-\{i\}$'' for the forward
selection $S\cup\stdacc{i}$ and backward elimination
$S\setminus\stdacc{i}$. We can then introduce the generic notation
$S\pm\{i\}$ for single replacements: $S\pm\{i\}$ stands for $S+\{i\}$
if $i\notin S$, and $S-\{i\}$ if $i\in S$. We
will frequently resort to the geometrical interpretation of
Fig.~\ref{fig:path}. With a slight abuse of terminology, the line
$\lambda\mapsto\Ec(S)+\lambda\stdbars{S}$ will be simply referred to
as ``the line $S$''.

Hereafter, we start by defining the $\ell_0$-regularized paths as the
set of supports of the solutions to problems~\eqref{eq:cstr_LS},
\eqref{eq:cstr_LS2} and~\eqref{eq:pen_LS} for varying
hyperparameters. As seen in Section~\ref{sec:intro}, the solutions may
differ whether the $\ell_0$-regularization takes the form of a bound
constraint or a penalty. This will lead us to distinguish the
``$\ell_0$-constrained path'' and the ``$\ell_0$-penalized path''. We
will keep the generic terminology ``$\ell_0$-regularization paths''
for statements that apply to both. The solutions delivered by our
greedy algorithms will be referred to as the ``approximate
$\ell_0$-penalized path'' since they are suboptimal algorithms.

\subsection{Definition and properties of the $\ell_0$-regularized paths}
\label{sec:optimalpath}
The continuous problems~\eqref{eq:cstr_LS}, \eqref{eq:cstr_LS2}
and~\eqref{eq:pen_LS} can be converted as the discrete problems:
\begin{align}
&\min_S\,{\Ec(S)}\;\;\;\textrm{subject to}\;\;\;
\stdbars{S}\leq k, \label{eq:Sc}\\
&\min_S\,{\stdbars{S}}\;\;\;\;\;\;\textrm{subject to}\;\;\;
\Ec(S)\leq\varepsilon,\label{eq:Sc2}\\
&\min_{S}\,\bigacc{
\hat{\Jc}(S;\lambda)\triangleq\Ec(S)+\lambda\stdbars{S}},
\label{eq:Sp}
\end{align}
where $S$ stands for the support of \xb.
The optimal solutions \xb to problems~\eqref{eq:cstr_LS},
\eqref{eq:cstr_LS2} and~\eqref{eq:pen_LS} can indeed be simply deduced
from those of~\eqref{eq:Sc},~\eqref{eq:Sc2} and~\eqref{eq:Sp},
respectively, \xb reading as the least-square minimizers among all
vectors supported by $S$. In the following, the
formulation~\eqref{eq:Sc2} will be omitted because it leads to the
same $\ell_0$-regularization path as
formulation~\eqref{eq:Sc}~\cite{Nikolova13}.

Let us first define the set of solutions to~\eqref{eq:Sc}
and~\eqref{eq:Sp} and the $\ell_0$-curve, related to the minimum value
in~\eqref{eq:Sp} for all $\lambda>0$.
\begin{defi}
  \label{def:l0curve}
  For $k\leq\min(m,n)$, let $\Sc^\star_{\mathrm{C}}(k)$ be the set of
  minimizers of the constrained problem~\eqref{eq:Sc}.

  For $\lambda>0$, let $\Sc^\star_{\mathrm{P}}(\lambda)$ be the set of
  minimizers of the penalized problem~\eqref{eq:Sp}. Additionally, we
  define the $\ell_0$-curve as the function
  $\lambda\mapsto\min_S\{\hat{\Jc}(S;\lambda)\}$. It is the concave
  envelope of a finite number of linear functions. Thus, it is concave
  and piecewise affine. Let $\lambda^\star_{I+1}\triangleq
  0<\lambda^\star_I<\ldots<\lambda^\star_1
  <\lambda^\star_0\triangleq+\infty$ delimit the affine intervals
  ($I+1$ contiguous intervals; see Fig.~\ref{fig:path} in the case
  where $I=2$).
\end{defi}

Each set $\Sc^\star_{\mathrm{C}}(k)$ or
$\Sc^\star_{\mathrm{P}}(\lambda)$ can be thought of as a single
support (\eg $\Sc^\star_{\mathrm{C}}(k)$ is reduced to the support
$S^{\star a}$ in the example of Fig.~\ref{fig:pareto}). They are
defined as sets of supports because the minimizers of~\eqref{eq:Sc}
and~\eqref{eq:Sp} might not be always unique.
Let us now provide a key property of the set
$\Sc^\star_{\mathrm{P}}(\lambda)$.
\begin{theorem}
  \label{th:1}
  $\Sc^\star_{\mathrm{P}}(\lambda)$ is a piecewise constant function
  of $\lambda$, being constant on each interval
  $\lambda\in(\lambda_{i+1}^\star,\lambda_{i}^\star)$.
\end{theorem}
\begin{IEEEproof}
  See Appendix~\ref{app:theory}.
\end{IEEEproof}
This property allows us to define the 
$\ell_0$-regularization paths in a simple way.
\begin{defi}
  \label{def:Cpath}
  The $\ell_0$-constrained path is the set (of sets)
  $\Sc^\star_{\mathrm{C}}=\stdacc{\Sc^\star_{\mathrm{C}}(k),\,k=0,\ldots,\min(m,n)}$.

  The $\ell_0$-penalized path is defined as
  $\Sc^\star_{\mathrm{P}}=\stdacc{\Sc^\star_{\mathrm{P}}(\lambda),\,\lambda>0}$.
  According to Theorem~\ref{th:1}, $\Sc^\star_{\mathrm{P}}$ is
  composed of $(I+1)$ distinct sets $\Sc^\star_{\mathrm{P}}(\lambda)$,
  one for each interval
  $\lambda\in(\lambda_{i+1}^\star,\lambda_{i}^\star)$.
\end{defi}
$\Sc^\star_{\mathrm{C}}$ gathers the solutions to~\eqref{eq:Sc} for
all $k$. As illustrated on Fig.~\ref{fig:pareto}, the elements of
$\Sc^\star_{\mathrm{C}}$ are the Pareto solutions whereas the elements
of $\Sc^\star_{\mathrm{P}}$ correspond to the convex envelope of the
Pareto frontier.
Therefore, both $\ell_0$-regularization paths may not
coincide~\cite{Das97,Nikolova13}. As stated in
Theorem~\ref{th:2complet},
$\Sc^\star_{\mathrm{P}}\subset\Sc^\star_{\mathrm{C}}$, but the reverse
inclusion is not guaranteed.
\begin{theorem}
\label{th:2complet}
$\Sc^\star_{\mathrm{P}}\subset\Sc^\star_{\mathrm{C}}$. Moreover, for any
$\lambda\notin\stdacc{\lambda^\star_I,\ldots,\lambda^\star_0}$, there
exists $k$ such that $\Sc^\star_{\mathrm{P}}(\lambda)=\Sc^\star_{\mathrm{C}}(k)$.
\end{theorem}
\begin{IEEEproof}
  See Appendix~\ref{app:theory}.
\end{IEEEproof}

\subsection{Approximate $\ell_0$-penalized regularization path}
Let us introduce notations for the \emph{approximate}
$\ell_0$-penalized path delivered by our heuristic search
algorithms. Throughout the paper, the $\star$ notation is reserved for
optimal solutions (\eg $\Sc_{\mathrm{P}}^\star$). It is removed when
dealing with numerical solutions. The outputs of our algorithms will
be composed of a list
$\lambdab=\stdacc{\lambda_1,\ldots,\lambda_{J+1}}$ of decreasing
$\lambda$-values, and a list $\Sc=\stdacc{S_0,\ldots,S_J}$ of
candidate supports, with $S_0=\emptyset$. $S_j$ is a suboptimal
solution to~\eqref{eq:Sp} for
$\lambda\in(\lambda_{j+1},\lambda_j)$. In the first interval
$\lambda>\lambda_1$, the solution is $S_0=\emptyset$. The reader shall
keep in mind that each output $S_j$ induces a suboptimal solution
$\xb_j$ to~\eqref{eq:pen_LS} for
$\lambda\in(\lambda_{j+1},\lambda_j)$. This vector is the least-square
solution supported by $S_j$. It can be computed using the
pseudo-inverse of the subdictionary indexed by the set of atoms in
$S_j$.

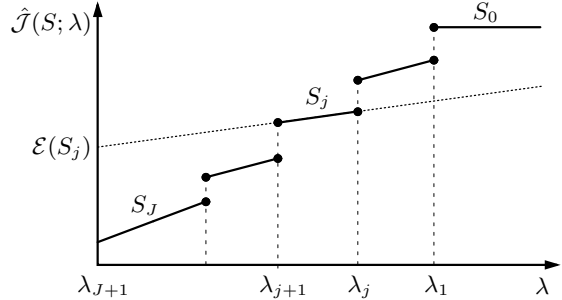
\begin{figure}[t]
\FINAL{\begin{center}}
\DRAFT{\centering}
\setlength{\unitlength}{0.00038in}
\begingroup\makeatletter\ifx\SetFigFont\undefined%
\gdef\SetFigFont#1#2#3#4#5{%
  \reset@font\fontsize{#1}{#2pt}%
  \fontfamily{#3}\fontseries{#4}\fontshape{#5}%
  \selectfont}%
\fi\endgroup%
{\renewcommand{\dashlinestretch}{30}
\begin{picture}(7485,3921)(0,250)
\thicklines
\put(5665,3739){\blacken\ellipse{90}{90}}
\put(5665,3739){\ellipse{90}{90}}
\put(5665,3289){\blacken\ellipse{90}{90}}
\put(5665,3289){\ellipse{90}{90}}
\put(4620,3011){\blacken\ellipse{90}{90}}
\put(4620,3011){\ellipse{90}{90}}
\put(4620,2582){\blacken\ellipse{90}{90}}
\put(4620,2582){\ellipse{90}{90}}
\put(3525,2432){\blacken\ellipse{90}{90}}
\put(3525,2432){\ellipse{90}{90}}
\put(3525,1939){\blacken\ellipse{90}{90}}
\put(3525,1939){\ellipse{90}{90}}
\put(2535,1679){\blacken\ellipse{90}{90}}
\put(2535,1679){\ellipse{90}{90}}
\put(2535,1344){\blacken\ellipse{90}{90}}
\put(2535,1344){\ellipse{90}{90}}
\thinlines
\dashline{60.000}(3525,2432)(3525,429)  
\thicklines
\path(5665,3739)(7125,3739)   
\thinlines
\dottedline{60.000}(1050,2093)(7170,2931)
\thicklines
\path(3525,2429)(4620,2582) 
\path(2535,1679)(3525,1939) 
\path(1815,1074)(2535,1344) 
\thinlines
\dashline{60.000}(5665,3739)(5665,474) 
\dashline{60.000}(2535,1679)(2535,474) 
\dashline{60.000}(4620,3011)(4620,474) 
\thicklines
\path(1050,787)(1815,1074)  
\path(4620,3011)(5665,3289) 
\blacken\path(1095.000,3894.000)(1050.000,4044.000)(1005.000,3894.000)(1095.000,3894.000)
\path(1050,4044)(1050,474)(7350,474)  
\blacken\path(7200.000,429.000)(7350.000,474.000)(7200.000,519.000)(7200.000,429.000)
\put(-150,3704){\makebox(0,0)[lb]{\smash{{\SetFigFont{9}{9}{\rmdefault}{\mddefault}{\updefault}$\hat{\Jc}(S;\lambda)$}}}}
\put(150,2000){\makebox(0,0)[lb]{\smash{{\SetFigFont{9}{9}{\rmdefault}{\mddefault}{\updefault}$\Ec(S_j)$}}}}
\put(6200,3850){\makebox(0,0)[lb]{\smash{{\SetFigFont{9}{9}{\rmdefault}{\mddefault}{\updefault}$S_0$}}}}
\put(3900,2689){\makebox(0,0)[lb]{\smash{{\SetFigFont{9}{9}{\rmdefault}{\mddefault}{\updefault}$S_j$}}}}
\put(1500,1200){\makebox(0,0)[lb]{\smash{{\SetFigFont{9}{9}{\rmdefault}{\mddefault}{\updefault}$S_{J}$}}}}
\put(7055,70){\makebox(0,0)[lb]{\smash{{\SetFigFont{9}{9}{\rmdefault}{\mddefault}{\updefault}$\lambda$}}}}
\put(5530,100){\makebox(0,0)[lb]{\smash{{\SetFigFont{9}{9}{\rmdefault}{\mddefault}{\updefault}$\lambda_1$}}}}
\put(4480,100){\makebox(0,0)[lb]{\smash{{\SetFigFont{9}{9}{\rmdefault}{\mddefault}{\updefault}$\lambda_j$}}}}
\put(3265,100){\makebox(0,0)[lb]{\smash{{\SetFigFont{9}{9}{\rmdefault}{\mddefault}{\updefault}$\lambda_{j+1}$}}}}
\put(770,100){\makebox(0,0)[lb]{\smash{{\SetFigFont{9}{9}{\rmdefault}{\mddefault}{\updefault}$\lambda_{J+1}$}}}}
\end{picture}
}
\FINAL{\end{center}}
\caption{ Notations relative to our heuristic search algorithms.
  Their outputs are: \emph{(i)} a sequence of values $\lambda_j$
  sorted in the decreasing order;
  \emph{(ii)} as many supports $S_j$, $S_j$ being the solution
  associated to all $\lambda\in(\lambda_{j+1},\lambda_{j})$. By
  extension, $S_0=\emptyset$ for $\lambda>\lambda_1$.  }
\label{fig:tracking}
\end{figure}

Geometrically, each support $S_j$ yields a line segment. Appending
these segments yields an approximate $\ell_0$-curve covering the
domain $(\lambda_{J+1},+\infty)$, as illustrated on
Fig.~\ref{fig:tracking}.

\section{Greedy continuation algorithm (CSBR)}
\label{sec:csbr}
Our starting point is the Single Best Replacement
algorithm~\cite{Soussen11c} dedicated to the minimization of
$\Jc(\xb;\lambda)$ with respect to \xb, or equivalently to
$\hat{\Jc}(S;\lambda)=\Ec(S)+\lambda\stdbars{S}$ with respect to
$S$. We first describe SBR for a given $\lambda$. Then, the CSBR
extension is presented for decreasing and adaptive $\lambda$'s.

\subsection{Single Best Replacement}
\label{sec:sbr}
SBR is a deterministic descent algorithm dedicated to the minimization
of $\hat{\Jc}(S;\lambda)$ with the initial solution $S=\emptyset$. An
SBR iteration consists of three steps:
\begin{enumerate}
\item Compute $\hat{\Jc}({S\pm\{i\}};\lambda)$ for all possible
  single replacements $S\pm\{i\}$ ($n$ insertion and removal trials);
\item Select the best replacement $S_{\mathrm{best}}=S\pm\{\ell\}$, with
  \begin{align}
  \ell&\in\argmin_{i\in\stdacc{1,\ldots,n}}\,\hat{\Jc}({S\pm\{i\}};\lambda);
  \label{eq:update_SBR}
  \end{align}
\item Update $S\leftarrow S_{\mathrm{best}}$.
\end{enumerate}
SBR terminates when
$\hat{\Jc}(S_{\mathrm{best}};\lambda)\geq\hat{\Jc}(S;\lambda)$, \ie
when no single replacement can decrease the cost function. This occurs
after a finite number of iterations because SBR is a descent algorithm
and there are a finite number of possible subsets
$S\subset\{1,\ldots,n\}$. In the limit case $\lambda=0$, we have
$\hat{\Jc}(S;0)=\Ec(S)$. Only insertions can be performed since any
removal increases the squared error $\Ec(S)$. SBR coincides with the
well-known OLS algorithm~\cite{Chen89}. Generally, the $n$ replacement
trials necessitate to compute $\Ec(S+\{i\})$ for all insertion trials
and $\Ec(S-\{i\})$ for all removals. In~\cite{Soussen11c}, we proposed
a fast and stable recursive implementation based on the Cholesky
factorization of the Gram matrix $\Ab_S^T\Ab_S$ when $S$ is modified
by one element (where $\Ab_S$ stands for the submatrix of \Ab
gathering the active columns).
SBR is summarized in Tab.~\ref{tab:SBR}. The optional output
parameters $\ell_\mathrm{add}$ and $\delta\Ec_\mathrm{add}$ are
unnecessary in the standard version. Their knowledge will be useful to
implement the extended CSBR algorithm.
\begin{table}[t]
  \caption{SBR algorithm for minimization of
    $\hat{\Jc}(S;\lambda)$ for fixed $\lambda$~\cite{Soussen11c}. By
    default, $S_{\mathrm{init}}=\emptyset$. The outputs 
    $\delta\Ec_\mathrm{add}$ and $\ell_\mathrm{add}$ are optional. 
    The single replacement tests appear in the for loop.
  }
\label{tab:SBR}
\centering
\begin{tabular}{l}
\hline
\textbf{inputs}~~: \Ab, \yb, $\lambda$, $S_{\mathrm{init}}$ \\
\textbf{outputs}: $S$, $\delta\Ec_\mathrm{add}$, $\ell_\mathrm{add}$ \\
\hline
$S_{\mathrm{best}}\leftarrow S_{\mathrm{init}}$;\\ 
\textbf{repeat}\\
\hspace*{.5cm}$S\leftarrow S_{\mathrm{best}}$;\\
\hspace*{.5cm}\textbf{for} $i=1$ \textbf{to} $n$ \textbf{do}\\
\hspace*{1.cm}Compute $\hat{\Jc}({S\pm\{i\}};\lambda)$;\\
\hspace*{.5cm}\textbf{end}\\
\hspace*{.5cm}$S_{\mathrm{best}}\leftarrow S\pm\{\ell\}$
with $\ell$ computed from~\eqref{eq:update_SBR};\\[.1cm]
\textbf{until} $\hat{\Jc}(S_{\mathrm{best}};\lambda)\geq\hat{\Jc}({S};\lambda)$;\\[.07cm]
Compute $\ell_\mathrm{add}$ according to~\eqref{eq:lplus};\\[.07cm]
Set $\delta\Ec_\mathrm{add}=\Ec(S)-\Ec(S+\{\ell_\mathrm{add}\})$;\\
\hline
\end{tabular}
\end{table}

\begin{figure}[t]
\centering
{
\setlength{\tabcolsep}{0.05cm}
\begin{tabular}{c}
\setlength{\unitlength}{0.00040in}
\begingroup\makeatletter\ifx\SetFigFont\undefined%
\gdef\SetFigFont#1#2#3#4#5{%
  \reset@font\fontsize{#1}{#2pt}%
  \fontfamily{#3}\fontseries{#4}\fontshape{#5}%
  \selectfont}%
\fi\endgroup%
{\renewcommand{\dashlinestretch}{30}
\begin{picture}(7350,5437)(150,1050)
\thicklines
\put(3530,5055){\blacken\ellipse{90}{90}}
\put(3530,5055){\ellipse{90}{90}}
\thinlines
\path(1050,5055)(6450,5055)  
\path(1050,2137)(6450,4605)  
\dashline{60.000}(3530,3270)(3530,1600)
\thicklines
\path(3530,5055)(3530,3270)
\blacken\path(1095.000,6150.000)(1050.000,6300.000)(1005.000,6150.000)(1095.000,6150.000)
\path(1050,6300)(1050,1600)(6900,1600)
\blacken\path(6750.000,1555)(6900.000,1600.000)(6750.000,1645)(6750.000,1555)
\thinlines
\dottedline{60.000}(1050,4370)(6450,5370) 
\dottedline{60.000}(1050,3220)(6450,5688) 
\dottedline{60.000}(1050,1600)(5505,5688) 
\blacken\path(3445.000,3420.000)(3530.000,3270.000)(3615.000,3420.000)(3445.000,3420.000)
\put(6670,1200){\makebox(0,0)[lb]{\smash{{\SetFigFont{9}{9}{\rmdefault}{\mddefault}{\updefault}$\lambda$}}}}
\put(-200,5975){\makebox(0,0)[lb]{\smash{{\SetFigFont{9}{9}{\rmdefault}{\mddefault}{\updefault}$\hat{\Jc}(S;\lambda)$}}}}
\put(3420,1200){\makebox(0,0)[lb]{\smash{{\SetFigFont{9}{9}{\rmdefault}{\mddefault}{\updefault}$\lambda$}}}}
\put(950,1200){\makebox(0,0)[lb]{\smash{{\SetFigFont{9}{9}{\rmdefault}{\mddefault}{\updefault}0}}}}
\put(6630,4980){\makebox(0,0)[lb]{\smash{{\SetFigFont{9}{9}{\rmdefault}{\mddefault}{\updefault}$S_{\mathrm{init}}=\emptyset$}}}}
\put(6630,4555){\makebox(0,0)[lb]{\smash{{\SetFigFont{9}{9}{\rmdefault}{\mddefault}{\updefault}$S=\{b,c\}$}}}}
\put(3260,4840){\makebox(0,0)[lb]{\smash{{\SetFigFont{9}{9}{\rmdefault}{\mddefault}{\updefault}+}}}}
\put(3260,4450){\makebox(0,0)[lb]{\smash{{\SetFigFont{9}{9}{\rmdefault}{\mddefault}{\updefault}+}}}}
\put(3260,3955){\makebox(0,0)[lb]{\smash{{\SetFigFont{9}{9}{\rmdefault}{\mddefault}{\updefault}+}}}}
\put(3245,3455){\makebox(0,0)[lb]{\smash{{\SetFigFont{9}{9}{\rmdefault}{\mddefault}{\updefault}--}}}}
\put(6630,5370){\makebox(0,0)[lb]{\smash{{\SetFigFont{9}{9}{\rmdefault}{\mddefault}{\updefault}$\{a\}$}}}}
\put(6630,5800){\makebox(0,0)[lb]{\smash{{\SetFigFont{9}{9}{\rmdefault}{\mddefault}{\updefault}$\{a,b\}$}}}}
\put(5005,5900){\makebox(0,0)[lb]{\smash{{\SetFigFont{9}{9}{\rmdefault}{\mddefault}{\updefault}$\{a,b,c\}$}}}}
\put(-250,3145){\makebox(0,0)[lb]{\smash{{\SetFigFont{9}{9}{\rmdefault}{\mddefault}{\updefault}$\Ec(\{a,b\})$}}}}
\put(-220,2130){\makebox(0,0)[lb]{\smash{{\SetFigFont{9}{9}{\rmdefault}{\mddefault}{\updefault}$\Ec(\{b,c\})$}}}}
\put(-540,1550){\makebox(0,0)[lb]{\smash{{\SetFigFont{9}{9}{\rmdefault}{\mddefault}{\updefault}$\Ec(\{a,b,c\})$}}}}
\put(30,4240){\makebox(0,0)[lb]{\smash{{\SetFigFont{9}{9}{\rmdefault}{\mddefault}{\updefault}$\Ec(\{a\})$}}}}
\put(250,5020){\makebox(0,0)[lb]{\smash{{\SetFigFont{9}{9}{\rmdefault}{\mddefault}{\updefault}$\|\yb\|_2^2$}}}}
\end{picture}
}
\end{tabular}
}
\caption{Step-by-step illustration of the call
  $S=\textrm{SBR}(\emptyset;\lambda)$. Each single replacement is
  represented by a vertical displacement (from top to bottom) from
  lines $S$ to $S\pm\{\ell\}$. The symbols `+' and `-' respectively
  refer to the selection and de-selection of atoms $a$, $b$ and $c$. 
  Four SBR iterations are done from the initial support
  $S_{\mathrm{init}}=\emptyset$: the selection of $a$, $b$
  and $c$, and the de-selection of $a$. The final output 
  $S\leftarrow\{b,c\}$ is of cardinality 2.  }
  \label{fig:sbr}
\end{figure}
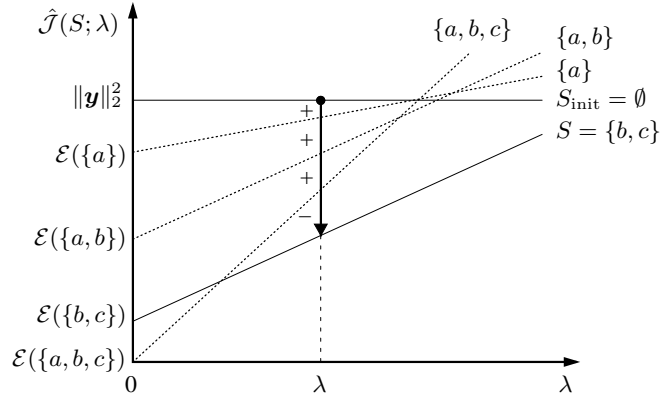

Let us illustrate the behavior of SBR on a simple example using the
geometrical interpretation of Fig.~\ref{fig:sbr}, where a single
replacement is represented by a vertical displacement (from top to
bottom) between the two lines $S$ and $S\pm\{\ell\}$.
$S_{\mathrm{init}}=\emptyset$ yields an horizontal line since
$\hat{\Jc}(\emptyset;\lambda)=\|\yb\|_2^2$ does not depend on
$\lambda$. At the first SBR iteration, a new dictionary atom $\ell=a$
is selected. The line related to the updated support
$S\leftarrow\{a\}$ is of slope $|S|=1$. Similarly, some new dictionary
atoms $b$ and $c$ are being selected in the next iterations, yielding
the supports $S\leftarrow\{a,b\}$ and $S\leftarrow\{a,b,c\}$.  On
Fig.~\ref{fig:sbr}, the dotted lines related to the latter supports
have slopes equal to 2 and 3. At iteration 4, the single best
replacement is the removal $\ell=a$. The resulting support
$S\leftarrow\{b,c\}$ is of cardinality 2, and the related line is
parallel to the line $\{a,b\}$ found at iteration 2. During the fifth
iteration, none of the $n$ single replacements decreases
$\hat{\Jc}(\{b,c\};\lambda)$. SBR stops with output $S=\{b,c\}$.

\subsection{Principle of the continuation search}
\label{sec:contin}
Our continuation strategy is inspired by $\ell_1$-homotopy which
recursively computes the minimizers of
$\|\yb-\Ab\xb\|_2^2+\lambda\|\xb\|_1$ when $\lambda$ is continuously
decreasing~\cite{Osborne00,Efron04,Donoho08}. An iteration of
$\ell_1$-homotopy consists in two steps:
\begin{itemize}
\item Find the next value
  $\lambda_{\mathrm{new}}<\lambda_{\mathrm{cur}}$ for which the
  $\ell_1$ optimality conditions are violated with the current active
  set $S$ ($\lambda_{\mathrm{cur}}$ denotes the current value);
\item Compute the single replacement
$S\leftarrow S\pm \{i\}$ allowing to fulfill the 
$\ell_1$ optimality conditions at $\lambda=\lambda_{\mathrm{new}}$.
\end{itemize}
CSBR follows the same principle. The first step is now related to some
\emph{local} $\ell_0$-optimality conditions, and the second step
consists in calling SBR at $\lambda_{\mathrm{new}}$ with the current
active set as initial solution; see Fig.~\ref{fig:csbr_path} for a
sketch. A main difference with $\ell_1$-homotopy is that the $\ell_0$
solutions are suboptimal, \ie they are local minimizers of
$\Jc(\xb;\lambda)$ with respect to \xb.

\subsubsection{Local optimality conditions}
Let us first reformulate the stopping conditions of SBR at a given
$\lambda$. SBR terminates when a local minimum of
$\hat{\Jc}({S};\lambda)$ has been found:
\begin{align}
\forall i\in\{1,\ldots,n\},\;\hat{\Jc}({S\pm\{i\}};\lambda)\geq
\hat{\Jc}({S};\lambda).
\label{eq:stop_SBR0}
\end{align}
This condition is illustrated on Fig.~\ref{fig:term_sbr}(a):
all lines related to single replacements $S\pm\{i\}$ lay above the
black point representing the value of $\hat{\Jc}(S;\lambda)$ for the
current $\lambda$.
By separating the conditions related to insertions $S+\{i\}$ and
removals $S-\{i\}$, \eqref{eq:stop_SBR0} rereads as the interval
condition:
\begin{align}
\lambda\in[\delta\Ec_\mathrm{add}(S),\delta\Ec_\mathrm{rmv}(S)],
\label{eq:stop_SBR}
\end{align}
where
\begin{subequations}
\label{eq:lambda_plusmoins}
\begin{align}
\delta\Ec_\mathrm{add}(S)&\triangleq\max_{i\notin S}\bigacc{\Ec(S)-\Ec(S+\{i\})}
\label{eq:lambda_plus}\\
\delta\Ec_\mathrm{rmv}(S)&\triangleq\min_{i\in S}\bigacc{\Ec(S-\{i\})-\Ec(S)}
\label{eq:lambda_moins}
\end{align}
\end{subequations}
refer to the maximum variation of the squared error when an atom is
added in the support $S$ (respectively, removed from $S$).
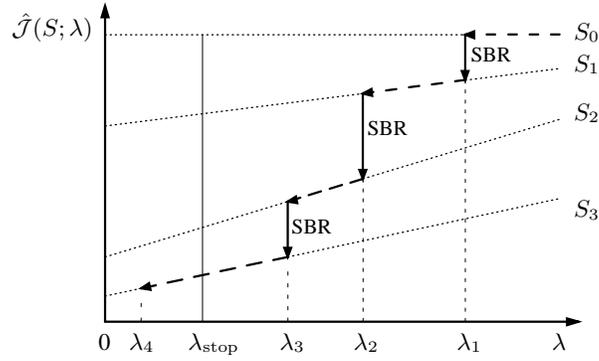
\begin{figure}
\begin{center}
\setlength{\unitlength}{0.00040in}
\begingroup\makeatletter\ifx\SetFigFont\undefined%
\gdef\SetFigFont#1#2#3#4#5{%
  \reset@font\fontsize{#1}{#2pt}%
  \fontfamily{#3}\fontseries{#4}\fontshape{#5}%
  \selectfonts3}%
\fi\endgroup%
{\renewcommand{\dashlinestretch}{30}
\begin{picture}(7930,4700)(100,1650)
\dottedline{70.000}(1205,5760)(5955,5760)  
\blacken\thicklines
\path(5910.000,5280)(5955.000,5160.000)(6000.000,5280.000)(5910.000,5280) 
\thinlines
\dashline{60.000}(5955,5160)(5955,1970) %

\dottedline{70.000}(1205,4554)(4605,4984) 
\dottedline{70.000}(5955,5160)(7200,5312) 
\thicklines
\path(5955,5160)(5955,5760) 
\blacken\path(4718.674,5043.187)(4605.000,4984)(4729.563,4953.848)(4718.674,5043.187)
\thicklines
\dashline{160.000}(4605,4984)(5955,5160)  
\thicklines
\blacken\path(6075.000,5805.000)(5955.000,5760.000)(6075.000,5715.000)(6075.000,5805.000)
\thicklines\dashline{140.000}(5955,5760)(7200,5760)  
\thinlines
\dashline{60.000}(1680,2211)(1680,1970) 
\dashline{60.000}(3615,2626)(3615,1970) 
\blacken\thicklines
\path(3570.000,2946.000)(3615.000,2826.000)(3660.000,2946.000)(3570.000,2946.000)
\path(3615,2826)(3615,3556)  
\thinlines
\dottedline{70.000}(1205,2826)(3615,3556) 
\dottedline{70.000}(4605,3856)(7200,4642) 
\blacken\thicklines
\path(3723.844,3623.660)(3615.000,3556.000)(3741.495,3535.408)(3723.844,3623.660)
\dashline{227.000}(3615,3556)(4605,3856)
\thinlines
\dashline{60.000}(4605,3856)(4605,1970) 

\dottedline{70.000}(1205,2309)(1680,2411)
\dottedline{70.000}(3615,2826)(7200,3596)

\blacken\thicklines

\path(1785.195,2481.244)(1680.000,2411)(1804.951,2392.439)(1785.195,2481.244)

\dashline{227.000}(1680,2411)(3615,2826) 
\thinlines
\path(2490,5760)(2490,1970) 
\blacken\thicklines
\path(1250.000,5990.000)(1205.000,6140.000)(1160.000,5990.000)(1250.000,5990.000)
\path(1205,6140)(1205,1970)(7450,1970)   
\blacken\path(7300.000,1925.000)(7450.000,1970.000)(7300.000,2015.000)(7300.000,1925.000)
\blacken\path(4560.000,3976.000)(4605.000,3856.000)(4650.000,3976.000)(4560.000,3976) 
\thicklines
\path(4605,3856)(4605,4984) 
\put(7390,3350){\makebox(0,0)[lb]{\smash{{\SetFigFont{9}{9}{\rmdefault}{\mddefault}{\updefault}$S_3$}}}}
\put(7390,4630){\makebox(0,0)[lb]{\smash{{\SetFigFont{9}{9}{\rmdefault}{\mddefault}{\updefault}$S_2$}}}}
\put(-30,5750){\makebox(0,0)[lb]{\smash{{\SetFigFont{9}{9}{\rmdefault}{\mddefault}{\updefault}$\hat{\Jc}(S;\lambda)$}}}}
\put(7100,1600){\makebox(0,0)[lb]{\smash{{\SetFigFont{9}{9}{\rmdefault}{\mddefault}{\updefault}$\lambda$}}}}
\put(7390,5280){\makebox(0,0)[lb]{\smash{{\SetFigFont{9}{9}{\rmdefault}{\mddefault}{\updefault}$S_1$}}}}
\put(7390,5700){\makebox(0,0)[lb]{\smash{{\SetFigFont{9}{9}{\rmdefault}{\mddefault}{\updefault}$S_0$}}}}
\put(6030,5400){\makebox(0,0)[lb]{\smash{{\SetFigFont{8}{8}{\rmdefault}{\mddefault}{\updefault}SBR}}}}
\put(4665,4430){\makebox(0,0)[lb]{\smash{{\SetFigFont{8}{8}{\rmdefault}{\mddefault}{\updefault}SBR}}}}
\put(3665,3130){\makebox(0,0)[lb]{\smash{{\SetFigFont{8}{8}{\rmdefault}{\mddefault}{\updefault}SBR}}}}
\put(5840,1600){\makebox(0,0)[lb]{\smash{{\SetFigFont{9}{9}{\rmdefault}{\mddefault}{\updefault}$\lambda_1$}}}}
\put(4480,1600){\makebox(0,0)[lb]{\smash{{\SetFigFont{9}{9}{\rmdefault}{\mddefault}{\updefault}$\lambda_2$}}}}
\put(3505,1600){\makebox(0,0)[lb]{\smash{{\SetFigFont{9}{9}{\rmdefault}{\mddefault}{\updefault}$\lambda_3$}}}}
\put(2300,1600){\makebox(0,0)[lb]{\smash{{\SetFigFont{9}{9}{\rmdefault}{\mddefault}{\updefault}$\lambda_\mathrm{stop}$}}}}
\put(1520,1600){\makebox(0,0)[lb]{\smash{{\SetFigFont{9}{9}{\rmdefault}{\mddefault}{\updefault}$\lambda_4$}}}}
\put(1115,1600){\makebox(0,0)[lb]{\smash{{\SetFigFont{9}{9}{\rmdefault}{\mddefault}{\updefault}0}}}}
\end{picture}
}
\end{center}
\caption{Step-by-step illustration of CSBR with the early stopping
  condition $\lambda_j\leq\lambda_{\mathrm{stop}}$. The initial
  support is $S_0=\emptyset$. SBR is called for three decreasing
  values (plain vertical arrows), with output $S_j$ at
  $\lambda_j$. The search for the next value $\lambda_{j+1}$ is
  represented by an oblique displacement along the line $S_j$.  }
\label{fig:csbr_path}
\end{figure}
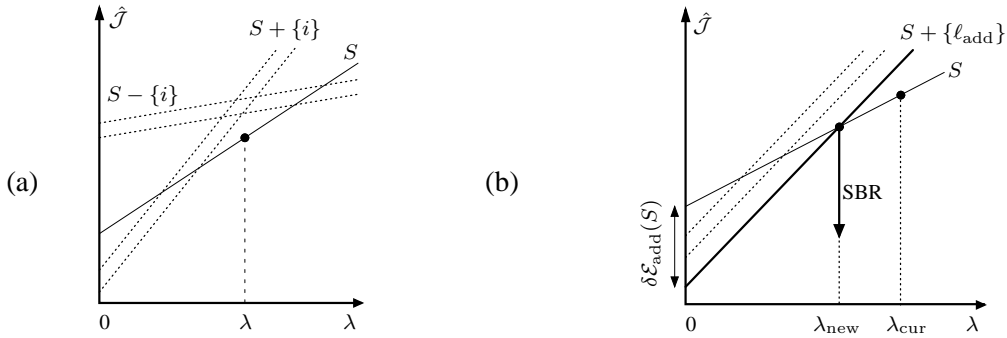
\begin{figure}
\centering
{
\FINAL{\setlength{\tabcolsep}{0.01cm}}
\DRAFT{\setlength{\tabcolsep}{0.3cm}}
\begin{tabular}{cc}
\DRAFT{\begin{tabular}{c}(a)\hspace*{0.1cm}\end{tabular}\begin{tabular}{c}}
\setlength{\unitlength}{0.00038in}
\begingroup\makeatletter\ifx\SetFigFont\undefined%
\gdef\SetFigFont#1#2#3#4#5{%
  \reset@font\fontsize{#1}{#2pt}%
  \fontfamily{#3}\fontseries{#4}\fontshape{#5}%
  \selectfont}%
\fi\endgroup%
{\renewcommand{\dashlinestretch}{30}
\begin{picture}(3750,4575)(1220,1760)
\thicklines
\put(3050,4521){\blacken\ellipse{90}{90}}
\put(3050,4521){\ellipse{90}{90}}
\blacken\path(1095.000,6150.000)(1050.000,6300.000)(1005.000,6150.000)(1095.000,6150.000)

\path(1050,6300)(1050,2250)(4647,2250) 

\blacken\path(4548.000,2205.000)(4647,2250)(4548.000,2295.000)(4548.000,2205.000)
\thinlines
\dottedline{70.000}(1050,4520)(4614,5120) 
\dottedline{70.000}(1050,4720)(4614,5320) 
\dottedline{70.000}(1050,2700)(3498,5750) 
\dottedline{70.000}(1050,2400)(3738,5750) 

\dashline{60.000}(3050,4521)(3050,2250)  
\thicklines
\thinlines
\path(1050,3202)(4599,5543) 
\blacken\thicklines
\thinlines
\put(4400,1900){\makebox(0,0)[lb]{\smash{{\SetFigFont{8}{8}{\rmdefault}{\mddefault}{\updefault}$\lambda$}}}}
\put(1180,6050){\makebox(0,0)[lb]{\smash{{\SetFigFont{8}{8}{\rmdefault}{\mddefault}{\updefault}$\hat{\Jc}$}}}}  
\put(2970,1900){\makebox(0,0)[lb]{\smash{{\SetFigFont{8}{8}{\rmdefault}{\mddefault}{\updefault}$\lambda$}}}}
\put(3100,5930){\makebox(0,0)[lb]{\smash{{\SetFigFont{8}{8}{\rmdefault}{\mddefault}{\updefault}$S+\stdacc{i}$}}}}
\put(4400,5600){\makebox(0,0)[lb]{\smash{{\SetFigFont{8}{8}{\rmdefault}{\mddefault}{\updefault}$S$}}}}
\put(1150,5030){\makebox(0,0)[lb]{\smash{{\SetFigFont{8}{8}{\rmdefault}{\mddefault}{\updefault}$S-\stdacc{i}$}}}}
\put(1050,1900){\makebox(0,0)[lb]{\smash{{\SetFigFont{8}{8}{\rmdefault}{\mddefault}{\updefault}0}}}}
\end{picture}
}
\DRAFT{\end{tabular}}&
\DRAFT{\begin{tabular}{c}(b)\hspace*{1.3cm}\end{tabular}\begin{tabular}{c}}
\setlength{\unitlength}{0.00038in}
\begingroup\makeatletter\ifx\SetFigFont\undefined%
\gdef\SetFigFont#1#2#3#4#5{%
  \reset@font\fontsize{#1}{#2pt}%
  \fontfamily{#3}\fontseries{#4}\fontshape{#5}%
  \selectfont}%
\fi\endgroup%
{\renewcommand{\dashlinestretch}{30}
\begin{picture}(4400,4625)(1005,1862)
\thicklines
\put(4000,5232){\blacken\ellipse{90}{90}}
\put(4000,5232){\ellipse{90}{90}}
\put(3157,4795){\blacken\ellipse{90}{90}}
\put(3157,4795){\ellipse{90}{90}}  
\path(1050,6400)(1050,2360)(5147,2360)
\blacken\path(1095,6250.000)(1050.000,6400.000)(1005,6250.000)(1095,6250.000)  
\blacken\path(5048,2305.000)(5147.000,2350.000)(5048,2395.000)(5048,2305.000) 

\path(1050,2600)(4170,5850) 
\thinlines
\path(1050,3702)(4599,5543) 
\blacken\thicklines
\path(3217,3446)(3157,3300)(3098,3446)(3217,3446)  
\path(3157,3300)(3157,4795)  
\thinlines
\blacken\path(950,3566)(900,3702)(850,3566)(950,3566)
\blacken\path(950,2736)(900,2600)(850,2736)(950,2736)
\path(900,3702)(900,2600)
\dottedline{70.000}(1050,3300)(3498,5850) 
\dottedline{70.000}(1050,3000)(3786,5850) 
\dottedline{70.000}(4000,5232)(4000,2350) 
\dottedline{70.000}(3157,3300)(3157,2350) 

\put(4895,2000){\makebox(0,0)[lb]{\smash{{\SetFigFont{8}{8}{\rmdefault}{\mddefault}{\updefault}$\lambda$}}}}
\put(1150,6100){\makebox(0,0)[lb]{\smash{{\SetFigFont{8}{8}{\rmdefault}{\mddefault}{\updefault}$\hat{\Jc}$}}}}  
\put(2798,2000){\makebox(0,0)[lb]{\smash{{\SetFigFont{8}{8}{\rmdefault}{\mddefault}{\updefault}$\lambda_{\mathrm{new}}$}}}}
\put(3202,3800){\makebox(0,0)[lb]{\smash{{\SetFigFont{8}{8}{\rmdefault}{\mddefault}{\updefault}SBR}}}}
\put(3800,2000){\makebox(0,0)[lb]{\smash{{\SetFigFont{8}{8}{\rmdefault}{\mddefault}{\updefault}$\lambda_{\mathrm{cur}}$}}}}
\put(3990,6000){\makebox(0,0)[lb]{\smash{{\SetFigFont{8}{8}{\rmdefault}{\mddefault}{\updefault}$S+\{\ell_{\mathrm{add}}\}$}}}}
\put(4650,5480){\makebox(0,0)[lb]{\smash{{\SetFigFont{8}{8}{\rmdefault}{\mddefault}{\updefault}$S$}}}}
\put(430,2600){\makebox(0,0)[lb]{{\SetFigFont{8}{8}{\rmdefault}{\mddefault}{\updefault}{\rotatebox{90}{$\delta\Ec_{\mathrm{add}}(S)$}}}}} %
\put(1050,2000){\makebox(0,0)[lb]{\smash{{\SetFigFont{8}{8}{\rmdefault}{\mddefault}{\updefault}0}}}}
\end{picture}
}
\DRAFT{\end{tabular}}
\FINAL{\\\small{(a)}&\small{(b)}}
\end{tabular}
}
\caption{Termination of SBR and next call to SBR.\; (a) When SBR
  terminates, no single replacement $S\pm\{i\}$ can decrease
  $\hat{\Jc}(S;\lambda)$. The dotted lines $S+\{i\}$ (of slope
  $\stdbars{S}+1$) lay above the black point
  $(\lambda,\hat{\Jc}(S;\lambda))$. Similarly, all lines $S-\{i\}$, of
  slope $\stdbars{S}-1$, lay above this point.\; (b) Here, $S$ is the
  SBR output at $\lambda_{\mathrm{cur}}$. The next call to SBR is done
  at $\lambda_{\mathrm{new}}=\delta\Ec_\mathrm{add}(S)$ with the
  initial subset $S+\{\ell_\mathrm{add}\}$. The line
  $S+\{\ell_\mathrm{add}\}$ lays below all other lines $S+\{i\}$
  (dotted lines). Here, the $\lambda$-axis has been stretched by an
  arbitrary factor for improved readability. The horizontal length
  $\lambda_{\mathrm{new}}$ does not match the vertical length
  $\delta\Ec_\mathrm{add}(S)$, as it should without any stretching.
  The same stretching process will be done in
  Fig.~\ref{fig:refinement}.  }
  \label{fig:term_sbr}
\end{figure}

\subsubsection{Violation of the local optimality conditions}
Consider the current output
$S=\mathrm{SBR}(S_{\mathrm{init}};\lambda_{\mathrm{cur}})$. The local
optimality condition~\eqref{eq:stop_SBR} is then met for
$\lambda=\lambda_{\mathrm{cur}}$, but also for any
$\lambda\in[\delta\Ec_\mathrm{add}(S), \lambda_{\mathrm{cur}}]$.  The
new value for which~\eqref{eq:stop_SBR} is violated is
$\lambda_{\mathrm{new}}=\delta\Ec_\mathrm{add}(S)-c$ where $c>0$ is
arbitrarily small. The violation occurs for $i=\ell_\mathrm{add}$,
with
\begin{align}
  \ell_\mathrm{add}\in\argmax_{i\notin S}\stdacc{\Ec(S)-\Ec(S+\{i\})}.
  \label{eq:lplus}
\end{align}
In practice, $\lambda_{\mathrm{new}}$ can be set to the limit value
\begin{align}
\lambda_{\mathrm{new}}=\delta\Ec_\mathrm{add}(S)
\label{eq:lambdanew}
\end{align}
provided that $S$ is replaced with $S+\{\ell_\mathrm{add}\}$.
\begin{table}[t]
  \caption{CSBR algorithm: SBR is called repeatedly for decreasing
    $\lambda_j$'s. At
    iteration $j$, both the next value $\lambda_{j+1}$ and the next
    initial subset $S_j+\{\ell_\mathrm{add}\}$ are provided as SBR outputs. 
}
\label{tab:CSBR}
\centering 
\begin{tabular}{l} 
\hline
\textbf{inputs}~~: \Ab, \yb \\
\textbf{outputs}: \Sc: list of supports $S_j$; \lambdab: list of $\lambda_j$ \\
\hline
$S_0\leftarrow\emptyset$;\\
$S_{\mathrm{init}}\leftarrow\{\ell_\mathrm{add}\}$
with $\ell_\mathrm{add}$ computed from~\eqref{eq:lplus_init};\\
Compute $\lambda_1$ according to~\eqref{eq:lplus_init};\\
$j\leftarrow1$;\\
\textbf{while} $\lambda_{j}>0$ \textbf{do}\\[.04cm]
\hspace*{.5cm}Call $\stdcro{S_j,\,\delta\Ec_\mathrm{add},\,\ell_\mathrm{add}}=
\textrm{SBR}(S_{\mathrm{init}};\lambda_{j})$;\\[.04cm]
\hspace*{.5cm}$S_{\mathrm{init}}\leftarrow S_j+\{\ell_\mathrm{add}\}$;\\[.04cm]
\hspace*{.5cm}$\lambda_{j+1}\leftarrow\delta\Ec_\mathrm{add}$;\\[.04cm]
\hspace*{.5cm}$j\leftarrow j+1$;\\
\textbf{end}\\
\hline
\end{tabular}
\end{table}

As illustrated on Fig.~\ref{fig:term_sbr}(b), the line
$S+\{\ell_\mathrm{add}\}$ lays below all other parallel lines
$S+\{i\}$. It intersects line $S$ at $\lambda_{\mathrm{new}}$. The
vertical arrow represents the new call to SBR with inputs
$S+\{\ell_\mathrm{add}\}$ and $\lambda_{\mathrm{new}}$. Because $S$
and $S+\{\ell_\mathrm{add}\}$ both lead to the same value of
$\hat{\Jc}(\,.\,;\lambda_{\mathrm{new}})$, the de-selection of
$\ell_\mathrm{add}$ is forbidden in the first iteration of SBR.

\subsection{CSBR algorithm}
\label{sec:csbr_impl}
CSBR is summarized in Tab.~\ref{tab:CSBR}. The repeated calls to SBR
deliver subsets $S_j$ for decreasing $\lambda_j$. As shown on
Fig.~\ref{fig:csbr_path}, the solution $S_j$ covers the interval
$(\lambda_{j+1},\lambda_j]$. At the very first iteration, we have
$S_0=\emptyset$, and \eqref{eq:lplus}-\eqref{eq:lambdanew} reread:
\begin{align}
  \ell_\mathrm{add}&\in\argmax_{i\in\{1,\ldots,n\}}\,
  \frac{\stdbars{\stdscal{\yb,\ab_i}}}{\|\ab_i\|_2}~~\mathrm{and}~~%
  \lambda_1=\frac{\stdscal{\yb,\ab_{\ell_\mathrm{add}}}^2}{\|\ab_{\ell_\mathrm{add}}\|_2^2}.
  \label{eq:lplus_init}
\end{align}
According to Tab.~\ref{tab:CSBR}, CSBR stops when $\lambda_j=0$, \ie
the whole domain $\lambda\in\Rbb_+$ has been scanned. However, this
choice may not be appropriate when dealing with noisy data and
overcomplete dictionaries. In such cases, \emph{ad hoc} early stopping
rules can be considered~\cite{Donoho07,Donoho08}. A natural rule takes
the form $\lambda_j\leq\lambda_{\mathrm{stop}}$ with
$\lambda_{\mathrm{stop}}>0$. Alternative rules involve a maximum
cardinality ($\stdbars{S_j}\geq k_{\mathrm{stop}}$) and/or a minimum
squared error ($\Ec(S_j)\leq\varepsilon_{\mathrm{stop}}$).

Fig.~\ref{fig:csbr_path} shows a step-by-step illustration with the
early stop $\lambda_j\leq\lambda_{\mathrm{stop}}$. The initial support
$S_{\mathrm{init}}=\{\ell_\mathrm{add}\}$ and $\lambda_1$ are
precomputed in~\eqref{eq:lplus_init}. In the first call
$S_1=\mathrm{SBR}(S_{\mathrm{init}};\lambda_1)$, a number of single
replacements updates $S\leftarrow S\pm\{\ell\}$ are carried out
leading to $S_1=S$. This process is represented by the plain vertical
arrow at $\lambda_1$ linking both lines $S_0$ and $S_1$ (the line
$S_{\mathrm{init}}$ is not shown for readability reasons). Once $S_1$
is obtained, the next value $\lambda_2$ is computed. This process is
represented by an oblique, dashed arrow joining $\lambda_1$ and
$\lambda_2$. These two processes are being repeated alternatively at
the second and third iterations of CSBR. Finally, CSBR terminates
after $\lambda_4$ has been computed because
$\lambda_4\leq\lambda_{\mathrm{stop}}$.

\section{$\ell_0$-regularization path descent ($\ell_0$-PD)}
\label{sec:track_l0}
On the theoretical side, the $\ell_0$-penalized regularization path is
piecewise constant (Theorem~\ref{th:1}). It yields the $\ell_0$ curve
which is piecewise affine, continuous and concave
(Fig.~\ref{fig:path}). The curve related to the CSBR outputs does not
fulfill this property since: \emph{(i)} there might be jumps in this
curve; and \emph{(ii)} the slope of the line $S_j$ is not necessarily
increasing with $j$ (see Fig.~\ref{fig:csbr_path}). This motivates us
to propose another algorithm whose outputs are consistent with the
structural properties of the $\ell_0$-curve.

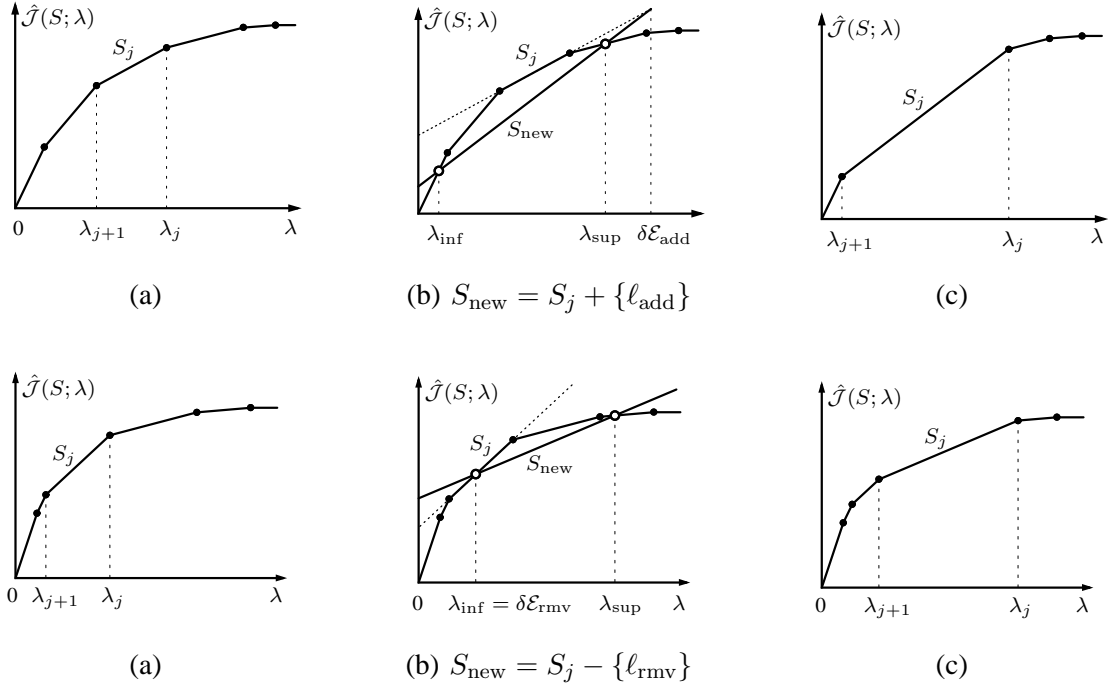
\begin{figure}[t]
\FINAL{
  \centering
  \setlength{\tabcolsep}{0.cm}
  \begin{tabular}{ccc}
  & 
  \begin{tabular}{c}\small{\fbox{\emph{$S_{\mathrm{new}}=S_j+\{\ell_\mathrm{add}\}$}}}\end{tabular} & 
  \begin{tabular}{c}\small{\fbox{\emph{$S_{\mathrm{new}}=S_j-\{\ell_\mathrm{rmv}\}$}}
  }
  \end{tabular}\\[.3cm]
  \begin{tabular}{c}\small{(a)}\end{tabular}&
  \begin{tabular}{c}
\setlength{\unitlength}{0.00029in}
\begingroup\makeatletter\ifx\SetFigFont\undefined%
\gdef\SetFigFont#1#2#3#4#5{%
  \reset@font\fontsize{#1}{#2pt}%
  \fontfamily{#3}\fontseries{#4}\fontshape{#5}%
  \selectfont}%
\fi\endgroup%
{\renewcommand{\dashlinestretch}{30}
\begin{picture}(5100,4681)(750,0)
\thicklines
\put(5550,3984){\blacken\ellipse{90}{90}}
\put(5550,3984){\ellipse{90}{90}}
\put(4965,3939){\blacken\ellipse{90}{90}}
\put(4965,3939){\ellipse{90}{90}}
\put(1375,1779){\blacken\ellipse{90}{90}}
\put(1375,1779){\ellipse{90}{90}}
\put(2319,2891){\blacken\ellipse{90}{90}}
\put(2319,2891){\ellipse{90}{90}}
\put(3579,3576){\blacken\ellipse{90}{90}}
\put(3579,3576){\ellipse{90}{90}}
\thinlines
\dashline{60.000}(3579,3576)(3579,675)
\blacken\thicklines
\path(895.000,4194.000)(850,4344.000)(805.000,4194.000)(895.000,4194.000)
\path(850,4344)(850,675)(5950,675)  
\blacken\path(5800.000,630.000)(5950.000,675.000)(5800.000,720.000)(5800.000,630.000)
\thinlines
\dashline{60.000}(2319,2891)(2319,675)
\thicklines
\path(850,675)(1375,1779)
	(2319,2891)(3579,3576)(4965,3939)
	(5550,3984)(5905,3984)
\thinlines
\thicklines
\put(810,200){\makebox(0,0)[lb]{\smash{{\SetFigFont{8}{8}{\rmdefault}{\mddefault}{\updefault}0}}}}
\put(2000,200){\makebox(0,0)[lb]{\smash{{\SetFigFont{8}{8}{\rmdefault}{\mddefault}{\updefault}$\lambda_{j+1}$}}}}
\put(3450,200){\makebox(0,0)[lb]{\smash{{\SetFigFont{8}{8}{\rmdefault}{\mddefault}{\updefault}$\lambda_j$}}}}
\put(950,4000){\makebox(0,0)[lb]{\smash{{\SetFigFont{8}{8}{\rmdefault}{\mddefault}{\updefault}$\hat{\Jc}(S;\lambda)$}}}}
\put(2600,3459){\makebox(0,0)[lb]{\smash{{\SetFigFont{8}{8}{\rmdefault}{\mddefault}{\updefault}$S_j$}}}}
\put(5700,200){\makebox(0,0)[lb]{\smash{{\SetFigFont{8}{8}{\rmdefault}{\mddefault}{\updefault}$\lambda$}}}}
\end{picture}
}
\end{tabular}&
  \;\begin{tabular}{c}
\setlength{\unitlength}{0.00023in}
\begingroup\makeatletter\ifx\SetFigFont\undefined%
\gdef\SetFigFont#1#2#3#4#5{%
  \reset@font\fontsize{#1}{#2pt}%
  \fontfamily{#3}\fontseries{#4}\fontshape{#5}%
  \selectfont}%
\fi\endgroup%
{\renewcommand{\dashlinestretch}{30}
\begin{picture}(6800,5943)(840,970)
\thicklines
\put(6518,5721){\blacken\ellipse{113}{113}}
\put(6518,5721){\ellipse{113}{113}}
\put(1635,3309){\blacken\ellipse{113}{113}}
\put(1635,3309){\ellipse{113}{113}}
\put(1838,3731){\blacken\ellipse{113}{113}}
\put(1838,3731){\ellipse{113}{113}}
\put(5288,5613){\blacken\ellipse{113}{113}}
\put(5288,5613){\ellipse{113}{113}}
\put(3300,5091){\blacken\ellipse{113}{113}}
\put(3300,5091){\ellipse{113}{113}}
\blacken\path(1185.000,6309.000)(1140.000,6459.000)(1095.000,6309.000)(1185.000,6309.000)
\path(1140,6459)(1140,1824)(7265,1824)  
\blacken\path(7115.000,1779.000)(7265.000,1824.000)(7115.000,1869.000)(7115.000,1779.000)
\thinlines
\dashline{90.000}(1838,3731)(1838,1824)  
\dashline{90.000}(3300,5091)(3300,1824)  
\thicklines
\path(1140,1824)(1635,3309)(1838,3731)
(3300,5091)(5288,5613)
(6518,5721)(7128,5721)
\put(950,1269){\makebox(0,0)[lb]{\smash{{\SetFigFont{8}{8.0}{\rmdefault}{\mddefault}{\updefault}0}}}}
\put(6950,1269){\makebox(0,0)[lb]{\smash{{\SetFigFont{8}{8.0}{\rmdefault}{\mddefault}{\updefault}$\lambda$}}}}
\put(1510,1269){\makebox(0,0)[lb]{\smash{{\SetFigFont{8}{8.0}{\rmdefault}{\mddefault}{\updefault}$\lambda_{j+1}$}}}}
\put(3100,1269){\makebox(0,0)[lb]{\smash{{\SetFigFont{8}{8.0}{\rmdefault}{\mddefault}{\updefault}$\lambda_j$}}}}
\put(1980,4524){\makebox(0,0)[lb]{\smash{{\SetFigFont{8}{8.0}{\rmdefault}{\mddefault}{\updefault}$S_j$}}}}
\put(1280,6020){\makebox(0,0)[lb]{\smash{{\SetFigFont{8}{8.0}{\rmdefault}{\mddefault}{\updefault}$\hat{\Jc}(S;\lambda)$}}}}
\end{picture}
}
\end{tabular}\\
  \begin{tabular}{c}\small{(b)}\end{tabular}&
  \begin{tabular}{c}
\setlength{\unitlength}{0.00029in}
\begingroup\makeatletter\ifx\SetFigFont\undefined%
\gdef\SetFigFont#1#2#3#4#5{%
  \reset@font\fontsize{#1}{#2pt}%
  \fontfamily{#3}\fontseries{#4}\fontshape{#5}%
  \selectfont}%
\fi\endgroup%
{\renewcommand{\dashlinestretch}{30}
\begin{picture}(5100,4881)(750,0)
\thicklines
\put(5550,3984){\blacken\ellipse{90}{90}}
\put(5550,3984){\ellipse{90}{90}}
\put(4965,3939){\blacken\ellipse{90}{90}}
\put(4965,3939){\ellipse{90}{90}}
\put(1375,1779){\blacken\ellipse{90}{90}}
\put(1375,1779){\ellipse{90}{90}}
\put(2319,2891){\blacken\ellipse{90}{90}}
\put(2319,2891){\ellipse{90}{90}}
\put(3579,3576){\blacken\ellipse{90}{90}}
\put(3579,3576){\ellipse{90}{90}}
\thinlines
\blacken\thicklines
\path(895.000,4194.000)(850,4344.000)(805.000,4194.000)(895.000,4194.000)
\path(850,4344)(850,675)(5950,675)  
\blacken\path(5800.000,630.000)(5950.000,675.000)(5800.000,720.000)(5800.000,630.000)
\thinlines

\dashline{60.000}(4225,3745)(4225,675)
\dashline{60.000}(1217,1446)(1217,675)
\thicklines
\path(850,675)(1375,1779)
	(2319,2891)(3579,3576)
(4965,3939)(5550,3984)(5905,3984)
\thinlines
\dottedline{90.000}(5050,4375)(850,2092)
\dashline{60.000}(5050,4375)(5050,675)
\thicklines
\path(5050,4375)(850,1166) 
\put(1000,200){\makebox(0,0)[lb]{\smash{{\SetFigFont{8}{8}{\rmdefault}{\mddefault}{\updefault}$\lambda_{\mathrm{inf}}$}}}}
\put(3750,200){\makebox(0,0)[lb]{\smash{{\SetFigFont{8}{8}{\rmdefault}{\mddefault}{\updefault}$\lambda_{\mathrm{sup}}$}}}}
\put(950,4050){\makebox(0,0)[lb]{\smash{{\SetFigFont{8}{8}{\rmdefault}{\mddefault}{\updefault}$\hat{\Jc}(S;\lambda)$}}}}
\put(2600,3459){\makebox(0,0)[lb]{\smash{{\SetFigFont{8}{8}{\rmdefault}{\mddefault}{\updefault}$S_j$}}}}
\put(4800,200){\makebox(0,0)[lb]{\smash{{\SetFigFont{8}{8}{\rmdefault}{\mddefault}{\updefault}$\delta\Ec_{\mathrm{add}}$}}}}
\put(2450,2089){\makebox(0,0)[lb]{\smash{{\SetFigFont{8}{8}{\rmdefault}{\mddefault}{\updefault}$S_{\mathrm{new}}$}}}}
\put(1217,1446){\whiten\ellipse{150}{150}}
\put(1217,1446){\ellipse{150}{150}}
\put(4225,3745){\whiten\ellipse{150}{150}}
\put(4225,3745){\ellipse{150}{150}}
\end{picture}
}
\end{tabular}&
  \;\begin{tabular}{c}
\setlength{\unitlength}{0.00023in}
\begingroup\makeatletter\ifx\SetFigFont\undefined%
\gdef\SetFigFont#1#2#3#4#5{%
  \reset@font\fontsize{#1}{#2pt}%
  \fontfamily{#3}\fontseries{#4}\fontshape{#5}%
  \selectfont}%
\fi\endgroup%
{\renewcommand{\dashlinestretch}{30}
\begin{picture}(6800,6143)(840,970)
\thicklines
\put(6518,5721){\blacken\ellipse{113}{113}}
\put(6518,5721){\ellipse{113}{113}}
\put(1635,3309){\blacken\ellipse{113}{113}}
\put(1635,3309){\ellipse{113}{113}}
\put(1838,3731){\blacken\ellipse{113}{113}}
\put(1838,3731){\ellipse{113}{113}}
\put(5288,5613){\blacken\ellipse{113}{113}}
\put(5288,5613){\ellipse{113}{113}}
\put(3300,5091){\blacken\ellipse{113}{113}}
\put(3300,5091){\ellipse{113}{113}}
\blacken\path(1185.000,6309.000)(1140.000,6459.000)(1095.000,6309.000)(1185.000,6309.000)
\path(1140,6459)(1140,1824)(7265,1824)  
\blacken\path(7115.000,1779.000)(7265.000,1824.000)(7115.000,1869.000)(7115.000,1779.000)
\thinlines
\dottedline{130.000}(4647,6344)(1140,3082)
\dashline{90.000}(5629,5643)(5629,1824)  
\dashline{90.000}(2451,4301)(2451,1824)
\thicklines
\path(7035,6237)(1140,3748)    
\path(1140,1824)(1635,3309)(1838,3731)
(3300,5091)
(5288,5613)
(6518,5721)(7128,5721)
\put(1000,1200){\makebox(0,0)[lb]{\smash{{\SetFigFont{8}{8.0}{\rmdefault}{\mddefault}{\updefault}0}}}}
\put(6950,1200){\makebox(0,0)[lb]{\smash{{\SetFigFont{8}{8.0}{\rmdefault}{\mddefault}{\updefault}$\lambda$}}}}
\put(1810,1200){\makebox(0,0)[lb]{\smash{{\SetFigFont{8}{8.0}{\rmdefault}{\mddefault}{\updefault}$\lambda_{\mathrm{inf}}=\delta\Ec_{\mathrm{rmv}}$}}}}
\put(5265,1200){\makebox(0,0)[lb]{\smash{{\SetFigFont{8}{8.0}{\rmdefault}{\mddefault}{\updefault}$\lambda_{\mathrm{sup}}$}}}}
\put(3600,4400){\makebox(0,0)[lb]{\smash{{\SetFigFont{8}{8.0}{\rmdefault}{\mddefault}{\updefault}$S_{\mathrm{new}}$}}}}
\put(2300,4850){\makebox(0,0)[lb]{\smash{{\SetFigFont{8}{8.0}{\rmdefault}{\mddefault}{\updefault}$S_j$}}}}
\put(1280,6020){\makebox(0,0)[lb]{\smash{{\SetFigFont{8}{8.0}{\rmdefault}{\mddefault}{\updefault}$\hat{\Jc}(S;\lambda)$}}}}
\put(2451,4301){\whiten\ellipse{189}{189}}
\put(2451,4301){\ellipse{189}{189}}
\put(5629,5643){\whiten\ellipse{189}{189}}
\put(5629,5643){\ellipse{189}{189}}
\end{picture}
}
\end{tabular}\\
  \begin{tabular}{c}\small{(c)}\end{tabular}&
  \begin{tabular}{c}
\setlength{\unitlength}{0.00029in}
\begingroup\makeatletter\ifx\SetFigFont\undefined%
\gdef\SetFigFont#1#2#3#4#5{%
  \reset@font\fontsize{#1}{#2pt}%
  \fontfamily{#3}\fontseries{#4}\fontshape{#5}%
  \selectfont}%
\fi\endgroup%
{\renewcommand{\dashlinestretch}{30}
\begin{picture}(5100,4881)(745,100)
\thicklines
\put(5550,3984){\blacken\ellipse{90}{90}}
\put(5550,3984){\ellipse{90}{90}}
\put(4965,3939){\blacken\ellipse{90}{90}}
\put(4965,3939){\ellipse{90}{90}}
\thinlines
\dashline{60.000}(1217,1446)(1217,675)
\blacken\thicklines
\path(895.000,4194.000)(850,4344.000)(805.000,4194.000)(895.000,4194.000)
\path(850,4344)(850,675)(5950,675)  
\blacken\path(5800.000,630.000)(5950.000,675.000)(5800.000,720.000)(5800.000,630.000)
\thinlines
\thicklines
\thinlines
\dashline{60.000}(4225,3745)(4225,675)
\thicklines
\path(850,675)(1217,1446)(4225,3745)(4965,3939)(5550,3984)(5905,3984)
\put(930,200){\makebox(0,0)[lb]{\smash{{\SetFigFont{8}{8}{\rmdefault}{\mddefault}{\updefault}$\lambda_{j+1}$}}}}
\put(4100,200){\makebox(0,0)[lb]{\smash{{\SetFigFont{8}{8}{\rmdefault}{\mddefault}{\updefault}$\lambda_j$}}}}
\put(950,4000){\makebox(0,0)[lb]{\smash{{\SetFigFont{8}{8}{\rmdefault}{\mddefault}{\updefault}$\hat{\Jc}(S;\lambda)$}}}}
\put(2300,2759){\makebox(0,0)[lb]{\smash{{\SetFigFont{8}{8}{\rmdefault}{\mddefault}{\updefault}$S_j$}}}}
\put(5700,250){\makebox(0,0)[lb]{\smash{{\SetFigFont{8}{8}{\rmdefault}{\mddefault}{\updefault}$\lambda$}}}}
\put(1217,1446){\blacken\ellipse{90}{90}}
\put(1217,1446){\ellipse{90}{90}}
\put(4225,3745){\blacken\ellipse{90}{90}}
\put(4225,3745){\ellipse{90}{90}}
\end{picture}
}
\end{tabular}&
  \;\begin{tabular}{c}
\setlength{\unitlength}{0.00023in}
\begingroup\makeatletter\ifx\SetFigFont\undefined%
\gdef\SetFigFont#1#2#3#4#5{%
  \reset@font\fontsize{#1}{#2pt}%
  \fontfamily{#3}\fontseries{#4}\fontshape{#5}%
  \selectfont}%
\fi\endgroup%
{\renewcommand{\dashlinestretch}{30}
\begin{picture}(6800,6200)(840,1060)
\thicklines
\put(6518,5721){\blacken\ellipse{113}{113}}
\put(6518,5721){\ellipse{113}{113}}
\put(1635,3309){\blacken\ellipse{113}{113}}
\put(1635,3309){\ellipse{113}{113}}
\put(1838,3731){\blacken\ellipse{113}{113}}
\put(1838,3731){\ellipse{113}{113}}
\blacken\path(1185.000,6309.000)(1140.000,6459.000)(1095.000,6309.000)(1185.000,6309.000)
\path(1140,6459)(1140,1824)(7265,1824)  
\blacken\path(7115.000,1779.000)(7265.000,1824.000)(7115.000,1869.000)(7115.000,1779.000)
\thinlines
\dashline{90.000}(5629,5643)(5629,1824) 
\dashline{90.000}(2451,4301)(2451,1824) 
\thicklines
\path(1140,1824)(1635,3309)(1838,3731)
(2451,4301)
(5629,5643) 
(6518,5721)(7128,5721)
\put(1000,1269){\makebox(0,0)[lb]{\smash{{\SetFigFont{8}{8.0}{\rmdefault}{\mddefault}{\updefault}0}}}}
\put(6950,1269){\makebox(0,0)[lb]{\smash{{\SetFigFont{8}{8.0}{\rmdefault}{\mddefault}{\updefault}$\lambda$}}}}
\put(2110,1269){\makebox(0,0)[lb]{\smash{{\SetFigFont{8}{8.0}{\rmdefault}{\mddefault}{\updefault}$\lambda_{j+1}$}}}}
\put(5430,1269){\makebox(0,0)[lb]{\smash{{\SetFigFont{8}{8.0}{\rmdefault}{\mddefault}{\updefault}$\lambda_j$}}}}
\put(3480,5124){\makebox(0,0)[lb]{\smash{{\SetFigFont{8}{8.0}{\rmdefault}{\mddefault}{\updefault}$S_j$}}}}
\put(1280,6020){\makebox(0,0)[lb]{\smash{{\SetFigFont{8}{8.0}{\rmdefault}{\mddefault}{\updefault}$\hat{\Jc}(S;\lambda)$}}}}
\put(2451,4301){\blacken\ellipse{113}{113}}
\put(2451,4301){\ellipse{113}{113}}
\put(5629,5643){\blacken\ellipse{113}{113}}
\put(5629,5643){\ellipse{113}{113}}
\end{picture}
}
\end{tabular}
  \end{tabular}
}
\DRAFT{
  \setlength{\tabcolsep}{0.25cm}
  \begin{center}
  \begin{tabular}{ccc}
  \begin{tabular}{c}
\setlength{\unitlength}{0.00029in}
\begingroup\makeatletter\ifx\SetFigFont\undefined%
\gdef\SetFigFont#1#2#3#4#5{%
  \reset@font\fontsize{#1}{#2pt}%
  \fontfamily{#3}\fontseries{#4}\fontshape{#5}%
  \selectfont}%
\fi\endgroup%
{\renewcommand{\dashlinestretch}{30}
\begin{picture}(5100,4681)(750,0)
\thicklines
\put(5550,3984){\blacken\ellipse{90}{90}}
\put(5550,3984){\ellipse{90}{90}}
\put(4965,3939){\blacken\ellipse{90}{90}}
\put(4965,3939){\ellipse{90}{90}}
\put(1375,1779){\blacken\ellipse{90}{90}}
\put(1375,1779){\ellipse{90}{90}}
\put(2319,2891){\blacken\ellipse{90}{90}}
\put(2319,2891){\ellipse{90}{90}}
\put(3579,3576){\blacken\ellipse{90}{90}}
\put(3579,3576){\ellipse{90}{90}}
\thinlines
\dashline{60.000}(3579,3576)(3579,675)
\blacken\thicklines
\path(895.000,4194.000)(850,4344.000)(805.000,4194.000)(895.000,4194.000)
\path(850,4344)(850,675)(5950,675)  
\blacken\path(5800.000,630.000)(5950.000,675.000)(5800.000,720.000)(5800.000,630.000)
\thinlines
\dashline{60.000}(2319,2891)(2319,675)
\thicklines
\path(850,675)(1375,1779)
	(2319,2891)(3579,3576)(4965,3939)
	(5550,3984)(5905,3984)
\thinlines
\thicklines
\put(810,200){\makebox(0,0)[lb]{\smash{{\SetFigFont{8}{8}{\rmdefault}{\mddefault}{\updefault}0}}}}
\put(2000,200){\makebox(0,0)[lb]{\smash{{\SetFigFont{8}{8}{\rmdefault}{\mddefault}{\updefault}$\lambda_{j+1}$}}}}
\put(3450,200){\makebox(0,0)[lb]{\smash{{\SetFigFont{8}{8}{\rmdefault}{\mddefault}{\updefault}$\lambda_j$}}}}
\put(950,4000){\makebox(0,0)[lb]{\smash{{\SetFigFont{8}{8}{\rmdefault}{\mddefault}{\updefault}$\hat{\Jc}(S;\lambda)$}}}}
\put(2600,3459){\makebox(0,0)[lb]{\smash{{\SetFigFont{8}{8}{\rmdefault}{\mddefault}{\updefault}$S_j$}}}}
\put(5700,200){\makebox(0,0)[lb]{\smash{{\SetFigFont{8}{8}{\rmdefault}{\mddefault}{\updefault}$\lambda$}}}}
\end{picture}
}
\end{tabular}&
  \begin{tabular}{c}
\setlength{\unitlength}{0.00029in}
\begingroup\makeatletter\ifx\SetFigFont\undefined%
\gdef\SetFigFont#1#2#3#4#5{%
  \reset@font\fontsize{#1}{#2pt}%
  \fontfamily{#3}\fontseries{#4}\fontshape{#5}%
  \selectfont}%
\fi\endgroup%
{\renewcommand{\dashlinestretch}{30}
\begin{picture}(5100,4881)(750,0)
\thicklines
\put(5550,3984){\blacken\ellipse{90}{90}}
\put(5550,3984){\ellipse{90}{90}}
\put(4965,3939){\blacken\ellipse{90}{90}}
\put(4965,3939){\ellipse{90}{90}}
\put(1375,1779){\blacken\ellipse{90}{90}}
\put(1375,1779){\ellipse{90}{90}}
\put(2319,2891){\blacken\ellipse{90}{90}}
\put(2319,2891){\ellipse{90}{90}}
\put(3579,3576){\blacken\ellipse{90}{90}}
\put(3579,3576){\ellipse{90}{90}}
\thinlines
\blacken\thicklines
\path(895.000,4194.000)(850,4344.000)(805.000,4194.000)(895.000,4194.000)
\path(850,4344)(850,675)(5950,675) 
\blacken\path(5800.000,630.000)(5950.000,675.000)(5800.000,720.000)(5800.000,630.000)
\thinlines
\dashline{60.000}(4225,3745)(4225,675)
\dashline{60.000}(1217,1446)(1217,675)
\thicklines
\path(850,675)(1375,1779)
	(2319,2891)(3579,3576)
(4965,3939)(5550,3984)(5905,3984)
\thinlines
\dottedline{90.000}(5050,4375)(850,2092)
\dashline{60.000}(5050,4375)(5050,675)
\thicklines
\path(5050,4375)(850,1166) 
\put(1000,200){\makebox(0,0)[lb]{\smash{{\SetFigFont{8}{8}{\rmdefault}{\mddefault}{\updefault}$\lambda_{\mathrm{inf}}$}}}}
\put(3750,200){\makebox(0,0)[lb]{\smash{{\SetFigFont{8}{8}{\rmdefault}{\mddefault}{\updefault}$\lambda_{\mathrm{sup}}$}}}}
\put(950,4050){\makebox(0,0)[lb]{\smash{{\SetFigFont{8}{8}{\rmdefault}{\mddefault}{\updefault}$\hat{\Jc}(S;\lambda)$}}}}
\put(2600,3459){\makebox(0,0)[lb]{\smash{{\SetFigFont{8}{8}{\rmdefault}{\mddefault}{\updefault}$S_j$}}}}
\put(4800,200){\makebox(0,0)[lb]{\smash{{\SetFigFont{8}{8}{\rmdefault}{\mddefault}{\updefault}$\delta\Ec_{\mathrm{add}}$}}}}
\put(2450,2089){\makebox(0,0)[lb]{\smash{{\SetFigFont{8}{8}{\rmdefault}{\mddefault}{\updefault}$S_{\mathrm{new}}$}}}}
\put(1217,1446){\whiten\ellipse{150}{150}}
\put(1217,1446){\ellipse{150}{150}}
\put(4225,3745){\whiten\ellipse{150}{150}}
\put(4225,3745){\ellipse{150}{150}}
\end{picture}
}

\end{tabular}&
  \begin{tabular}{c}
\setlength{\unitlength}{0.00029in}
\begingroup\makeatletter\ifx\SetFigFont\undefined%
\gdef\SetFigFont#1#2#3#4#5{%
  \reset@font\fontsize{#1}{#2pt}%
  \fontfamily{#3}\fontseries{#4}\fontshape{#5}%
  \selectfont}%
\fi\endgroup%
{\renewcommand{\dashlinestretch}{30}
\begin{picture}(5100,4881)(745,100)
\thicklines
\put(5550,3984){\blacken\ellipse{90}{90}}
\put(5550,3984){\ellipse{90}{90}}
\put(4965,3939){\blacken\ellipse{90}{90}}
\put(4965,3939){\ellipse{90}{90}}
\thinlines
\dashline{60.000}(1217,1446)(1217,675)
\blacken\thicklines
\path(895.000,4194.000)(850,4344.000)(805.000,4194.000)(895.000,4194.000)
\path(850,4344)(850,675)(5950,675)  
\blacken\path(5800.000,630.000)(5950.000,675.000)(5800.000,720.000)(5800.000,630.000)
\thinlines
\thicklines
\thinlines
\dashline{60.000}(4225,3745)(4225,675)
\thicklines
\path(850,675)(1217,1446)(4225,3745)(4965,3939)(5550,3984)(5905,3984)
\put(930,200){\makebox(0,0)[lb]{\smash{{\SetFigFont{8}{8}{\rmdefault}{\mddefault}{\updefault}$\lambda_{j+1}$}}}}
\put(4100,200){\makebox(0,0)[lb]{\smash{{\SetFigFont{8}{8}{\rmdefault}{\mddefault}{\updefault}$\lambda_j$}}}}
\put(950,4000){\makebox(0,0)[lb]{\smash{{\SetFigFont{8}{8}{\rmdefault}{\mddefault}{\updefault}$\hat{\Jc}(S;\lambda)$}}}}
\put(2300,2759){\makebox(0,0)[lb]{\smash{{\SetFigFont{8}{8}{\rmdefault}{\mddefault}{\updefault}$S_j$}}}}
\put(5700,250){\makebox(0,0)[lb]{\smash{{\SetFigFont{8}{8}{\rmdefault}{\mddefault}{\updefault}$\lambda$}}}}
\put(1217,1446){\blacken\ellipse{90}{90}}
\put(1217,1446){\ellipse{90}{90}}
\put(4225,3745){\blacken\ellipse{90}{90}}
\put(4225,3745){\ellipse{90}{90}}
\end{picture}
}
\end{tabular}\\
  (a) & (b) $S_{\mathrm{new}}=S_j+\{\ell_\mathrm{add}\}$ & (c) \\[.4cm]
  \begin{tabular}{c}
\setlength{\unitlength}{0.00023in}
\begingroup\makeatletter\ifx\SetFigFont\undefined%
\gdef\SetFigFont#1#2#3#4#5{%
  \reset@font\fontsize{#1}{#2pt}%
  \fontfamily{#3}\fontseries{#4}\fontshape{#5}%
  \selectfont}%
\fi\endgroup%
{\renewcommand{\dashlinestretch}{30}
\begin{picture}(6800,5943)(840,970)
\thicklines
\put(6518,5721){\blacken\ellipse{113}{113}}
\put(6518,5721){\ellipse{113}{113}}
\put(1635,3309){\blacken\ellipse{113}{113}}
\put(1635,3309){\ellipse{113}{113}}
\put(1838,3731){\blacken\ellipse{113}{113}}
\put(1838,3731){\ellipse{113}{113}}
\put(5288,5613){\blacken\ellipse{113}{113}}
\put(5288,5613){\ellipse{113}{113}}
\put(3300,5091){\blacken\ellipse{113}{113}}
\put(3300,5091){\ellipse{113}{113}}
\blacken\path(1185.000,6309.000)(1140.000,6459.000)(1095.000,6309.000)(1185.000,6309.000)
\path(1140,6459)(1140,1824)(7265,1824)  
\blacken\path(7115.000,1779.000)(7265.000,1824.000)(7115.000,1869.000)(7115.000,1779.000)
\thinlines
\dashline{90.000}(1838,3731)(1838,1824)  
\dashline{90.000}(3300,5091)(3300,1824)  
\thicklines
\path(1140,1824)(1635,3309)(1838,3731)
(3300,5091)(5288,5613)
(6518,5721)(7128,5721)
\put(950,1269){\makebox(0,0)[lb]{\smash{{\SetFigFont{8}{8.0}{\rmdefault}{\mddefault}{\updefault}0}}}}
\put(6950,1269){\makebox(0,0)[lb]{\smash{{\SetFigFont{8}{8.0}{\rmdefault}{\mddefault}{\updefault}$\lambda$}}}}
\put(1510,1269){\makebox(0,0)[lb]{\smash{{\SetFigFont{8}{8.0}{\rmdefault}{\mddefault}{\updefault}$\lambda_{j+1}$}}}}
\put(3100,1269){\makebox(0,0)[lb]{\smash{{\SetFigFont{8}{8.0}{\rmdefault}{\mddefault}{\updefault}$\lambda_j$}}}}
\put(1980,4524){\makebox(0,0)[lb]{\smash{{\SetFigFont{8}{8.0}{\rmdefault}{\mddefault}{\updefault}$S_j$}}}}
\put(1280,6020){\makebox(0,0)[lb]{\smash{{\SetFigFont{8}{8.0}{\rmdefault}{\mddefault}{\updefault}$\hat{\Jc}(S;\lambda)$}}}}
\end{picture}
}
\end{tabular}&
  \begin{tabular}{c}
\setlength{\unitlength}{0.00023in}
\begingroup\makeatletter\ifx\SetFigFont\undefined%
\gdef\SetFigFont#1#2#3#4#5{%
  \reset@font\fontsize{#1}{#2pt}%
  \fontfamily{#3}\fontseries{#4}\fontshape{#5}%
  \selectfont}%
\fi\endgroup%
{\renewcommand{\dashlinestretch}{30}
\begin{picture}(6800,6143)(840,970)
\thicklines
\put(6518,5721){\blacken\ellipse{113}{113}}
\put(6518,5721){\ellipse{113}{113}}
\put(1635,3309){\blacken\ellipse{113}{113}}
\put(1635,3309){\ellipse{113}{113}}
\put(1838,3731){\blacken\ellipse{113}{113}}
\put(1838,3731){\ellipse{113}{113}}
\put(5288,5613){\blacken\ellipse{113}{113}}
\put(5288,5613){\ellipse{113}{113}}
\put(3300,5091){\blacken\ellipse{113}{113}}
\put(3300,5091){\ellipse{113}{113}}
\blacken\path(1185.000,6309.000)(1140.000,6459.000)(1095.000,6309.000)(1185.000,6309.000)
\path(1140,6459)(1140,1824)(7265,1824)  
\blacken\path(7115.000,1779.000)(7265.000,1824.000)(7115.000,1869.000)(7115.000,1779.000)
\thinlines
\dottedline{130.000}(4647,6344)(1140,3082)   
\dashline{90.000}(5629,5643)(5629,1824)  
\dashline{90.000}(2451,4301)(2451,1824)
\thicklines
\path(7035,6237)(1140,3748)    
\path(1140,1824)(1635,3309)(1838,3731)
(3300,5091)
(5288,5613)
(6518,5721)(7128,5721)
\put(1000,1200){\makebox(0,0)[lb]{\smash{{\SetFigFont{8}{8.0}{\rmdefault}{\mddefault}{\updefault}0}}}}
\put(6950,1200){\makebox(0,0)[lb]{\smash{{\SetFigFont{8}{8.0}{\rmdefault}{\mddefault}{\updefault}$\lambda$}}}}
\put(1810,1200){\makebox(0,0)[lb]{\smash{{\SetFigFont{8}{8.0}{\rmdefault}{\mddefault}{\updefault}$\lambda_{\mathrm{inf}}=\delta\Ec_{\mathrm{rmv}}$}}}}
\put(5265,1200){\makebox(0,0)[lb]{\smash{{\SetFigFont{8}{8.0}{\rmdefault}{\mddefault}{\updefault}$\lambda_{\mathrm{sup}}$}}}}
\put(3600,4400){\makebox(0,0)[lb]{\smash{{\SetFigFont{8}{8.0}{\rmdefault}{\mddefault}{\updefault}$S_{\mathrm{new}}$}}}}
\put(2300,4850){\makebox(0,0)[lb]{\smash{{\SetFigFont{8}{8.0}{\rmdefault}{\mddefault}{\updefault}$S_j$}}}}
\put(1280,6020){\makebox(0,0)[lb]{\smash{{\SetFigFont{8}{8.0}{\rmdefault}{\mddefault}{\updefault}$\hat{\Jc}(S;\lambda)$}}}}
\put(2451,4301){\whiten\ellipse{189}{189}}
\put(2451,4301){\ellipse{189}{189}}
\put(5629,5643){\whiten\ellipse{189}{189}}
\put(5629,5643){\ellipse{189}{189}}
\end{picture}
}
\end{tabular}&
  \begin{tabular}{c}
\setlength{\unitlength}{0.00023in}
\begingroup\makeatletter\ifx\SetFigFont\undefined%
\gdef\SetFigFont#1#2#3#4#5{%
  \reset@font\fontsize{#1}{#2pt}%
  \fontfamily{#3}\fontseries{#4}\fontshape{#5}%
  \selectfont}%
\fi\endgroup%
{\renewcommand{\dashlinestretch}{30}
\begin{picture}(6800,6200)(840,1060)
\thicklines
\put(6518,5721){\blacken\ellipse{113}{113}}
\put(6518,5721){\ellipse{113}{113}}
\put(1635,3309){\blacken\ellipse{113}{113}}
\put(1635,3309){\ellipse{113}{113}}
\put(1838,3731){\blacken\ellipse{113}{113}}
\put(1838,3731){\ellipse{113}{113}}
\blacken\path(1185.000,6309.000)(1140.000,6459.000)(1095.000,6309.000)(1185.000,6309.000)
\path(1140,6459)(1140,1824)(7265,1824)  
\blacken\path(7115.000,1779.000)(7265.000,1824.000)(7115.000,1869.000)(7115.000,1779.000)
\thinlines
\dashline{90.000}(5629,5643)(5629,1824) 
\dashline{90.000}(2451,4301)(2451,1824) 
\thicklines
\path(1140,1824)(1635,3309)(1838,3731)
(2451,4301)
(5629,5643) 
(6518,5721)(7128,5721)
\put(1000,1269){\makebox(0,0)[lb]{\smash{{\SetFigFont{8}{8.0}{\rmdefault}{\mddefault}{\updefault}0}}}}
\put(6950,1269){\makebox(0,0)[lb]{\smash{{\SetFigFont{8}{8.0}{\rmdefault}{\mddefault}{\updefault}$\lambda$}}}}
\put(2110,1269){\makebox(0,0)[lb]{\smash{{\SetFigFont{8}{8.0}{\rmdefault}{\mddefault}{\updefault}$\lambda_{j+1}$}}}}
\put(5430,1269){\makebox(0,0)[lb]{\smash{{\SetFigFont{8}{8.0}{\rmdefault}{\mddefault}{\updefault}$\lambda_j$}}}}
\put(3480,5124){\makebox(0,0)[lb]{\smash{{\SetFigFont{8}{8.0}{\rmdefault}{\mddefault}{\updefault}$S_j$}}}}
\put(1280,6020){\makebox(0,0)[lb]{\smash{{\SetFigFont{8}{8.0}{\rmdefault}{\mddefault}{\updefault}$\hat{\Jc}(S;\lambda)$}}}}
\put(2451,4301){\blacken\ellipse{113}{113}}
\put(2451,4301){\ellipse{113}{113}}
\put(5629,5643){\blacken\ellipse{113}{113}}
\put(5629,5643){\ellipse{113}{113}}
\end{picture}
}
\end{tabular}\\
  (a) & (b) $S_{\mathrm{new}}=S_j-\{\ell_\mathrm{rmv}\}$& (c) 
  \end{tabular}
  \end{center}
}
\caption{$\ell_0$-PD algorithm: descent of the concave polygon when a
  new support $S_{\mathrm{new}}=S_j+\{\ell_\mathrm{add}\}$
  (\FINAL{left}\DRAFT{top}) or
  $S_{\mathrm{new}}=S_j-\{\ell_\mathrm{rmv}\}$
  (\FINAL{right}\DRAFT{bottom}) is included. \;(a)~Initial
  configuration.\;(b)~The intersection with line $S_{\mathrm{new}}$ is
  computed. This yields an interval
  $\stdcro{\lambda_{\mathrm{inf}},\lambda_{\mathrm{sup}}}$ for which
  $S_{\mathrm{new}}$ lays below the concave polygon.\;(c) When this
  interval is non-empty, the supports $S_j$ whose related edges lay
  above the line $S_{\mathrm{new}}$ are removed while
  $S_{\mathrm{new}}$ is included in \Sc. The list of values
  $\lambda_j$ (corresponding to the vertices of the new concave
  polygon) is being updated. Their number may either increase or
  decrease.  }
\label{fig:refinement}
\end{figure}

We propose to gradually update a list \Sc of candidate subsets $S_j$
while imposing that the related curve is a concave polygon, obtained
as the concave envelope of the set of lines $S_j$ (see
Fig.~\ref{fig:refinement}(a)). The subsets in \Sc are updated so as to
decrease at most the concave polygonal curve.  In particular, we
impose that the least value is $\lambda_{J+1}=0$, so that the concave
envelope is computed over the whole domain $\lambda\in\Rbb_+$.

\subsection{Descent of the concave polygon}
The principle of $\ell_0$-PD is to perform a series of descent steps,
where a new candidate subset $S_{\mathrm{new}}$ is considered and
included in the list \Sc only if the resulting concave polygon can be
decreased. This descent test is illustrated on
Fig.~\ref{fig:refinement} for two examples (\FINAL{left and right
  columns}\DRAFT{top and bottom subfigures}). For each example, the
initial polygon is represented in~(a). It is updated when its
intersection with the line $S_{\mathrm{new}}$ is non-empty (b). The
new concave polygon~(c) is obtained as the concave envelope of the
former polygon and the line $S_{\mathrm{new}}$. All subsets in \Sc
whose edges lay above the line $S_{\mathrm{new}}$ are removed from
\Sc.

This procedure is formally presented in Tab.~\ref{tab:ccv}. Let us now
specify how the new candidate subsets $S_{\mathrm{new}}$ are built.
\begin{table}[t]
  \caption{Concave polygon descent procedure. When a new subset
      is included, both lists \Sc and \lambdab are updated. The
      function $\mathtt{intersect}$ computes the intersection
      between a line and a concave polygon. This yields an interval
      $\stdcro{\lambda_{\mathrm{inf}},\lambda_{\mathrm{sup}}}$. By convention,
      $\lambda_{\mathrm{inf}}>\lambda_{\mathrm{sup}}$ when the
      intersection is empty.  
}
\label{tab:ccv}
\centering 
\begin{tabular}{l} 
\hline
\textbf{Procedure}: $\mathtt{CCV\_Descent}(\Sc$, $S_\mathrm{new}$, $\lambdab$)\\ 
\hline
Call $\stdcro{\lambda_{\mathrm{inf}},\lambda_{\mathrm{sup}}}=
\mathtt{intersect}(\Sc$, $S_\mathrm{new})$;\\
\textbf{if} $\lambda_{\mathrm{inf}}<\lambda_{\mathrm{sup}}$ \textbf{then}\\[.04cm]
\hspace*{.4cm}Include $S_{\mathrm{new}}$ as an unexplored support in \Sc;\\
\hspace*{.4cm}Remove any subset $S_j$ from \Sc
such that $[\lambda_{j+1},\lambda_j]\subset[\lambda_{\mathrm{inf}},
\lambda_{\mathrm{sup}}]$;\\
\hspace*{.4cm}Sort the subsets in \Sc by increasing cardinality;\\
\hspace*{.4cm}Sort $\lambdab$ in the decreasing order;\\
\textbf{end}\\
\hline
\end{tabular}
\end{table}

\subsection{Selection of the new candidate support}
\label{sec:selection_rule}
We first need to assign a Boolean label $ S_j.\mathtt{expl}$ to
each subset $S_j$. It equals 1 if $S_j$ has already been ``explored''
and 0 otherwise. The following exploration process is being carried
out given a subset $S=S_j$: all the possible single
replacements $S\pm\{i\}$ are tested. The best insertion
$\ell_\mathrm{add}$ and removal $\ell_\mathrm{rmv}$ are both
kept in memory, with $\ell_\mathrm{add}$ defined in~\eqref{eq:lplus}
and similarly,
\begin{align}
\ell_\mathrm{rmv}&\in\argmin_{i\in S}\stdacc{\Ec(S-\{i\})-\Ec(S)}.
\label{eq:lmoins}
\end{align}

At any $\ell_0$-PD iteration, the unexplored subset $S_j$ of lowest
cardinality (\ie of lowest index $j$) is selected. $\ell_0$-PD
attempts to include $S_{\mathrm{add}}=S_j+\{\ell_\mathrm{add}\}$ and
$S_{\mathrm{rmv}}=S_j-\{\ell_\mathrm{rmv}\}$ into \Sc, so that the
concave polygon can be decreased at most. The $\mathtt{CCV\_Descent}$
procedure (Tab.~\ref{tab:ccv}) is first called with
$S_{\mathrm{new}}\leftarrow S_{\mathrm{add}}$ leading to possible
updates of \Sc and \lambdab. It is called again with
$S_{\mathrm{new}}\leftarrow
S_{\mathrm{rmv}}$. Fig.~\ref{fig:refinement} illustrates each of these
calls: the slope of $S_{\mathrm{new}}$ is $\stdbars{S_j}+1$ and
$\stdbars{S_j}-1$, respectively. When a support $S_j$ has been
explored, the new supports that have been included in \Sc (if any) are
tagged as unexplored.

\subsection{$\ell_0$-PD algorithm}
$\ell_0$-PD is stated in Tab.~\ref{tab:l0PD}. Initially, \Sc is formed
of the empty support $S_0=\emptyset$. The resulting concave polygon is
reduced to a single horizontal edge. The corresponding endpoints are
$\lambda_1=0$ and (by extension) $\lambda_0\triangleq+\infty$. In the
first iteration, $S_0$ is explored: the best insertion
$S_{\mathrm{add}}=\{\ell_\mathrm{add}\}$ is computed
in~\eqref{eq:lplus_init}, and included in \Sc during the call to
$\mathtt{CCV\_Descent}$. The updated set \Sc is now composed of
$S_0=\emptyset$ (explored) and $S_1=S_{\mathrm{add}}$
(unexplored). The new concave polygon has two edges delimited by
$\lambda_2=0$, $\lambda_1$ and $\lambda_0=+\infty$, with $\lambda_1$
given in~\eqref{eq:lplus_init}. Generally, either 0, 1, or 2 new
unexplored supports $S_{\mathrm{add}}$ and $S_{\mathrm{rmv}}$ may be
included in \Sc at a given iteration while a variable number of
supports may be removed from \Sc.
\begin{table}[t]
  \caption{$\ell_0$-PD algorithm. The algorithm maintains a
    list $\Sc$ of supports $S_j$ whose cardinality is increasing
    with $j$. The unexplored support having the lowest cardinality
    is explored at each iteration. The lists \Sc and
    $\lambdab$ are updated during the calls to
    $\mathtt{CCV\_Descent}$; \lambdab is sorted in the
    decreasing order, with $\lambda_{J+1}=0$. During the first 
    iteration, $j=0$ leads to $S_{\mathrm{rmv}}\leftarrow \emptyset$.
  }
  \label{tab:l0PD}
\centering 
\begin{tabular}{l} 
\hline
\textbf{inputs}~~: \Ab, \yb \\
\textbf{outputs}: \Sc, \lambdab \\
\hline
$\lambdab\leftarrow\stdacc{\lambda_1}$ with $\lambda_1\leftarrow 0$;\\[.04cm]
$S_0\leftarrow \emptyset$, $S_0.\mathtt{expl}\leftarrow 0$;\\[.04cm]
$\Sc\leftarrow \{S_0\}$;\\[.04cm]
\textbf{while} $\{\exists j:\,S_j.\mathtt{expl}=0\}$ \textbf{do}\\[.04cm]
\hspace*{.5cm}Set $j$ as the lowest index such that $S_j.\mathtt{expl}=0$;\\
\hspace*{.5cm}$S_j.\mathtt{expl}\leftarrow 1$;\\
\hspace*{.5cm}Compute $S_{\mathrm{add}}\leftarrow S_j+\{\ell_{\mathrm{add}}\}$
from~\eqref{eq:lplus};\\
\hspace*{.5cm}\textbf{if} $j=0$ \textbf{then}\\[.04cm]
\hspace*{1.0cm}$S_{\mathrm{rmv}}\leftarrow \emptyset$;\\
\hspace*{.5cm}\textbf{else}\\
\hspace*{1.0cm}Compute $S_{\mathrm{rmv}}\leftarrow S_j-\{\ell_{\mathrm{rmv}}\}$ 
from \eqref{eq:lmoins};\\
\hspace*{.5cm}\textbf{end}\\
\hspace*{.5cm}Call $\mathtt{CCV\_Descent}(\Sc$, $S_{\mathrm{add}}$, $\lambdab)$;\\
\hspace*{.5cm}Call $\mathtt{CCV\_Descent}(\Sc$, $S_{\mathrm{rmv}}$, $\lambdab)$;\\
\textbf{end}\\
\hline
\end{tabular}
\end{table}

$\ell_0$-PD terminates when all supports in \Sc have been
explored. When this occurs, the concave polygon cannot decrease
anymore with any single replacement $S_j\pm\{i\}$, with $S_j\in\Sc$.
Practically, the early stopping rule
$\lambda_j\leq\lambda_{\mathrm{stop}}$ can be adopted, where $j$
denotes the unexplored subset having the least cardinality. This rule
ensures that all candidate subsets $S_j$ corresponding to the interval
$(\lambda_{\mathrm{stop}},+\infty)$ have been explored. Similar to
CSBR, alternative stopping conditions of the form $\stdbars{S_j}\geq
k_{\mathrm{stop}}$ or $\Ec(S_j)\leq\varepsilon_{\mathrm{stop}}$ can be
adopted.

\subsection{Fast implementation}
The $\mathtt{CCV\_Descent}$ procedure calls the function
$\mathtt{intersect}$ to compute the intersection between a concave
polygon $\Sc$ and a line $S_\mathrm{new}$. Lemma~\ref{lem:0} states
that this intersection is empty in two simple situations. Hence, the
call to $\mathtt{intersect}$ is not needed in these situations. This
implementation detail is omitted in Tab.~\ref{tab:ccv} for brevity
reasons.
\begin{lemma}
\label{lem:0}
Let $\Sc=\{S_j,\,j=0,\ldots,J\}$ be a list of supports associated to a
continuous, concave polygon
$\lambda\mapsto\min_{j}\hat{\Jc}(S_j;\lambda)$ with $J+1$ edges,
delimited by $\lambdab=\{\lambda_0,\ldots,\lambda_{J+1}\}$. The
following properties hold for all $j$:
\begin{itemize}
\item If $\delta\Ec_\mathrm{add}(S_j)<\lambda_{j+1}$, then the line
$S_{\mathrm{add}}=S_j+\{\ell_\mathrm{add}\}$ lays above the current concave polygon.
\item If $\delta\Ec_\mathrm{rmv}(S_j)>\lambda_{j}$, then the line
$S_{\mathrm{rmv}}=S_j-\{\ell_\mathrm{rmv}\}$ lays above the current concave polygon.
\end{itemize}
\end{lemma}
\begin{IEEEproof}
  We give a sketch of proof using geometrical arguments. Firstly,
  $\delta\Ec_\mathrm{add}(S_j)$ is the $\lambda$-value of the
  intersection point between lines $S_j$ and
  $S_\mathrm{new}=S_j+\{\ell_\mathrm{add}\}$; see
  Fig.~\ref{fig:refinement}(b). Secondly, we notice that
  $\stdbars{S_{j}}\leq \stdbars{S_{\mathrm{add}}}\leq
  \stdbars{S_{j+1}}$ because the concave polygon is concave and
  $\stdbars{S_{\mathrm{add}}}=\stdbars{S_j}+1$.  It follows from these
  two facts that if $\delta\Ec_\mathrm{add}(S_j)<\lambda_{j+1}$, the
  line $S_{\mathrm{add}}$ lays above $S_{j+1}$ for
  $\lambda\leq\lambda_{j+1}$, and above $S_{j}$ for
  $\lambda\geq\lambda_{j+1}$.

  This proves the first result. A similar sketch applies to the second
  result.
\end{IEEEproof}

\subsection{Main differences between CSBR and $\ell_0$-PD}
\label{sec:algo}
First, we stress that contrary to CSBR, the index $j$ in $\lambda_j$
does not identify with the iteration number anymore for
$\ell_0$-PD. Actually, the current iteration of $\ell_0$-PD is related
to an edge of the concave polygon, \ie a whole interval
$(\lambda_{j+1},\lambda_j)$, whereas the current iteration of CSBR is
dedicated to a single value $\lambda_j$ which is decreasing when the
iteration number $j$ increases.

Second, the computation of the next value $\lambda_{j+1}\leq\lambda_j$
in CSBR is only based on the violation of the lower bound
of~\eqref{eq:stop_SBR}, corresponding to atom selections. In
$\ell_0$-PD, the upper bound is considered as well.  This is the
reason why the $\lambda$-values are not scanned in a decreasing order
anymore. This may improve the very sparse solutions found in the early
iterations within an increased computation time, as we will see
hereafter.
\DRAFT{ 
\begin{figure*}[!t]
\begin{center}
{
  \setlength{\tabcolsep}{.8cm}
  \begin{tabular}{cc}
    \small{(a)}\;\;\;\;\figc[height=40mm]{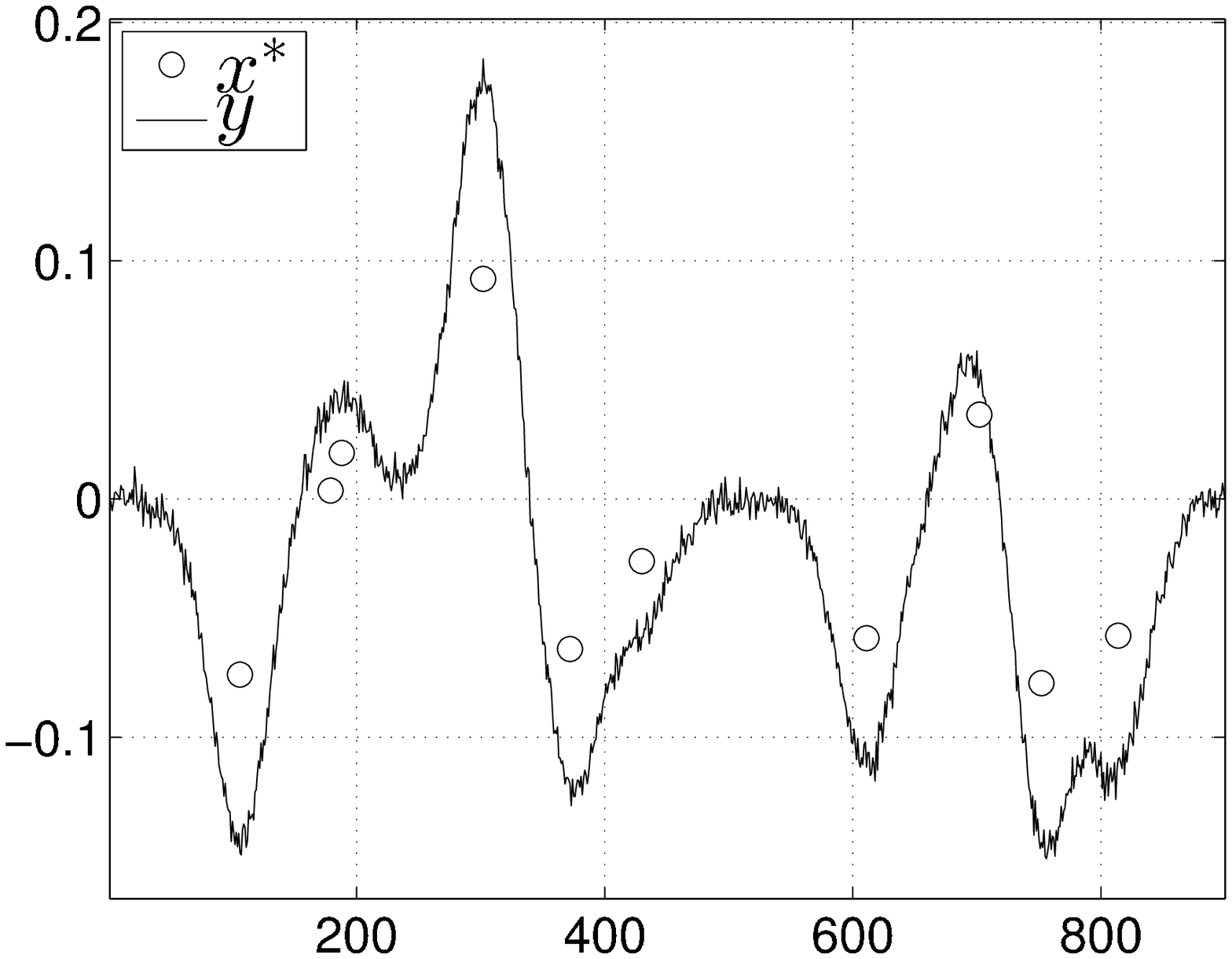}&
    \small{(b)}\;\;\;\;\figc[height=40mm]{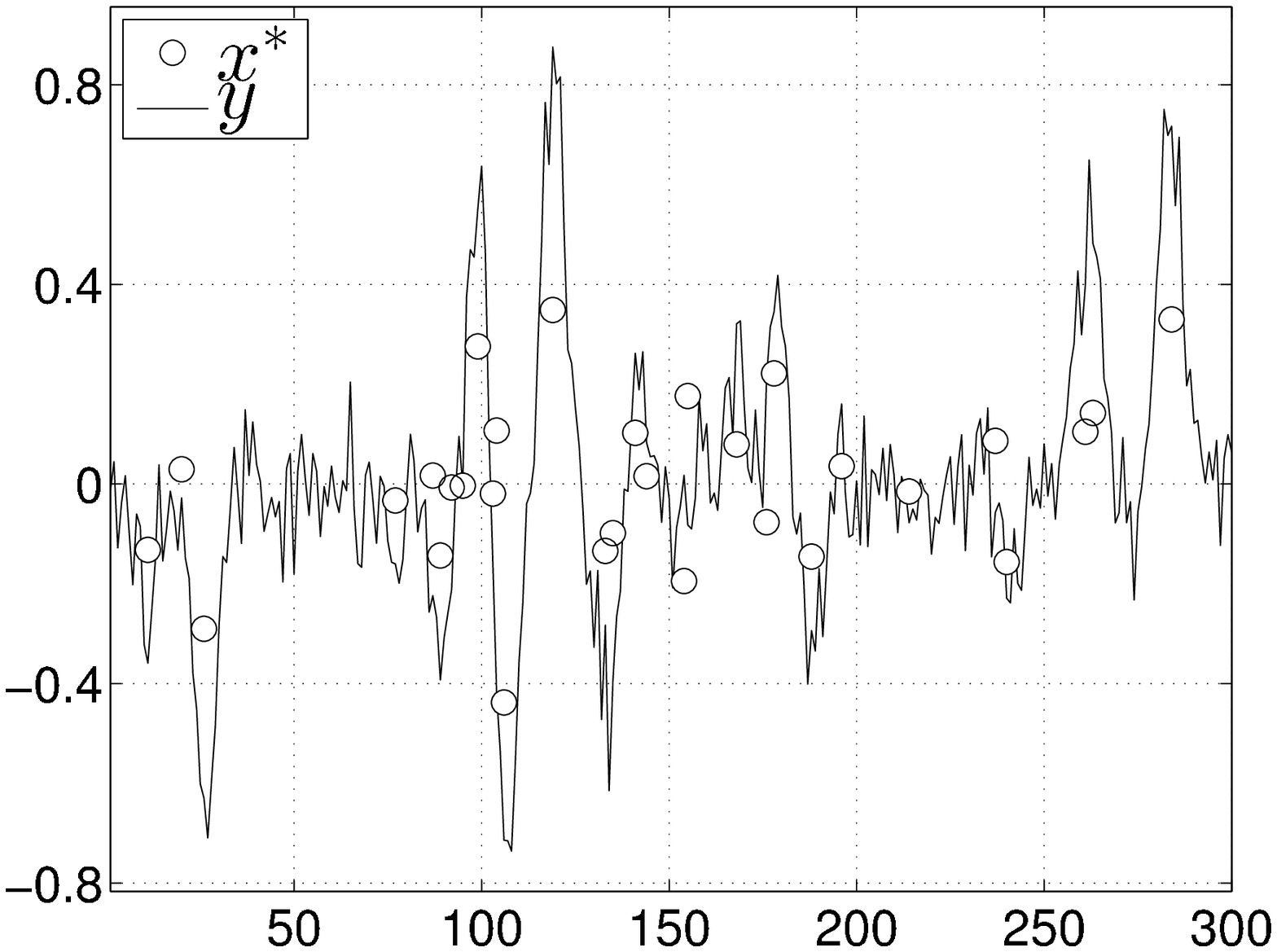}\\[2.1cm]
    \small{(c)}\;\;\;\;\figc[height=40mm]{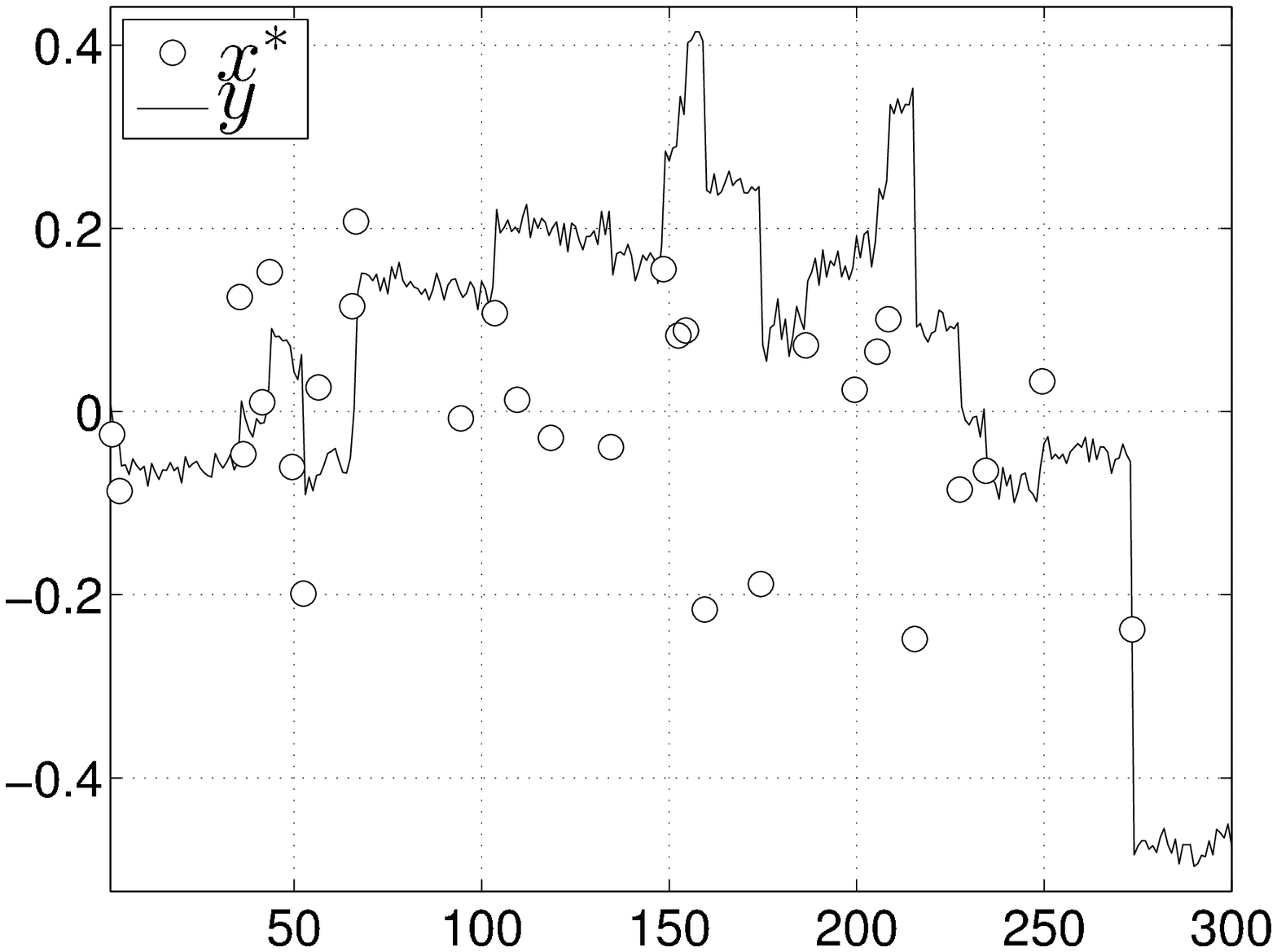}&
    \small{(d)}\;\;\;\;\figc[height=40mm]{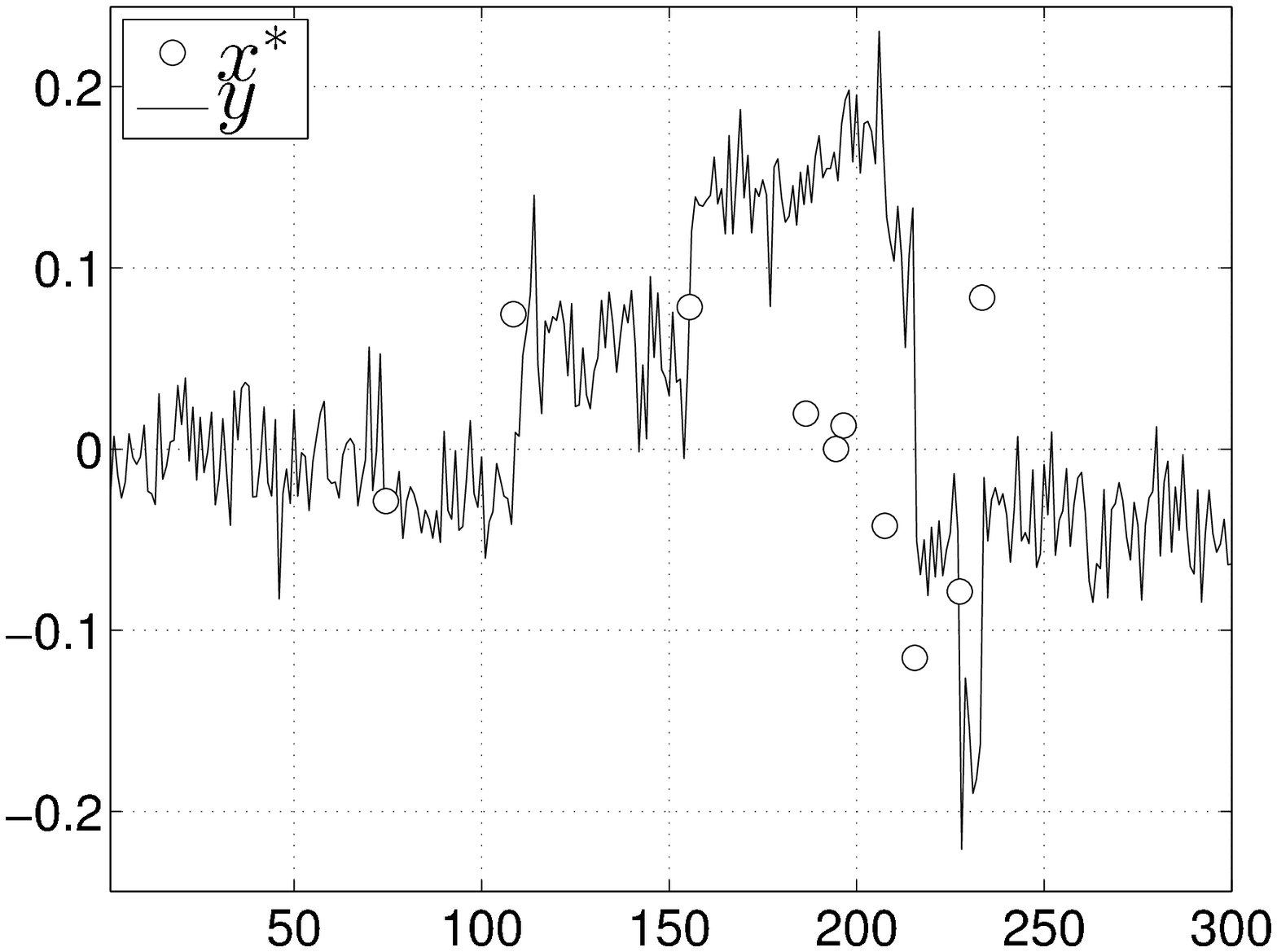}
  \end{tabular}
}
\end{center}
\caption[Generic deconvolution and jump detection problems]{Generic
  deconvolution (a,b) and jump detection (c,d) problems. The data
  vectors \yb and the $k$ nonzero amplitudes of $\xb^\star$ are
  represented in plain lines and with small circles,
  respectively.\quad (a) Sparse deconvolution problem with $k=10$
  spikes, $\mathrm{SNR}=25$ dB, $\sigma=24$ (Gaussian impulse
  response), and $m=900$, $n=756$ (size of dictionary \Ab).\quad (b)
  Sparse deconvolution problem with $k=30$, $\mathrm{SNR}=10$ dB,
  $\sigma=3$, $m=300$, $n=252$. \quad (c) Jump detection problem with
  $k=30$, $\mathrm{SNR}=25$ dB, $m=n=300$.  \quad (d) Jump detection
  problem with $k=10$, $\mathrm{SNR}=10$ dB, $m=n=300$.  }
\label{fig:problems}
\end{figure*}
}

\section{Numerical results}
\label{sec:simuls}
The algorithms are evaluated on two kinds of problems involving
ill-conditioned dictionaries. The behavior of CSBR and $\ell_0$-PD is
first analyzed for simple examples. Then, we provide a detailed
comparison with other nonconvex algorithms for many scenarii.
\begin{figure*}[t]
\centering
{
\setlength{\tabcolsep}{0.01cm}
{
\begin{tabular}{ccc}
\figc[width=58mm]{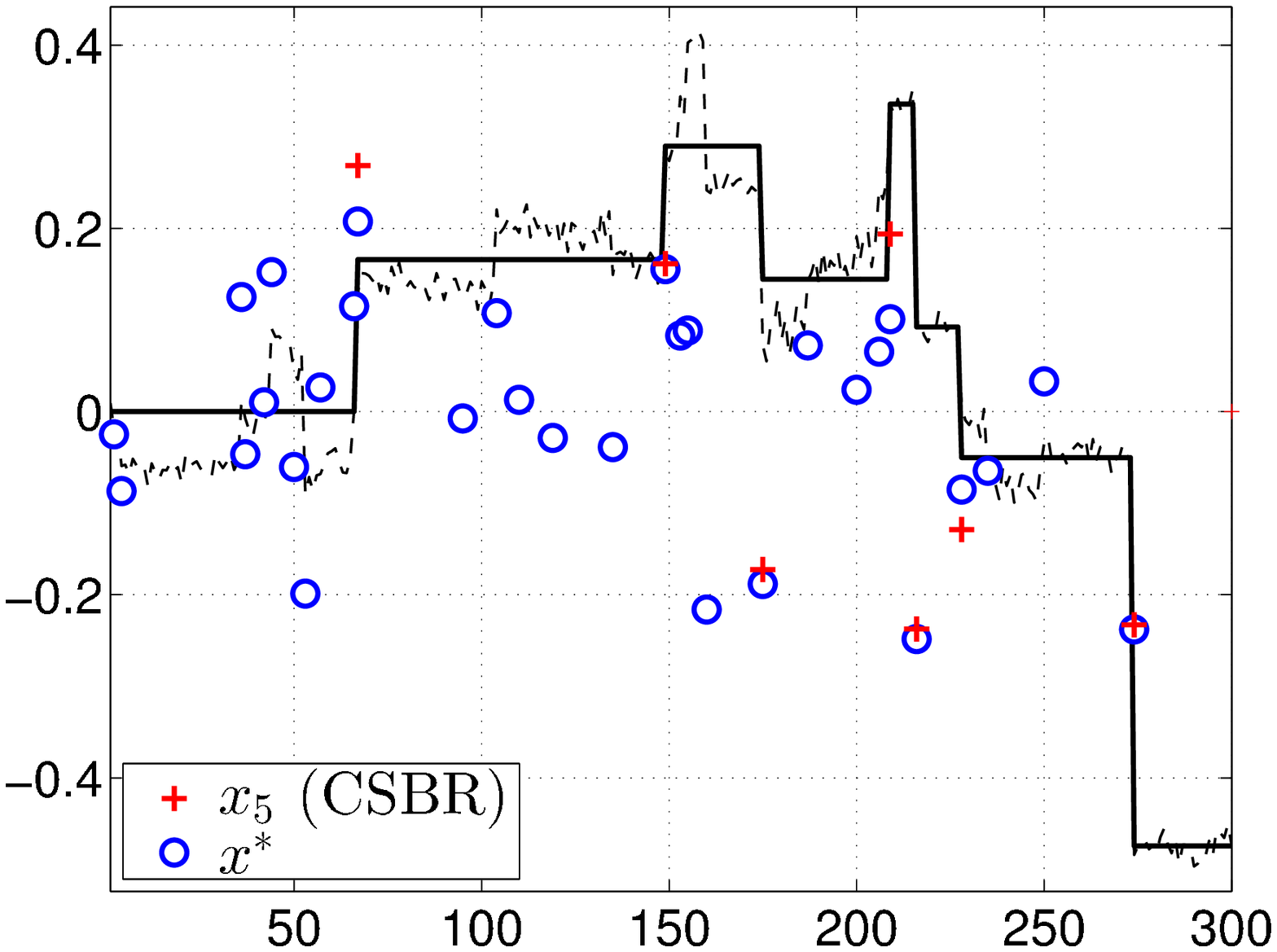}&
\figc[width=58mm]{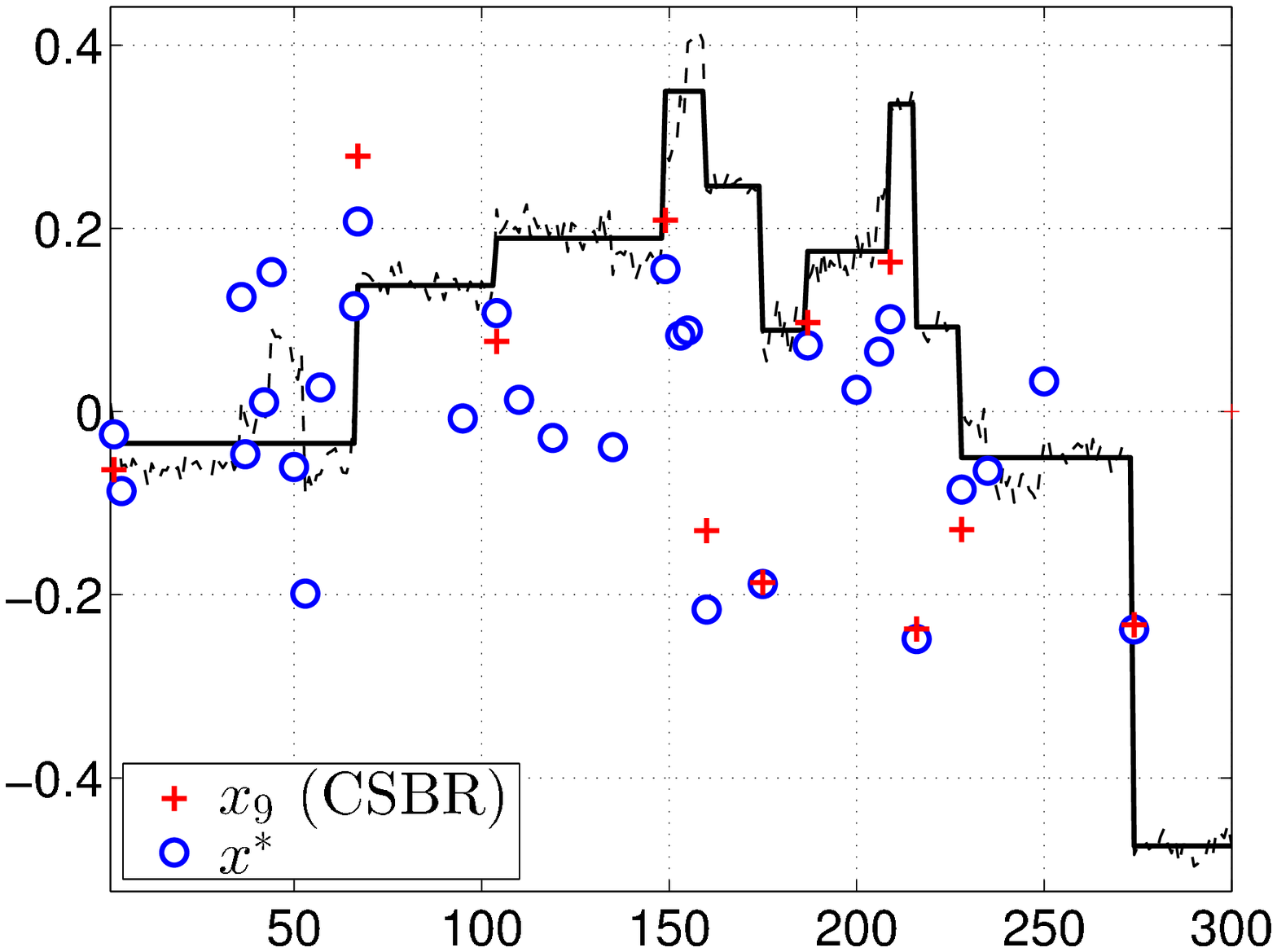}&
\figc[width=58mm]{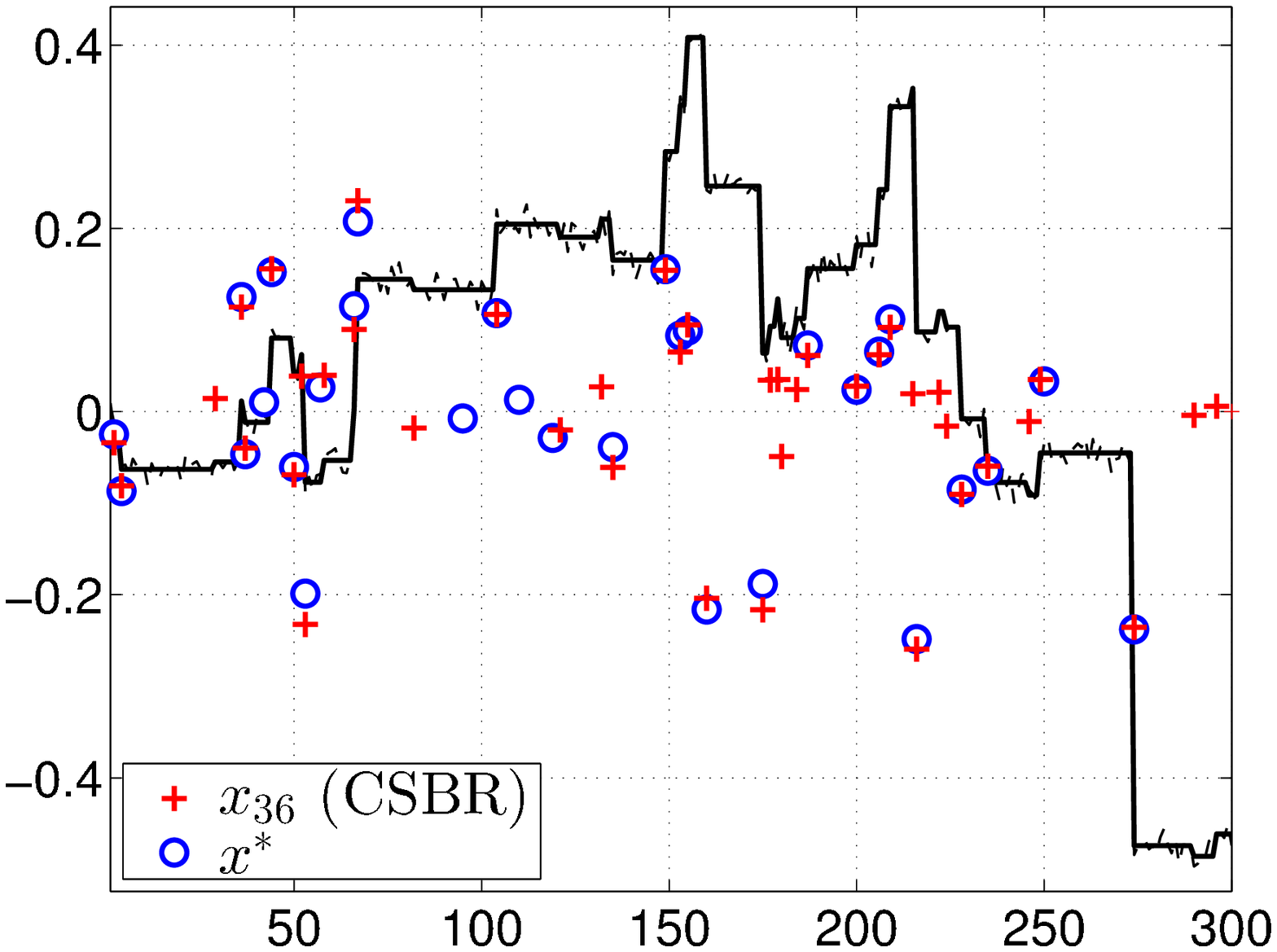}\\[2.2cm]
\small{(a)} $\stdbars{S_5}=7$, $\lambda_5=7.8e^{-2}$ & 
\small{(b)} $\stdbars{S_9}=11$, $\lambda_9=4.1e^{-2}$ & 
\small{(c)} $\stdbars{S_{36}}=41$, $\lambda_{36}=3.4e^{-4}$ 
\end{tabular}
}
}
\caption[Jump detection example]{Jump detection example: processing of
  the data of Fig.~\ref{fig:problems}(c) using CSBR. Three sparse
  solutions $\xb_j$ are shown, each being related to some CSBR output
  $S_j$, with $\|\xb_j\|_0=\stdbars{S_j}$.  The original vector \yb is
  represented in dashed lines and the approximation $\Ab\xb_j$ is in
  solid line.  }
\label{fig:two_examples}
\end{figure*}

\subsection{Two generic problems}
\label{sec:simuls_problems}
The sparse deconvolution problem takes the form
$\yb=\hb\ast\xb^\star+\nb$ where the impulse response \hb is a
Gaussian filter of standard deviation $\sigma$, and the noise \nb is
assumed \iid and Gaussian. The problem rereads $\yb=\Ab\xb^\star+\nb$
where \Ab is a convolution matrix. In the default setting, \yb and \xb
are sampled at the same frequency. \hb is approximated by a finite
impulse response of length $6\sigma$ by thresholding the smallest
values. \Ab is a Toeplitz matrix of dimensions chosen so that any
Gaussian feature $\hb\ast\xb^\star$ if fully contained within the
observation window $\stdacc{1,\ldots,m}$. This implies that \Ab is
slightly undercomplete: $m>n$ with $m\approx n$. Two simulated data
vectors \yb are represented in Fig.~\ref{fig:problems}(a,b) where
$\xb^\star$ are $k$-sparse vectors with $k=10$ and 30, and the
signal-to-noise ratio (SNR) is equal to 25 and 10 dB, respectively. It
is defined by
$\textrm{SNR}=10\,\log\stdpth{\froc{\|\Ab\xb^\star\|_2^2}{(m\sigma_n^2)}}$
where $\sigma_n^2$ is the variance of the noise process \nb.

The jump detection problem is illustrated on
Fig.~\ref{fig:problems}(c,d). Here, \Ab is the squared dictionary
($m=n$) defined by $A_{i,j}=1$ if $i\geq j$, and 0 otherwise. The atom
$\ab_j$ codes for a jump at location $j$, and $x_j^\star$ matches the
height of the jump. When $\xb^\star$ is $k$-sparse, $\Ab\xb^\star$
yields a piecewise constant signal with $k$ pieces, $\xb^\star$ being
the first-order derivative of the signal $\Ab\xb^\star$.

Both generic problems involve either square or slightly undercomplete
dictionaries. The case of overcomplete dictionaries will be discussed
as well, \eg by considering the deconvolution problem with
undersampled observations \yb. The generic problems are already
difficult because neighboring columns of \Ab are highly correlated,
and a number of fast algorithms that are efficient for
well-conditioned dictionaries may fail to recover the support of
$\xb^\star$. The degree of difficulty of the deconvolution problem is
controlled by the width $\sigma$ of the Gaussian impulse response and
the sparsity $k$: for large values of $k$ and/or $\sigma$, the
Gaussian features resulting from the convolution $\hb\ast\xb^\star$
strongly overlap. For the jump detection problem, all the step signals
related to the atoms $\ab_j$ have overlapping supports.

\subsection{Empirical behavior of CSBR and $\ell_0$-PD}
\label{sec:simuls_examples}
\subsubsection{Example}
Consider the problem shown on Fig.~\ref{fig:problems}(c).  Because
CSBR and $\ell_0$-PD provide very similar results, we only show the
CSBR results. CSBR delivers sparse solutions $\xb_j$ for decreasing
$\lambda_j$, $\xb_j$ being the least-square solution supported by the
$j$-th output of CSBR ($S_j$). Three sparse solutions $\xb_j$ are
represented on Fig.~\ref{fig:two_examples}. For the first solution
(lowest value of $\stdbars{S_j}$, largest $\lambda_j$), only the seven
main jumps are being detected (Fig.~\ref{fig:two_examples}(a)). The
cardinality of $S_j$ increases with $j$, and some other jumps are
obtained together with possible false detections
(Figs.~\ref{fig:two_examples}(b,c)).

\subsubsection{Model order selection}
\label{sec:model_sel}
It may often be useful to select a single solution $\xb_j$. The
proposed algorithms are compatible with most classical methods of
model order selection~\cite{Stoica04,Wang04} because they are greedy
algorithms. Assuming that the variance of the observation noise is
unknown, we distinguish two categories of cost functions for
estimation of the order $\|\xb_j\|_0=\stdbars{S_j}$. The first take
the form $\min_{j}\stdacc{m\log\Ec(S_j)+\alpha \stdbars{S_j}}$ where
$\alpha$ equals 2, $\log m$, and $2\log\log m$ for the Akaike, Minimum
Description Length (MDL) and Hannan and Quinn criteria,
respectively~\cite{Stoica04}. The second are cross-validation
criteria~\cite{Wahba77,Golub79}. The sparse approximation framework
allows one to derive simplified expressions of the latter up to the
storage of intermediate solutions of greedy algorithms for
\emph{consecutive} cardinalities~\cite{Miller02,Wang04,Austin10}.
\begin{figure}[t]
\centering
\DRAFT{
  \setlength{\tabcolsep}{0.8cm}
  \begin{tabular}{cc}
    (a)\;\;\;\;\figc[width=58mm]{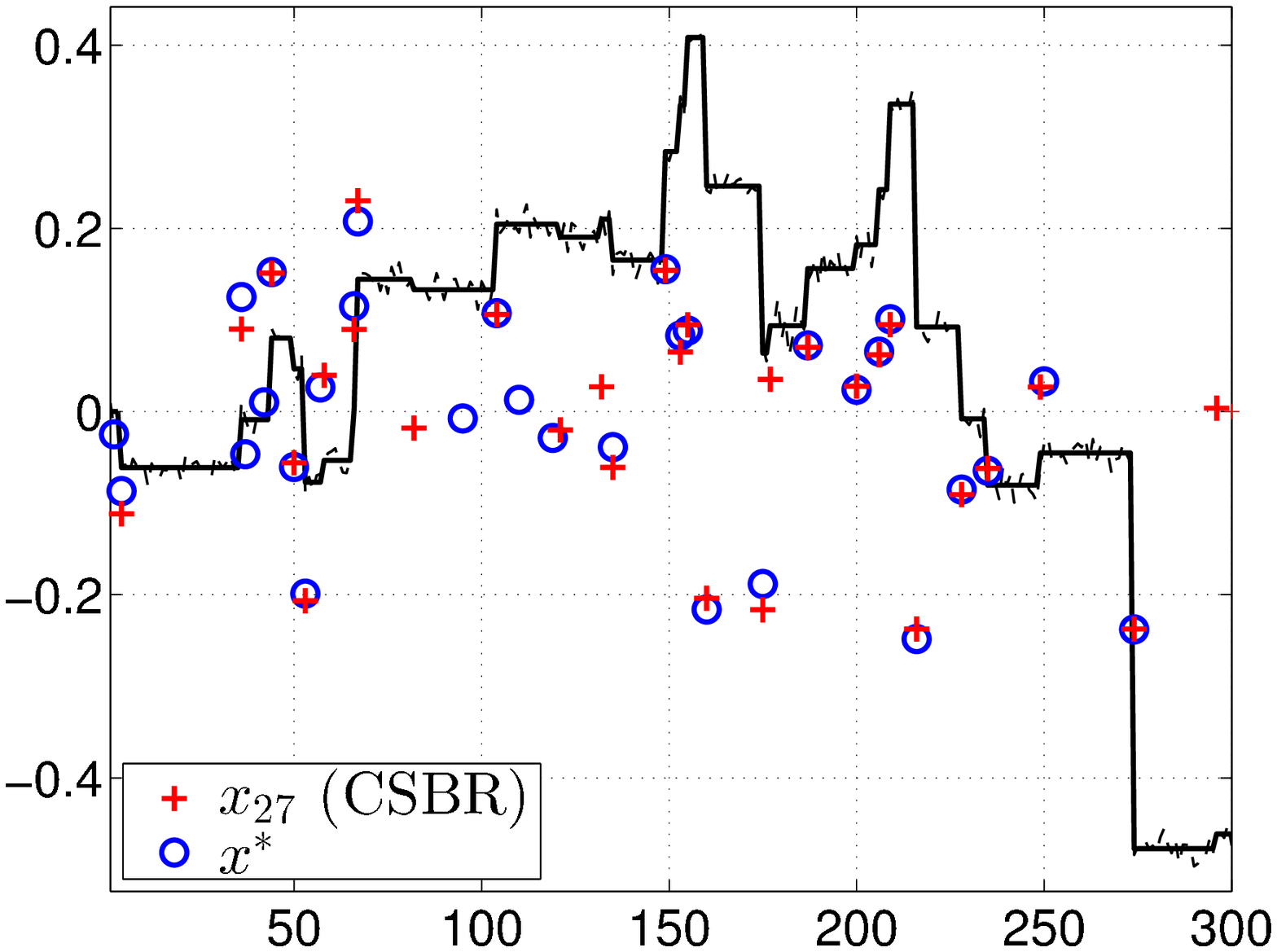}&
    (b)\;\;\;\;\figc[width=58mm]{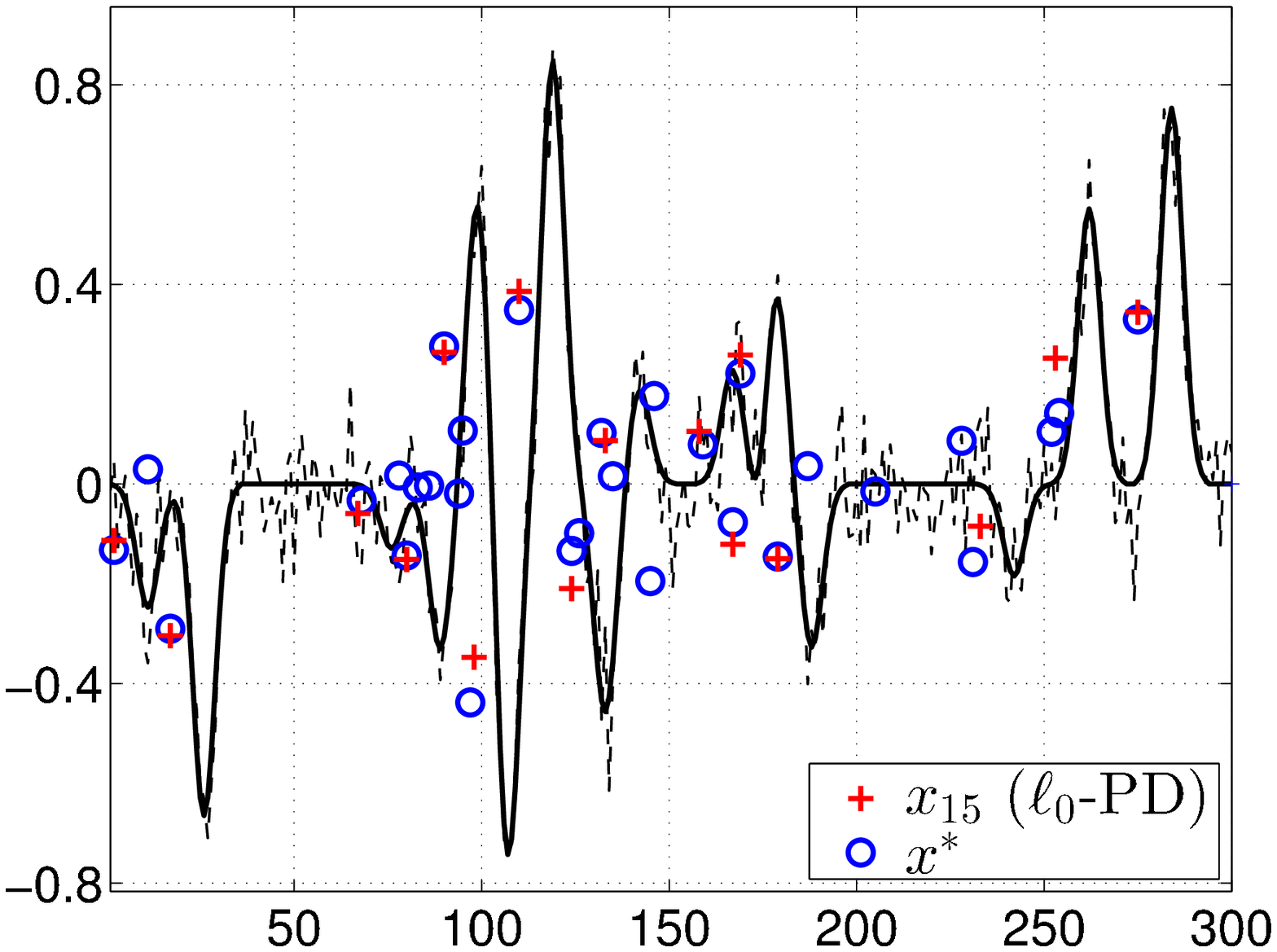}
  \end{tabular}
}
\caption[Model order selection using MDLc]{ Model order selection
  using MDLc: display of the selected sparse solution $\xb_j$ and the
  related data approximation signal. The data of
  Fig.~\ref{fig:problems}(c,b) ($k=30$ true spikes) are processed
  using CSBR and $\ell_0$-PD, respectively. (a) corresponds to the
  simulation shown on Fig.~\ref{fig:two_examples}. The MDLc solution
  is the CSBR output support $S_{27}$ of cardinality 27. (b) is
  related to the $\ell_0$-PD output $S_{15}$, with
  $\stdbars{S_{15}}=16$.  
}
\label{fig:two_examples_mdl}
\end{figure}
\begin{figure}[t]
\centering
{
\DRAFT{\setlength{\tabcolsep}{0.3cm}}
\begin{tabular}{cc}
\DRAFT{
(a)\;\;\;\figc[height=52mm]{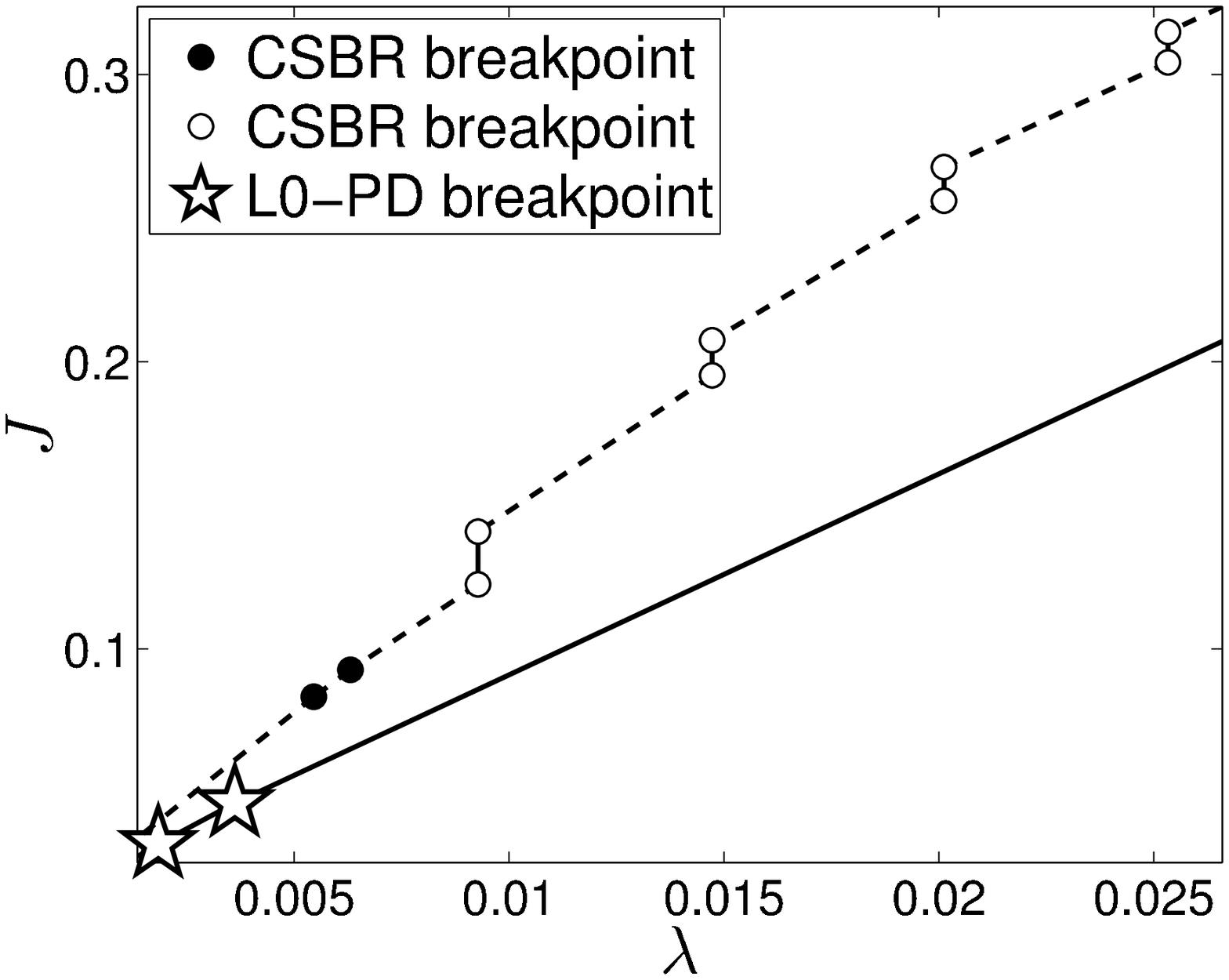}&
(b)\;\;\;\figc[height=52mm]{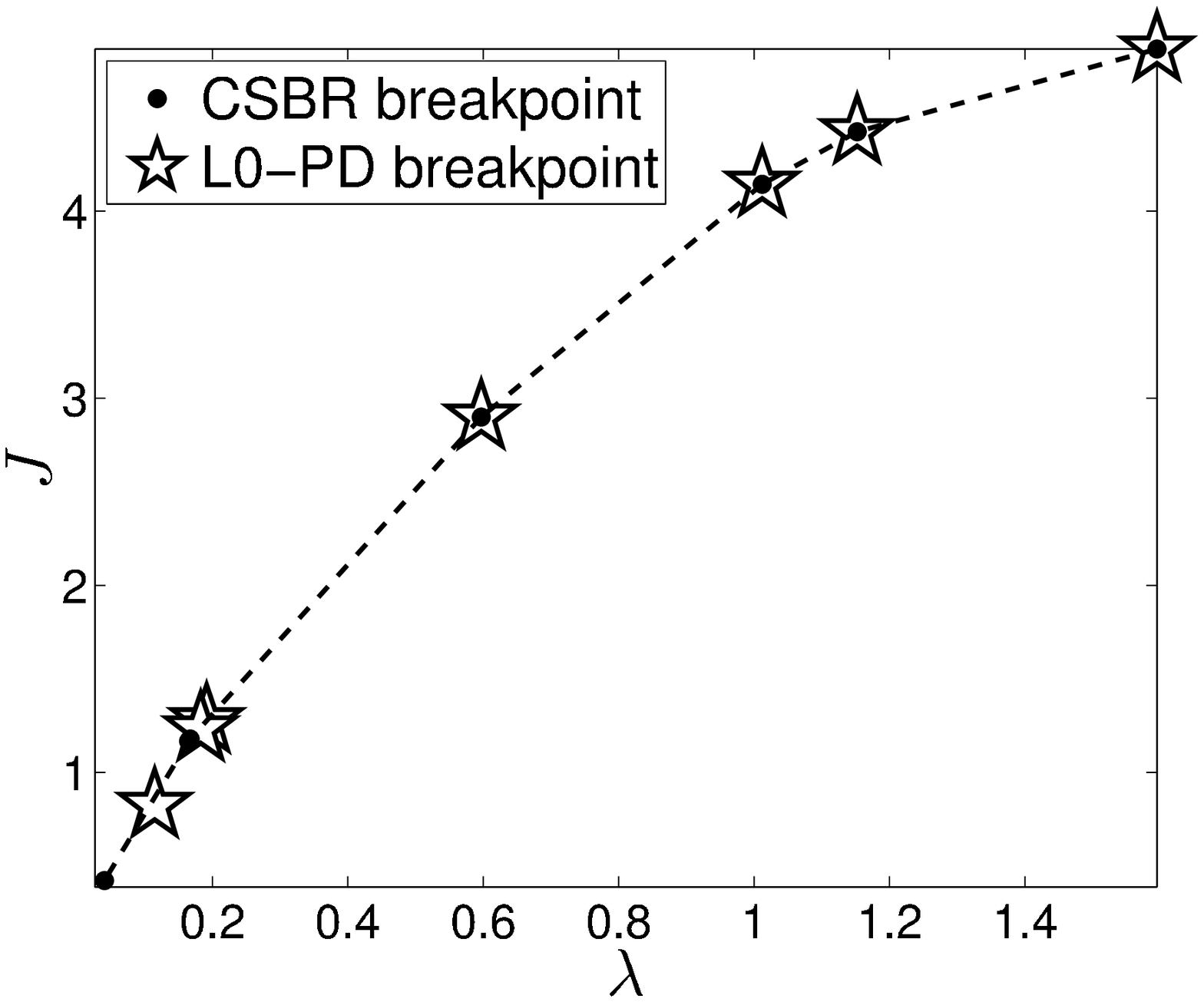}
}
\end{tabular}
}
\caption[Typical approximate $\ell_0$-curves]{Typical approximate 
  $\ell_0$-curves found for the deconvolution problem of 
  Fig.~\ref{fig:problems}(a): zoom in for small and large $\lambda$'s. 
  The $\ell_0$-PD curve is concave and continuous on $\lambda\in\Rbb_+$. 
  The CSBR curve is continuous only for large $\lambda$-values
  (black circles). For low $\lambda$'s, there can be 
  discontinuities at breakpoint locations (white circles). Here, both
  curves almost coincide for large $\lambda$'s. The $\ell_0$-PD curve lays below the CSBR curve for low $\lambda$'s.  
}
\label{fig:simul_rpath}
\end{figure}

For the sparse deconvolution and jump detection problems, we found
that the Akaike and cross validation criteria severely over-estimate
the expected number of spikes. On the contrary, the MDL criterion
yields quite accurate results. We found that the modified MDLc version
dedicated to short data records (\ie when the number of observations
is moderately larger than the model order)~\cite{Ridder05} yields the
best results for all the scenarii we have tested. It reads:
\begin{align}
\min_{j}\Bigacc{\log\Ec(S_j)+
   \frac{\log(m)(\stdbars{S_j}+1)}{m-\stdbars{S_j}-2}}.
\label{eq:mdlc}
\end{align}
Fig.~\ref{fig:two_examples_mdl}(a) illustrates that the number of
spikes found using MDLc is very accurate for high SNRs (27 spikes are
found, the unknown order being 30). It is underestimated for low SNRs:
16 spikes are found (instead of 30) for the simulation of
Fig.~\ref{fig:two_examples_mdl}(b) where $\textrm{SNR}=10$ dB. This
behavior is relevant because for noisy data, the spikes of smallest
amplitudes are drowned in the noise. One cannot expect to detect them.

\subsubsection{Further empirical observations}
Fig.~\ref{fig:simul_rpath} is a typical display of the approximate
$\ell_0$-curves yielded by CSBR and $\ell_0$-PD. The $\ell_0$-PD curve
is structurally continuous and concave whereas for the CSBR curve,
there are two kinds of breakpoints depicted with black and white
circles. The former are ``continuous'' breakpoints. This occurs when
no single replacement is done during the call to SBR
($\textrm{SBR}(S_{\mathrm{init}};\lambda_{j})$ returns
$S_j=S_{\mathrm{init}}$; see Tab.~\ref{tab:CSBR}). Otherwise, a
discontinuity breakpoint (white circle) appears. In
Fig.~\ref{fig:simul_rpath}, the CSBR and $\ell_0$-PD curves almost
coincide for large $\lambda$'s, where only continuous breakpoints can
be observed. For low $\lambda$'s, the $\ell_0$-PD curve lays below the
CSBR curve, and discontinuity breakpoints appear in the latter curve.
\DRAFT{ 
  \begin{figure*}[t]
    \setlength{\tabcolsep}{.3cm}
    \begin{center}
      \begin{tabular}{cc}
        \includegraphics[height=52mm]{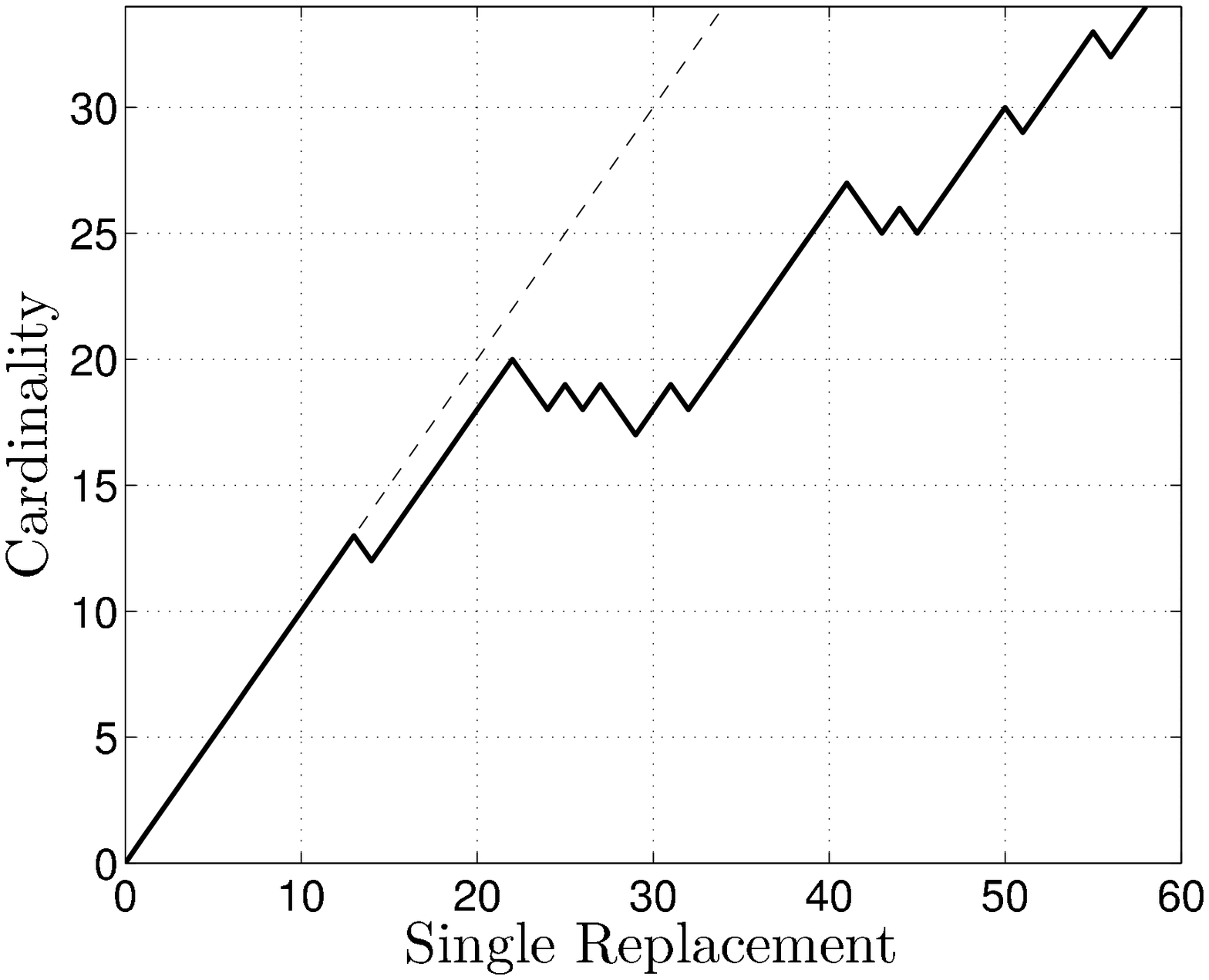}&
        \includegraphics[height=52mm]{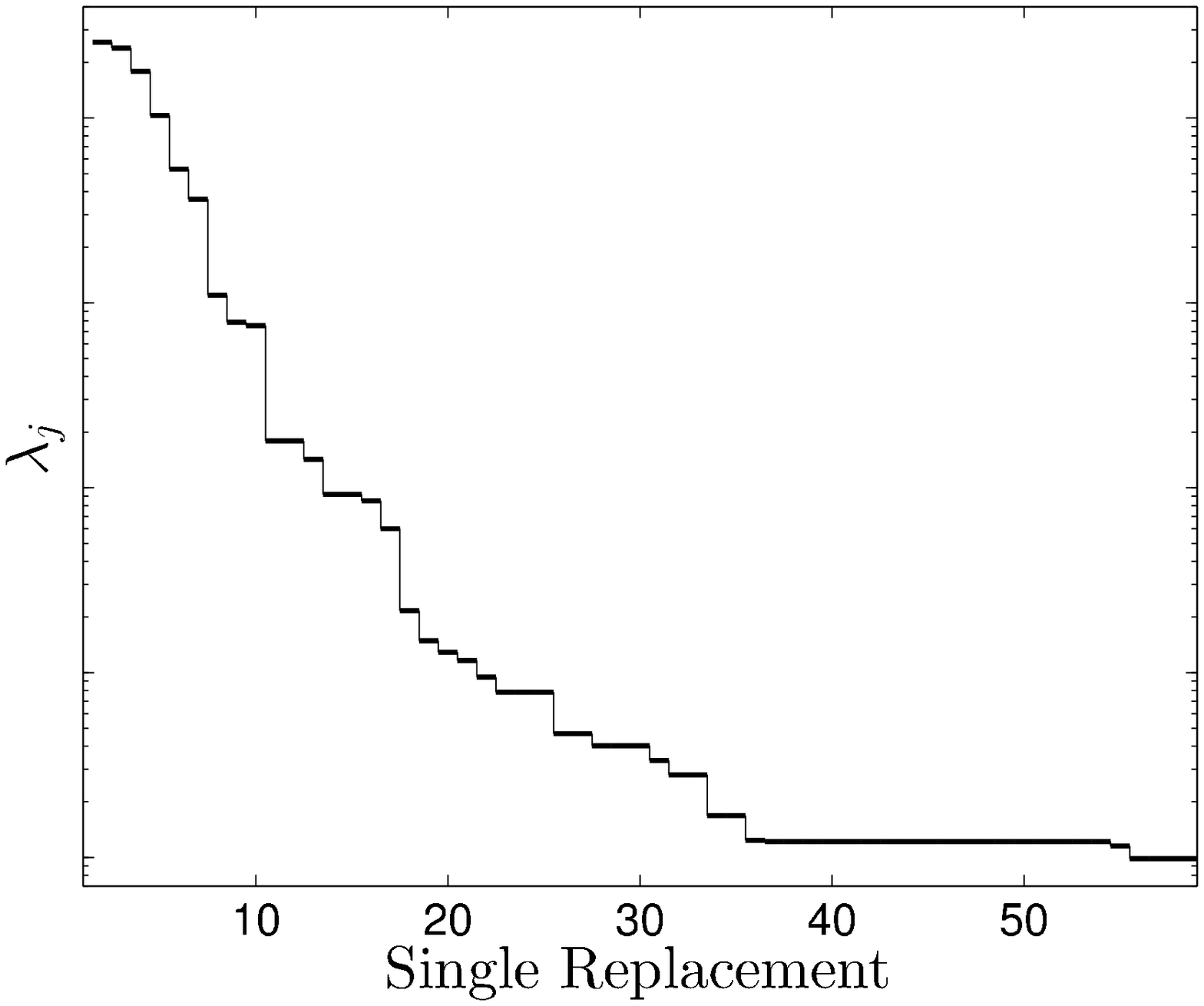}\\
        \small{(a)} CSBR&\small{(b)} CSBR\\      
        \includegraphics[height=52mm]{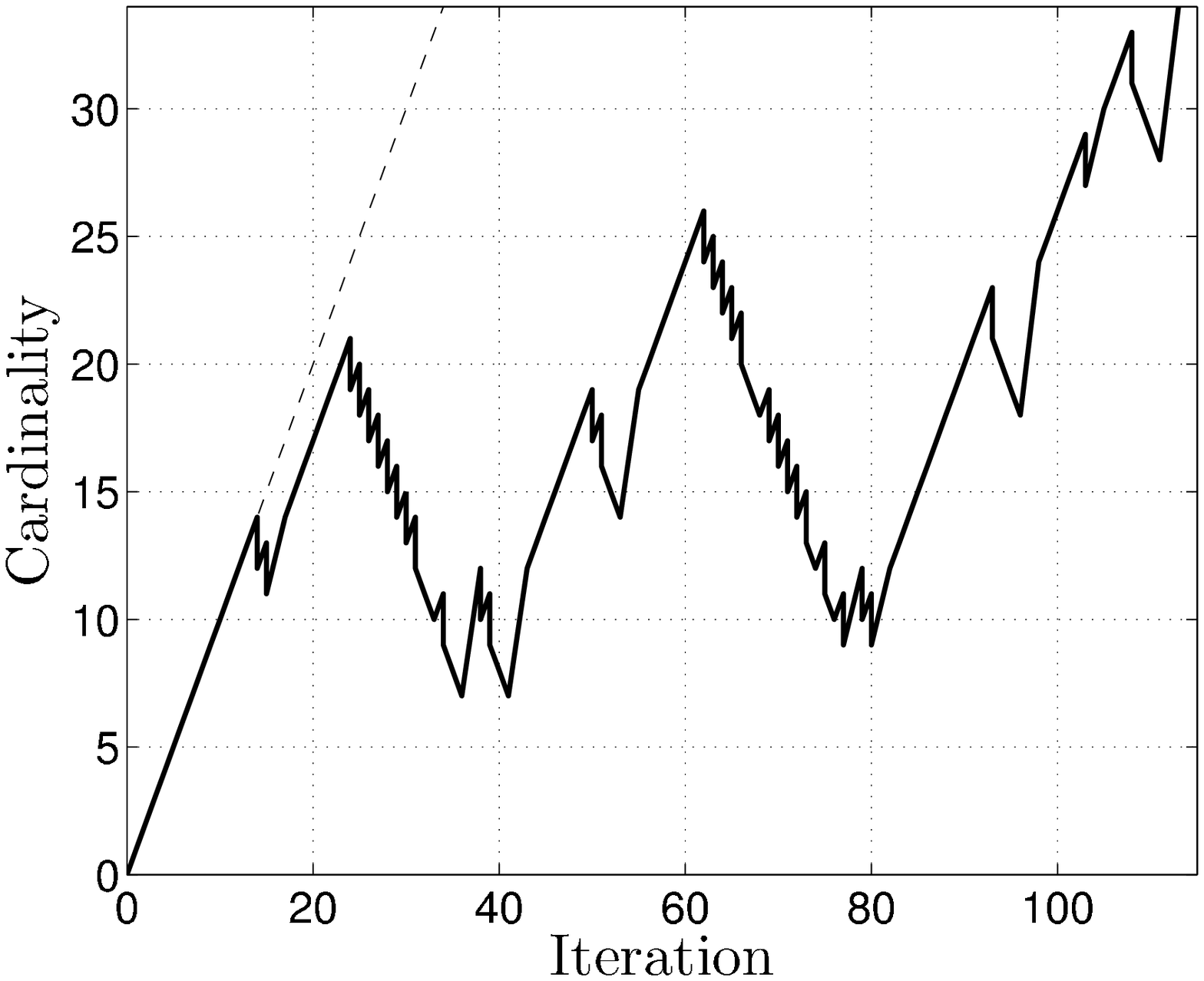}&
        \includegraphics[height=52mm]{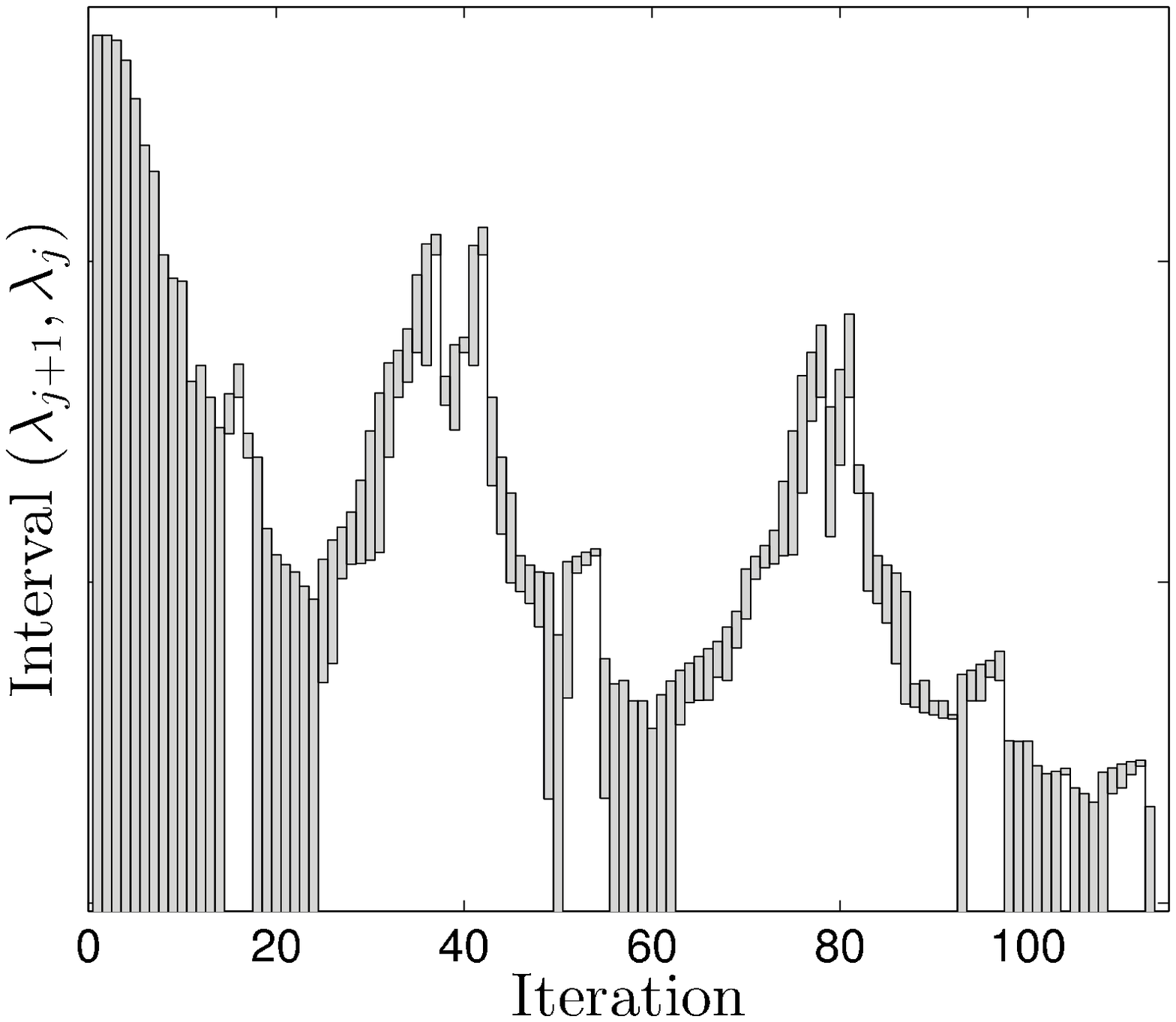}\\
        \small{(c)} $\ell_0$-PD&\small{(d)} $\ell_0$-PD
      \end{tabular}
    \end{center}
    \caption[Series of single replacements performed by CSBR and
    $\ell_0$-PD]{Series of single replacements performed by CSBR and
      $\ell_0$-PD.\quad(a) CSBR: cardinality of the current support
      found after each single replacement.\quad(b) Breakpoints
      $\lambda_j$ found by CSBR, represented in log-scale. SBR is
      executed for each $\lambda_j$, and the number of single
      replacements for fixed $\lambda_j$ matches the length of the
      horizontal steps in the figure.\quad(c) $\ell_0$-PD: cardinality
      of the supports appended to the regularization path during the
      iterations. At each iteration, 0, 1 or 2 supports are
      included. Vertical steps appear whenever two supports are
      simultaneously included.\quad (d)~$\ell_0$-PD: representation in
      log-scale of the current interval $(\lambda_{j+1},\lambda_j)$
      (grey color). When the grey bars reach the bottom of the image,
      the lower bound equals $\lambda_{j+1}=0$.}
  \label{fig:empirical_l0homot}
\end{figure*}
}

Fig.~\ref{fig:empirical_l0homot} provides some insight on the CSBR and
$\ell_0$-PD iterations for a sparse deconvolution problem with
$\|\xb^\star\|_0=17$ and $\textrm{SNR}=20$ dB.  In the CSBR
subfigures, the horizontal axis represents the number of single
replacements: 60 replacements are being performed from the initial
empty support during the successive calls to SBR. For $\ell_0$-PD, the
horizontal axis shows the iteration number. At most two new supports
are being included in the list of candidate subsets at each
iteration. The number of effective single replacements is therefore
increased by 0, 1 or 2. During the first 25 iterations, $\ell_0$-PD
mainly operates atom selections similar to CSBR. The explored subsets
are thus of increasing cardinality and $\lambda$ is decreasing
(Figs.~\ref{fig:empirical_l0homot}(c,d)). From iterations 25 to 40,
the very sparse solutions previously found ($k\leq 20$) are improved
as a series of atom de-selections is performed. They are being
improved again around iteration 80. On the contrary, the sparsest
solutions are never improved with CSBR, which works for decreasing
$\lambda$'s (Figs.~\ref{fig:empirical_l0homot}(a,b)). For $\ell_0$-PD,
the early stopping parameter $\lambda_{\mathrm{stop}}$ may have a
strong influence on the improvement of the sparsest solutions and the
overall computation time. This point will be further discussed below.

\subsection{Extensive comparisons}
\label{sec:simuls_compar}
The proposed algorithms are compared with popular nonconvex algorithms
for both problems introduced in subsection~\ref{sec:simuls_problems}
with various parameter settings: problem dimension ($m,n$), ratio
$m/n$, signal-to-noise ratio, cardinality of $\xb^\star$, and width
$\sigma$ of the Gaussian impulse response for the deconvolution
problem. The settings are listed on Table~\ref{tab:scenarios} for 10
scenarii. Because the proposed algorithms are orthogonal greedy
algorithms, they are better suited to problems in which the level of
sparsity is moderate to high. We therefore restrict ourselves to the
case where $k=\|\xb^\star\|_0\leq 30$.
\begin{table}[t]
  \caption{Settings related to each scenario:
    $k$ is the sparsity. $f$ controls the dictionary size: 
    $m=f\;m_{\mathrm{DEF}},\,n=f\;n_{\mathrm{DEF}}$
    with $n_{\mathrm{DEF}}\approx m_{\mathrm{DEF}}=300$. By default, $f=1$. The 
    undersampling parameter $\Delta$ equals 1 by default
    ($m_{\mathrm{DEF}}\geq n_{\mathrm{DEF}}$). It is increased to generate
    problems with overcomplete dictionaries ($m\approx n/\Delta$). 
    The Gaussian impulse response width is set to 
    $\sigma=f\;\sigma_{\mathrm{DEF}}$ with $\sigma_{\mathrm{DEF}}=3$ 
    or 8.
  }
\label{tab:scenarios}
\setlength{\arraycolsep}{5pt}
\centering
{
\setlength{\tabcolsep}{0.2cm}
\begin{tabular}{|c|c|r|r|r|r|r|r|r|}
\hline
Scenario & Type &SNR & $k$ & $f$ & $\Delta$ & $m$ & $n$ & $\sigma$ \\
\hline
A        & Deconv. & 25 & 30 & 1   &   1      & 300 & 282 & 3 \\ 
\hline
B        & Deconv. & 10 & 10 & 1   &   1      & 300 & 252 & 8 \\ 
\hline
C        & Deconv. & 25 & 10 & 3   &   1      & 900 & 756 & 24 \\  
\hline
D	      & Deconv. & 25 & 30 & 6   &   1      & 1800 & 1692 & 18 \\  
\hline
\hline
E	      & Jumps  & 25 & 10 & 1   &   1      & 300 & 300 & $\emptyset$ \\  
\hline
F	      & Jumps  & 25 & 30 & 1   &   1      & 300 & 300 & $\emptyset$ \\  
\hline
G	      & Jumps  & 10 & 10 & 1   &   1      & 300 & 300 & $\emptyset$ \\  
\hline
\hline
H	      & Deconv.  & $+\infty$ & 10 & 3   &   2      & 450 & 756 & 24 \\  
\hline
I	      & Deconv.  & $+\infty$ & 30 & 3   &   2      & 450 & 756 & 24 \\  
\hline
J	      & Deconv.  & $+\infty$ & 10 & 1   &   4      & 75 & 252 & 8 \\  
\hline
\end{tabular}
}
\end{table}

\subsubsection{Competing algorithms}
We focus on the comparison with algorithms based on nonconvex
penalties. It is indeed increasingly acknowledged that the BPDN
estimates are less accurate than sparse approximation estimates based
on nonconvex penalties. We do not consider forward greedy algorithms
either; we already showed that SBR is (unsurprisingly) more efficient
than the simpler OMP and OLS algorithms~\cite{Soussen11c}. Among the
popular nonconvex algorithms, we consider:
\begin{enumerate}
\item Iterative Reweighted Least Squares (IRLS) for $\ell_q$
  minimization, $q<1$~\cite{Lai13};
\item Iterative Reweighted $\ell_1$ (IR$\ell_1$) coupled with
  the penalty $\log(|x_i|+\varepsilon)$~\cite{Zou06,Candes08,Wipf10};
\item $\ell_0$ penalized least squares for cyclic descent
  (L0LS-CD)~\cite{Seneviratne12};
\item Smoothed $\ell_0$ (SL0)~\cite{Mohimani09,Eftekhari09}.
\end{enumerate}
We resort to a penalized least-square implementation for all
algorithms, the only algorithm directly working with the $\ell_0$
penalty being L0LS-CD. We do not consider simpler thresholding
algorithms (Iterative Hard Thresholding, CoSaMP, Subspace Pursuit)
proposed in the context of compressive sensing since we found that SBR
behaves much better than these algorithms for ill-conditioned
dictionaries~\cite{Soussen11c}. We found that L0LS-CD is more
efficient than thresholding algorithms. Moreover, the cyclic descent
approach is becoming very popular in the recent sparse approximation
literature~\cite{Mazumder11,Marjanovic14a} although its speed of
convergence is sensitive to the quality of the initial solution. Here,
we use the BPDN initial solution
$\argmin_\xb\,\{\|\yb-\Ab\xb\|_2^2+\mu\|\xb\|_1\}$ where $\mu$ is set
to half of the maximum tested $\lambda$-value (more details will be
given hereafter). This simple \emph{ad hoc} setting allows us to get a
rough initial solution that is nonzero and very sparse within a fast
computation time.

The three other considered algorithms work with sparsity measures
depending on an arbitrary parameter. Regarding IRLS, we set $q=0.5$ or
0.1 as suggested in~\cite{Lai13}. We chose to run IRLS twice, with
$q=0.5$ and then $q=0.1$ (with the previous output at $q=0.5$ as
initial solution) so that IRLS is less sensitive to local solutions at
$q=0.1$. SL0 is a GNC-like algorithm working for increasingly
non-convex penalties (\ie Gaussian functions of decreasing
widths). For simplicity reasons, we set the lowest width relative to
the knowledge of the smallest nonzero amplitude of the ground truth
solution $\xb^\star$. The basic SL0 implementation is dedicated to
noise-free problems~\cite{Mohimani09}. There exist several adaptations
in the noisy setting~\cite{Eftekhari09,Ye13} including the precursory
work~\cite{Saito90}. We chose the efficient implementation
of~\cite{Ye13} in which the original pseudo-inverse calculations are
replaced by a quasi-Newton strategy using limited memory BFGS updates.
Finally, the IR$\ell_1$ implementation depends on both the choice of
parameter $\varepsilon$ (which controls the degree of nonconvexity)
and the $\ell_1$ solver. We have tested two $\ell_1$ solvers: the
in-crowd algorithm~\cite{Gill11} together with an empirical setting of
$\varepsilon>0$, and $\ell_1$ homotopy in the limit case
$\varepsilon\rightarrow 0$, following~\cite{Zou06}. We found that
$\ell_1$ homotopy is faster than in-crowd, mainly because the Matlab
implementation of in-crowd (provided by the authors) makes calls to
the \texttt{quadprog} built-in function, which is computationally
expensive for large dimension problems.

\subsubsection{Numerical protocol}
Because the competing algorithms work for a single $\lambda$ value, we
need to define a grid, denoted by
$\{\lambda_i^\mathrm{G},\,i=1,\ldots,N_\lambda\}$, for comparison
purposes. Such grid is defined in logscale for each of the 10 scenarii
$(k,\Ab,\mathrm{SNR})$ defined in Table~\ref{tab:scenarios}. The
number of grid points is $N_\lambda=11$. For a given scenario, $T=30$
trials are being performed in which $k$-sparse vectors $\xb^\star(t)$
and noise vector $\nb(t)$ are randomly drawn. This leads us to
simulate $T$ observation vectors $\yb(t)=\Ab\xb^\star(t)+\nb(t)$ with
$t\in\{1,\ldots,T\}$. Specifically, the location of the nonzero
amplitudes in $\xb^\star(t)$ are uniformly distributed and the
amplitude values are drawn according to an \iid Gaussian
distribution. For each trial $t$, all competing algorithms need to be
run $N_\lambda$ times with $\yb(t)$ and $\lambda_i^\mathrm{G}$ as
inputs whereas CSBR and $\ell_0$-PD are run only once since they
deliver estimates for a continuum of values of $\lambda$. Their
solution for each $\lambda_i^\mathrm{G}$ directly deduces from their
set of output supports and the knowledge of both breakpoints
surrounding $\lambda_i^\mathrm{G}$.
\begin{figure}[!t]
{
\setlength{\tabcolsep}{.0cm}
\DRAFT{
\begin{tabular}{cccc}
$\min_\xb\Jc(\xb;\lambda)$ & CPU Time (seconds)&
$\min_\xb\Jc(\xb;\lambda)$ & CPU Time (seconds)\\
\figc[width=40mm]{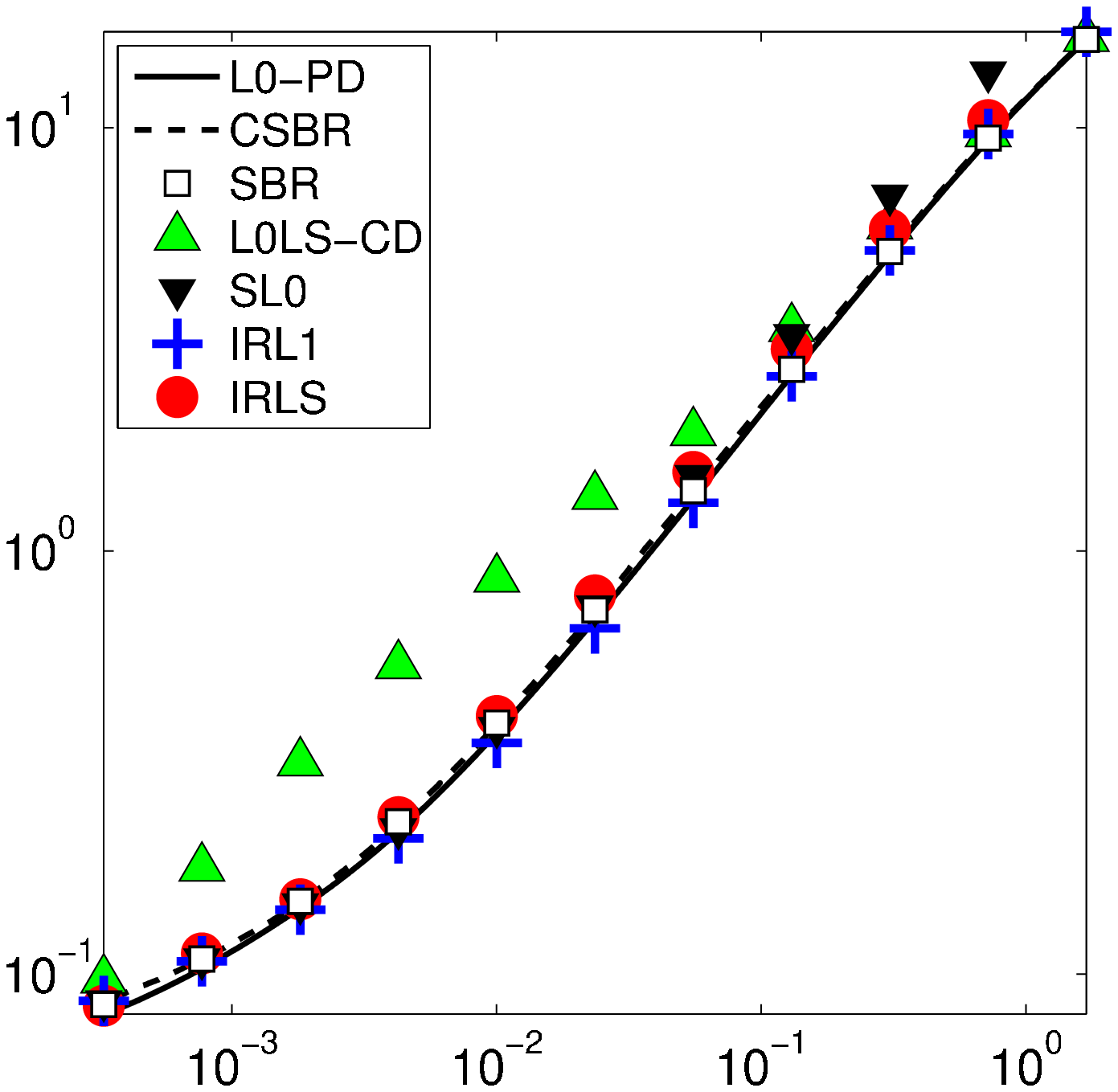}&
\figc[width=40mm]{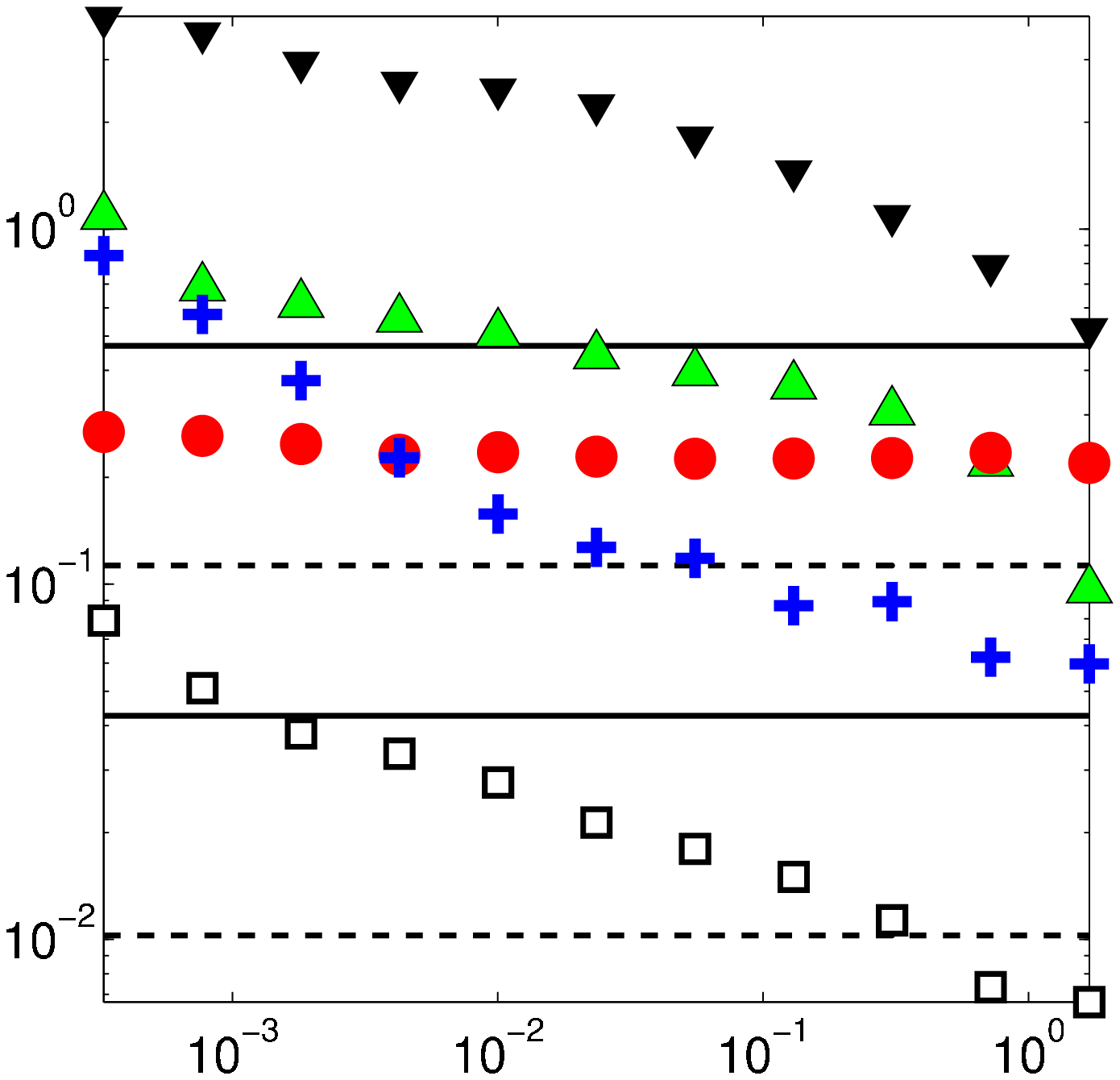}&
\figc[width=40mm]{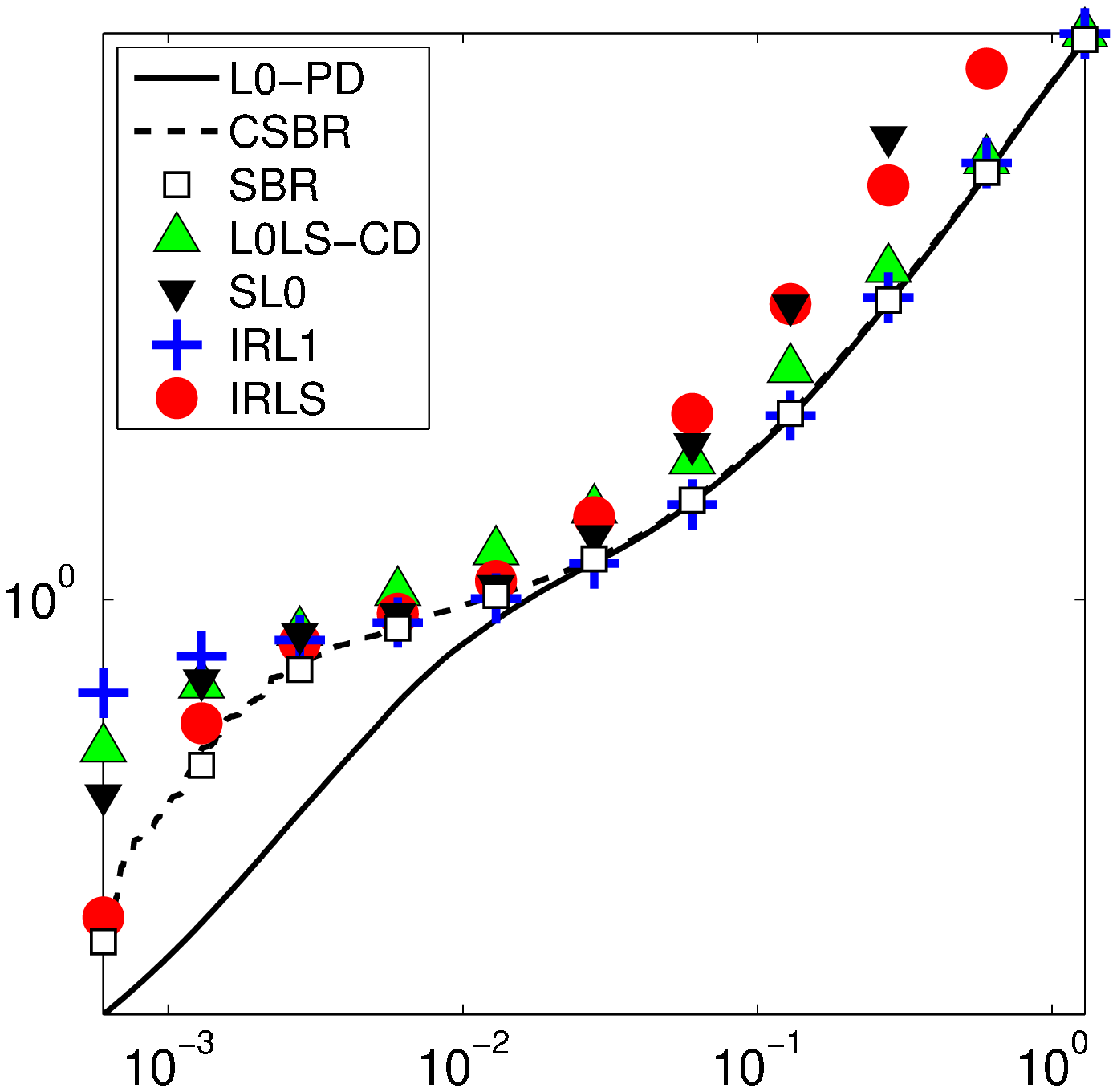}&
\figc[width=40mm]{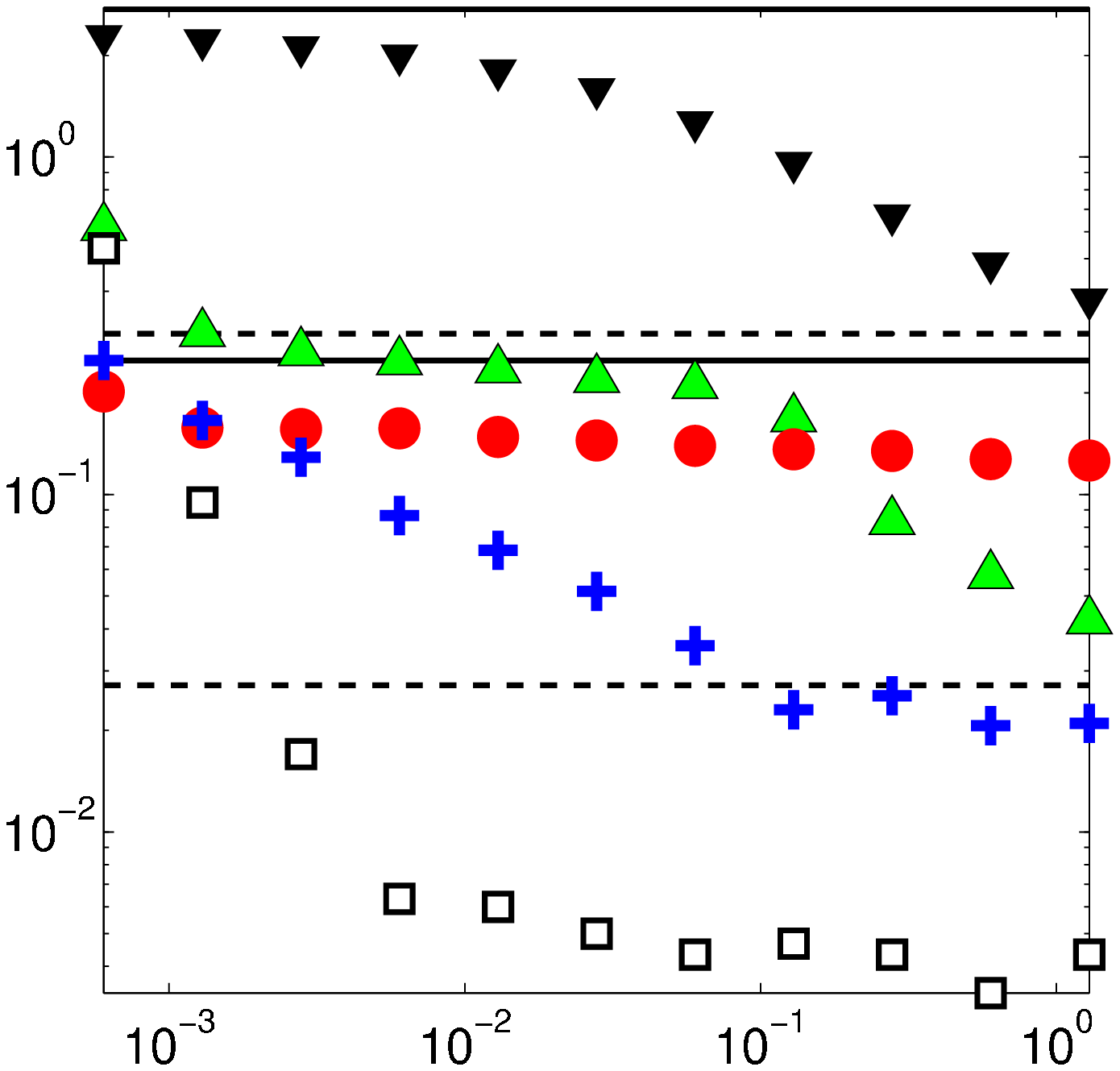}\\[-.2cm]
\begin{tabular}{c}\small{(A)}\end{tabular}&
\begin{tabular}{c}\small{(A)}\end{tabular}&
\begin{tabular}{c}\small{(B)}\end{tabular}&
\begin{tabular}{c}\small{(B)}\end{tabular}\\
\figc[width=40mm]{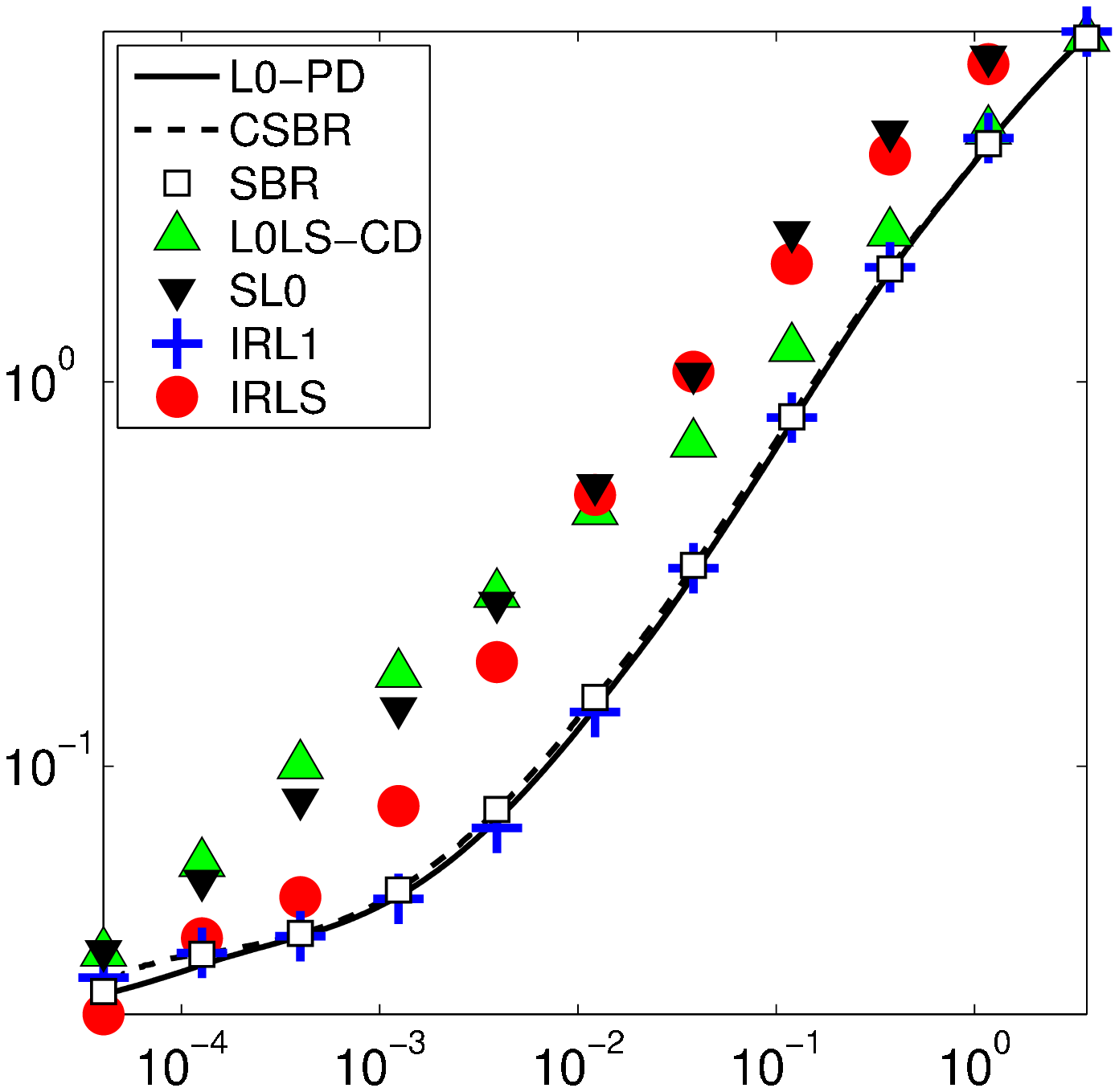}&
\figc[width=40mm]{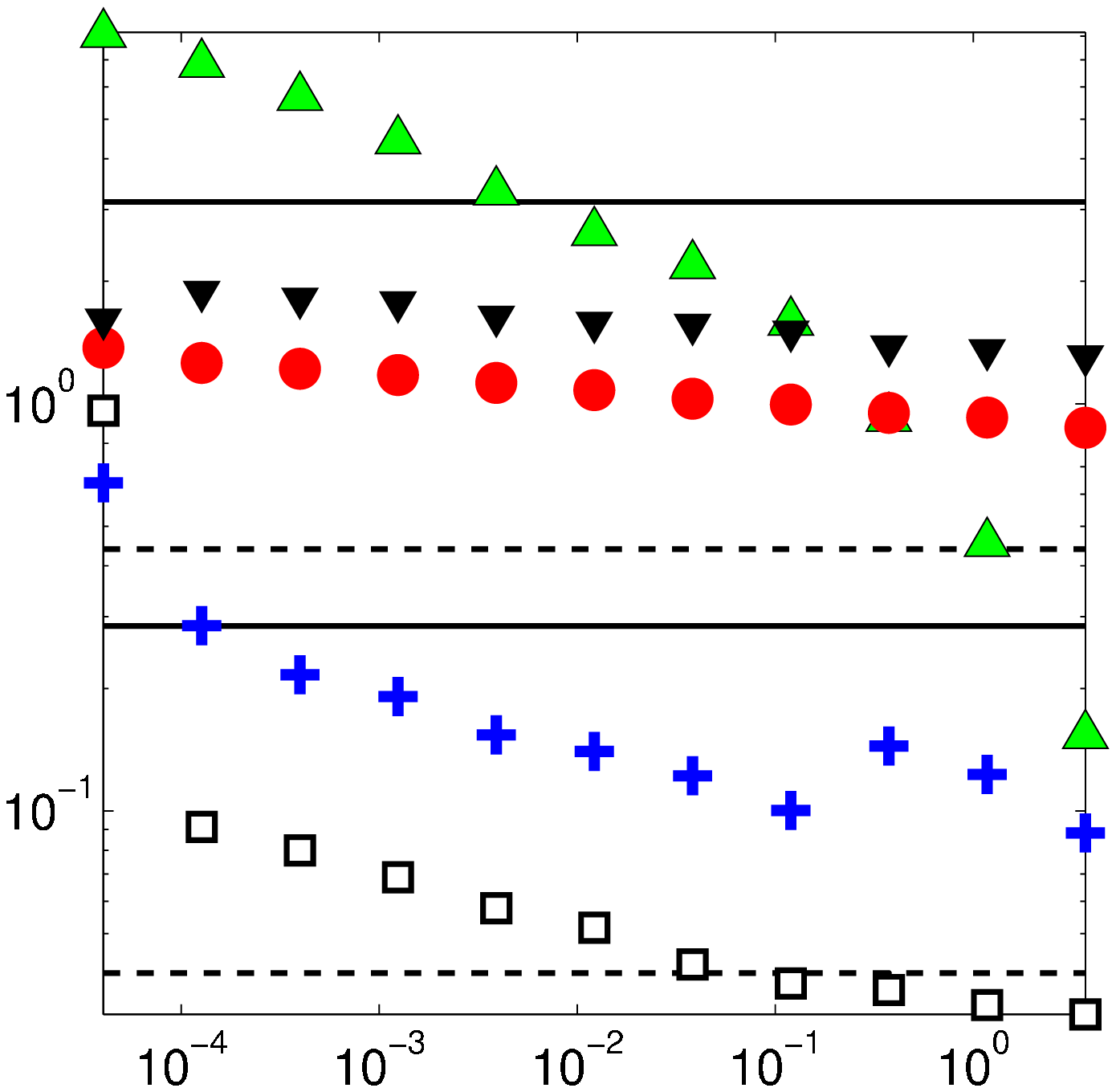}&
\figc[width=41mm]{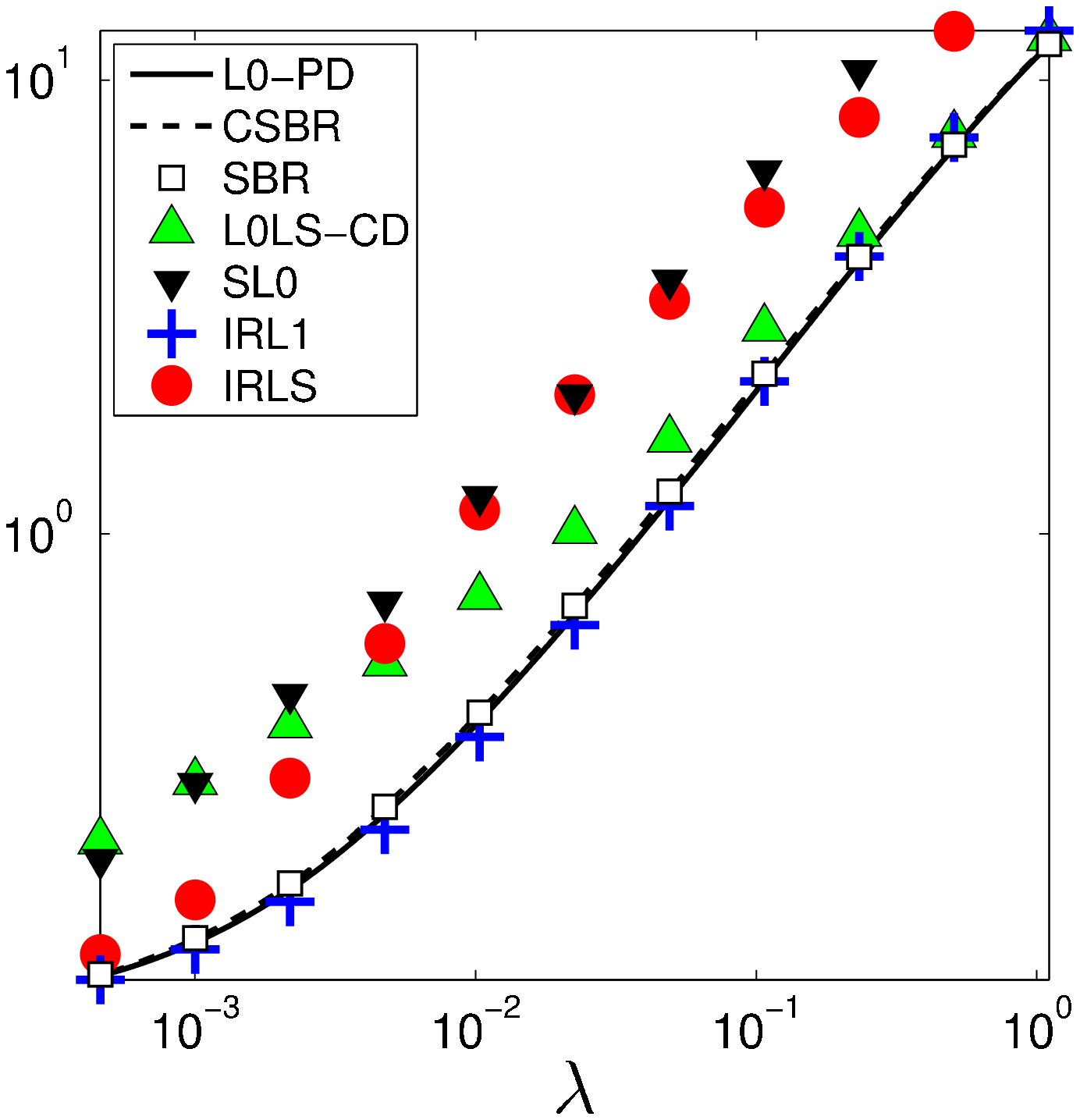}&
\figc[width=41mm]{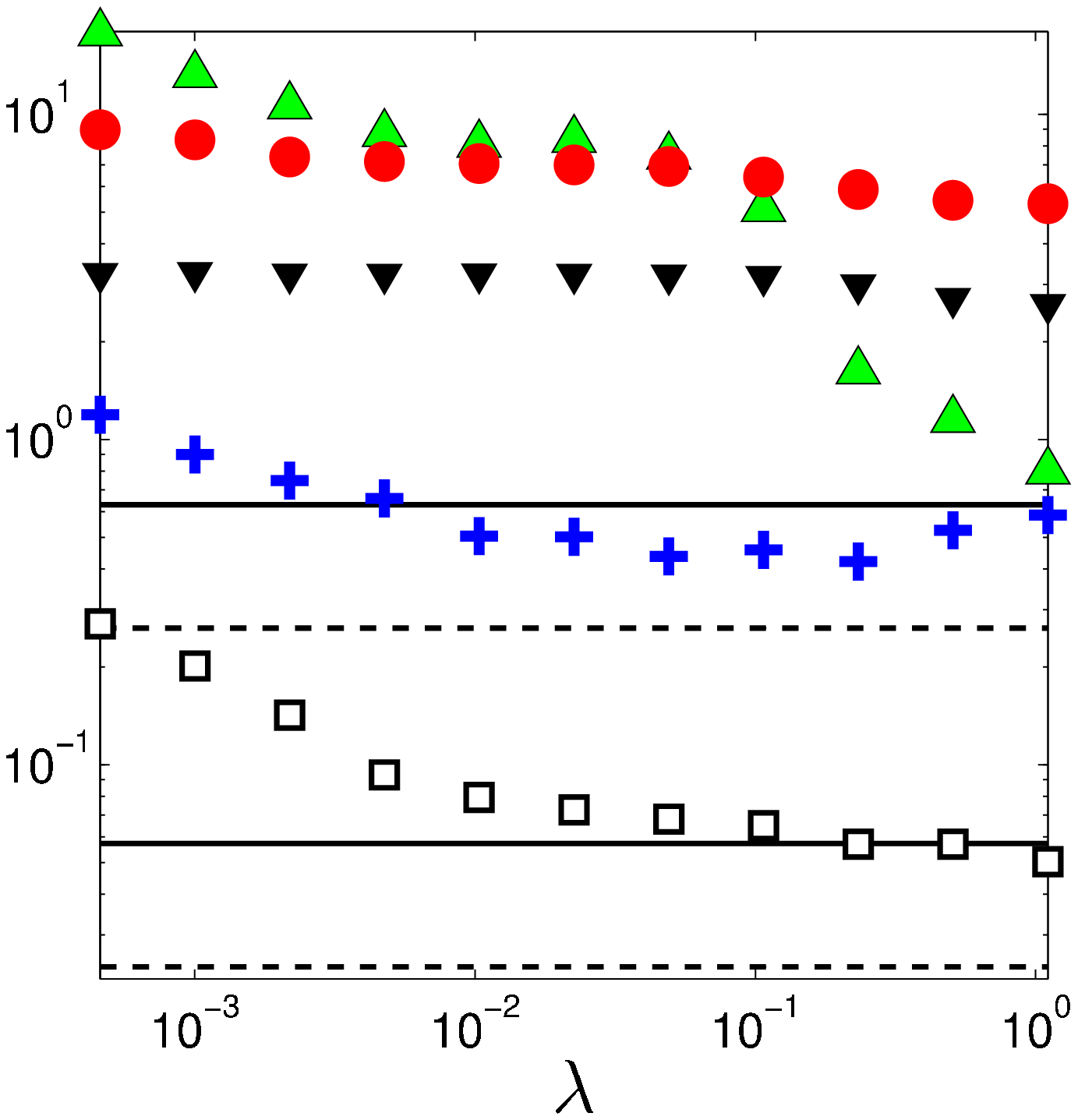}\\[-.2cm]
\begin{tabular}{c}\small{(C)}\end{tabular}&
\begin{tabular}{c}\small{(C)}\end{tabular}&
\begin{tabular}{c}\small{(D)}\end{tabular}&
\begin{tabular}{c}\small{(D)}\end{tabular}
\end{tabular}
}
\FINAL{
\begin{tabular}{ccc}
& $\min_\xb\Jc(\xb;\lambda)$ & CPU Time (seconds)\\
\begin{tabular}{c}\small{(A)}\end{tabular}&
\figc[width=40mm]{S8_Pb2N282fact1SNR25_figJJ_bis}&
\figc[width=40mm]{S8_Pb2N282fact1SNR25_figTT_bis}\\
\begin{tabular}{c}\small{(B)}\end{tabular}&\figc[width=40mm]{S8_Pb2N252fact1SNR10_figJJ_bis}&
\figc[width=40mm]{S8_Pb2N252fact1SNR10_figTT_bis}\\
\begin{tabular}{c}\small{(C)}\end{tabular}&
\figc[width=40mm]{S8_Pb2N756fact3SNR25K10_figJJ_bis}&
\figc[width=40mm]{S8_Pb2N756fact3SNR25K10_figTT_bis}\\
\begin{tabular}{c}\small{(D)}\end{tabular}&
\figc[width=41mm]{S8_Pb2N1692fact6SNR25_figJJ_bis}&
\figc[width=41mm]{S8_Pb2N1692fact6SNR25_figTT_bis}
\end{tabular}
}
}
\caption[Comparison of algorithms for the noisy deconvolution
problem]{Comparison of algorithms for the noisy deconvolution problem,
  \ie for the first scenarii reported on Table~\ref{tab:scenarios}.
  For each scenario, the algorithms are being evaluated in terms of
  $\Jc$-value and of CPU time for $N_\lambda=11$ values
  $\lambda_i^\mathrm{G}$. Evaluations are averaged over 30 trials. The
  overall and mean (normalization by $N_\lambda=11$) CPU times related
  to CSBR (respectively, $\ell_0$-PD) are shown as two parallel
  horizontal lines.  }
  \label{fig:bigcompar_deconv_noisy}
\end{figure}

The algorithms are first evaluated in the optimization viewpoint: the
related criteria are their capacity to reach a low value of
$\Jc(\xb;\lambda)$ and the corresponding CPU time. In this viewpoint,
the proposed methods might be somehow favored since they are more
directly designed with the criterion $\Jc(\xb;\lambda)$ in mind. On
the other hand, $\Jc(\xb;\lambda)$ appears to be a natural indicator
because solving either $\ell_0$-minimization
problem~\eqref{eq:cstr_LS}, \eqref{eq:cstr_LS2} or~\eqref{eq:pen_LS}
is the ultimate goal of any sparse approximation method. As detailed
below, some post-processing will be applied to the outputs of
algorithms that do not rely on the $\ell_0$-norm so that they are not
strongly disadvantaged. Practically, we store the value of
$\Jc(\xb;\lambda_i^\mathrm{G})$ found for each trial and each
$\lambda_i^\mathrm{G}$. Averaging this value over the trials $t$
yields a table $\textrm{TabJ}(a,\lambda_i^\mathrm{G})$ where $a$
denotes a candidate algorithm. Similarly, the CPU time is averaged
over the trials $t$, leading to another table
$\textrm{TabCPU}(a,\lambda_i^\mathrm{G})$. Each table is represented
separately as a 2D plot with a specific color for each algorithm: see,
\eg Fig.~\ref{fig:bigcompar_deconv_noisy}. CSBR and $\ell_0$-PD are
represented with continuous curves because $\Jc(\xb;\lambda)$ is
computed for a continuum of $\lambda$'s, and the CPU time is computed
only once.

The algorithms are also evaluated in terms of support recovery
accuracy. For this purpose, let us first define the ``support error''
as the minimum over $i$ of the distance
\begin{align}
\stdbars{S^\star(t)\backslash S(t,a,\lambda_i^\mathrm{G})}+
\stdbars{S(t,a,\lambda_i^\mathrm{G})\backslash S^\star(t)}
\label{eq:SE}
\end{align}
between the support $S^\star(t)$ of the unknown sparse vector
$\xb^\star(t)$ and the support $S(t,a,\lambda_i^\mathrm{G})$ of the
sparse reconstruction at $\lambda_i^\mathrm{G}$ with algorithm
$a$. \eqref{eq:SE} takes into account both numbers of false negatives
$\stdbars{S^\star(t)\backslash S(t,a,\lambda_i^\mathrm{G})}$ and of
false positives $\stdbars{S(t,a,\lambda_i^\mathrm{G})\backslash
  S^\star(t)}$. Denoting by
$S(t,a,\lambda_{\mathrm{opt}}^\mathrm{G})\leftarrow
S(t,a,\lambda_i^\mathrm{G})$ the solution support that is the closest
to $S^\star(t)$ according to~\eqref{eq:SE}, we further consider the
number of true positives in
$S(t,a,\lambda_{\mathrm{opt}}^\mathrm{G})$, defined as
$\stdbars{S^\star(t)\cap S(t,a,\lambda_{\mathrm{opt}}^\mathrm{G})}$.
We will thus report:
\begin{itemize}
\item the support error; 
\item the corresponding number of true positives;
\item the corresponding model order $\stdbars{S(t,a,\lambda_{\mathrm{opt}}^\mathrm{G})}$.
\end{itemize}
Averaging these measures over $T$ trials yields the support error
score SE($a$), the true positive score TP($a$) and the model order,
denoted by Order($a$). The numbers of false positives (FP) and of
true/false negatives can be directly deduced, \eg FP($a$) $=$
Order($a$) $-$ TP($a$).

The underlying idea in this analysis is that when SE is small
(respectively, TP is high), the algorithms are likely to perform well
provided that $\lambda$ is appropriately chosen.  However, in
practical applications, only one estimate is selected using a suitable
model selection criterion. We therefore provide additional evaluations
of the MDLc estimate accuracy. For CSBR and $\ell_0$-PD, all output
supports are considered to compute the MDLc estimate as described in
subsection~\ref{sec:simuls_examples}. For other algorithms, it is
equal to one of the sparse reconstructions obtained at
$\lambda_i^\mathrm{G}$ for $i\in\{1,\ldots,N_\lambda\}$. The same
three measures as above are computed for the MDLc estimate and
averaged over $T$ trials. They are denoted by MDLc-SE($a$),
MDLc-TP($a$) and MDLc-Order($a$).

\subsubsection{Technical adaptations for comparison purposes}
Because IRLS and SL0 do not deliver sparse vectors in the strict
sense, it is necessary to sparsify their outputs before computing
their SE($a$) score. This is done by running one iteration of cyclic
descent (L0LS-CD): most small nonzero amplitudes are then thresholded
to 0. Regarding the values of $\Jc(\xb;\lambda)$, a post-processing is
performed for algorithms that do not rely on the $\ell_0$-norm. This
post-processing can be interpreted as a local descent of
$\Jc(\xb;\lambda)$. It consists in: \emph{(i)} running one iteration
of cyclic descent (L0LS-CD); \emph{(ii)} computing the squared error
related to the output support.  L0LS-CD is indeed a local descent
algorithm dedicated to $\Jc(\xb;\lambda)$ but the convergence towards
a least-square minimizer is not reached in one iteration.
\begin{figure}[t]
{
\centering
\DRAFT{\setlength{\tabcolsep}{.0cm}
\begin{tabular}{cccc}
$\min_\xb\Jc(\xb;\lambda)$ & CPU Time (seconds)&
$\min_\xb\Jc(\xb;\lambda)$ & CPU Time (seconds)\\
\figc[width=40mm]{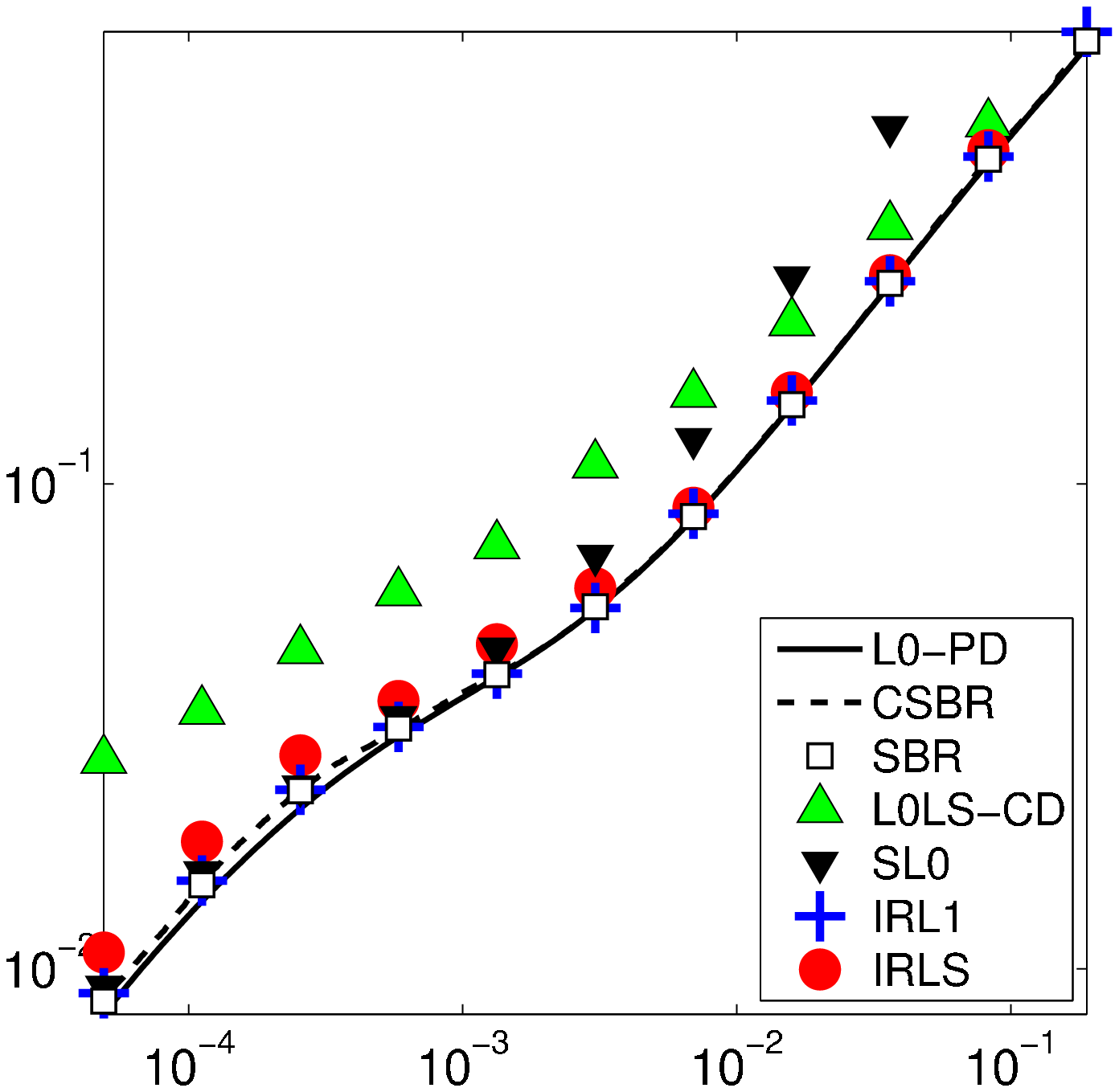}&
\figc[width=40mm]{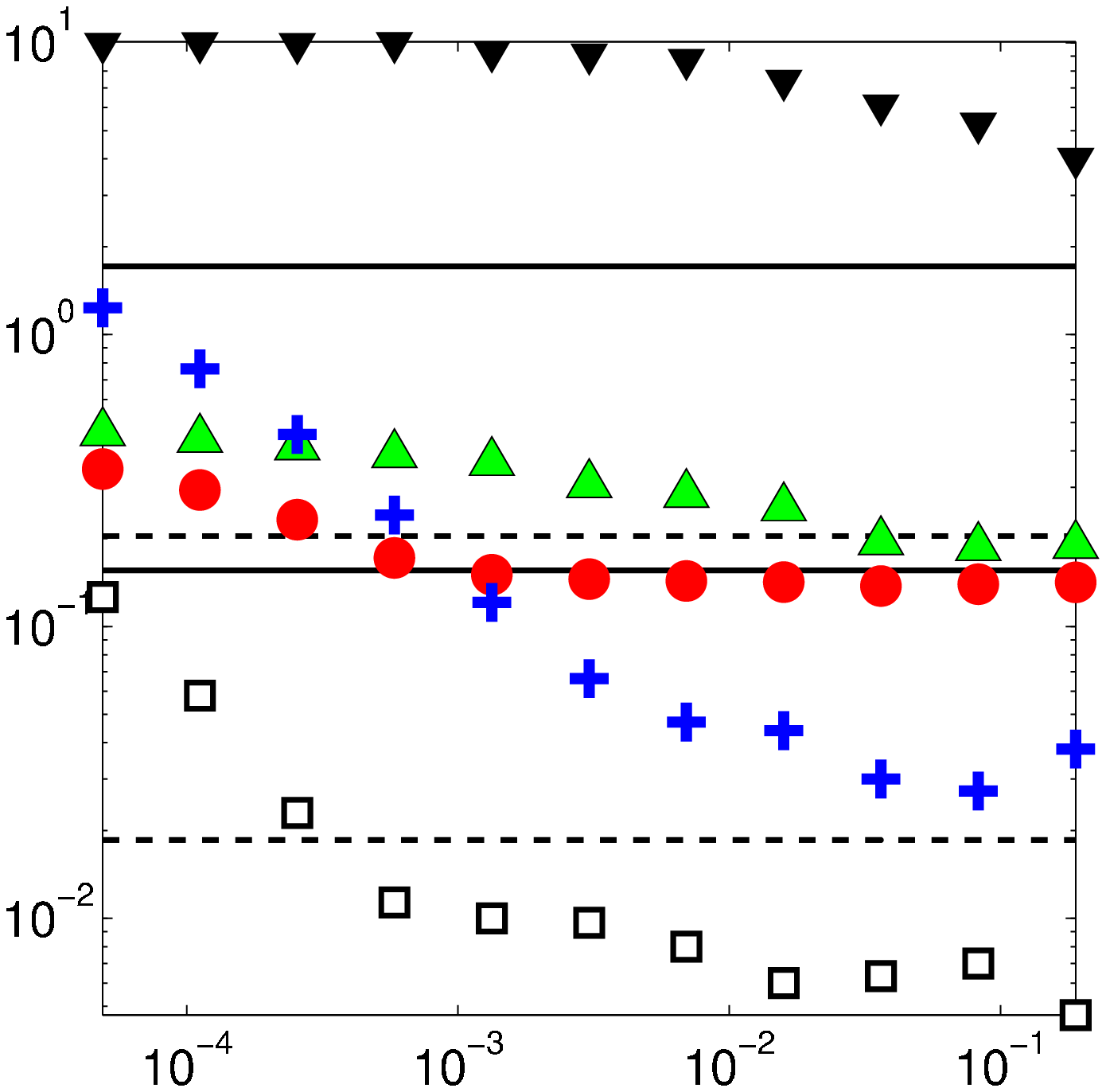}&
\figc[width=40mm]{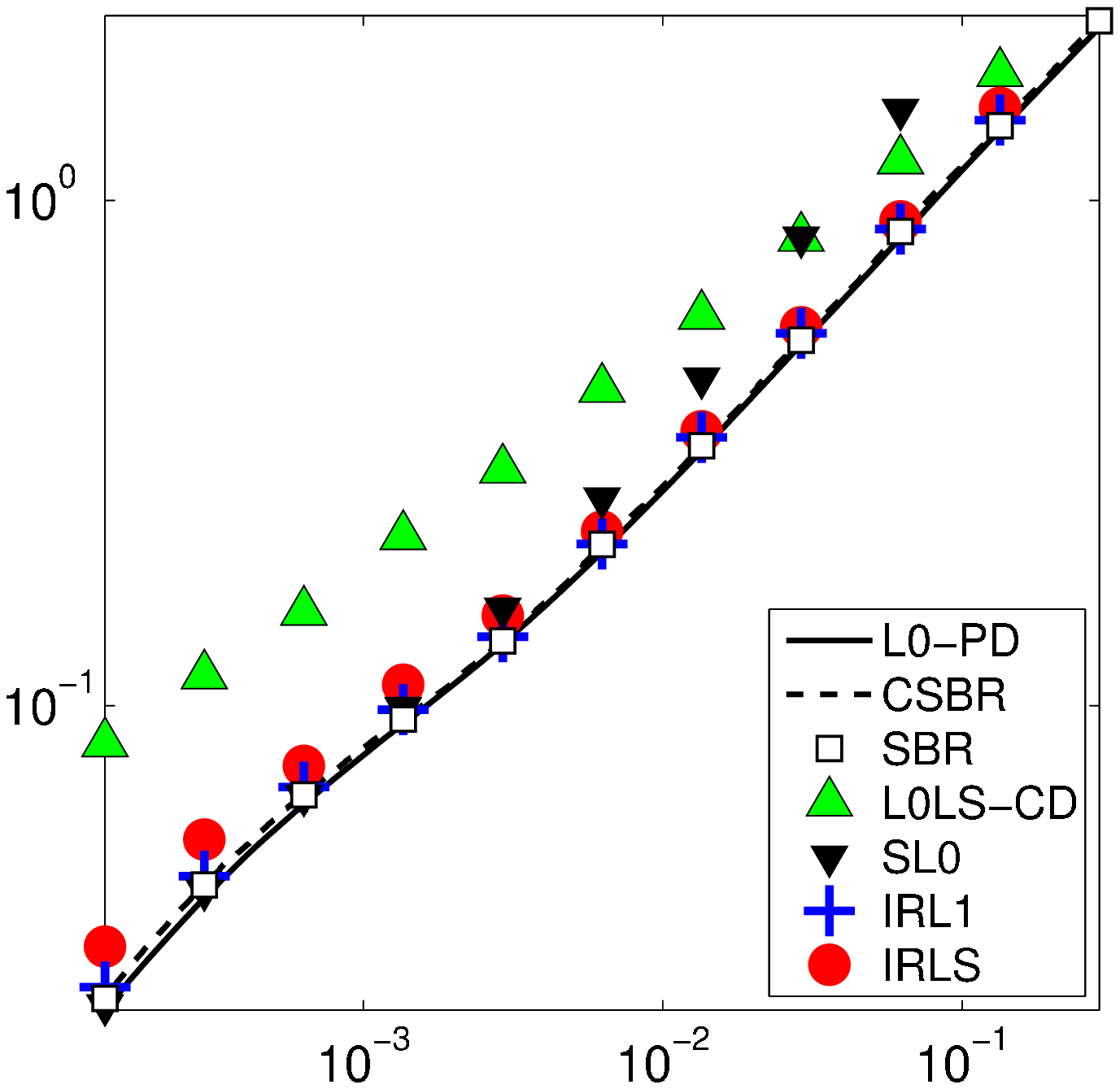}&
\figc[width=40mm]{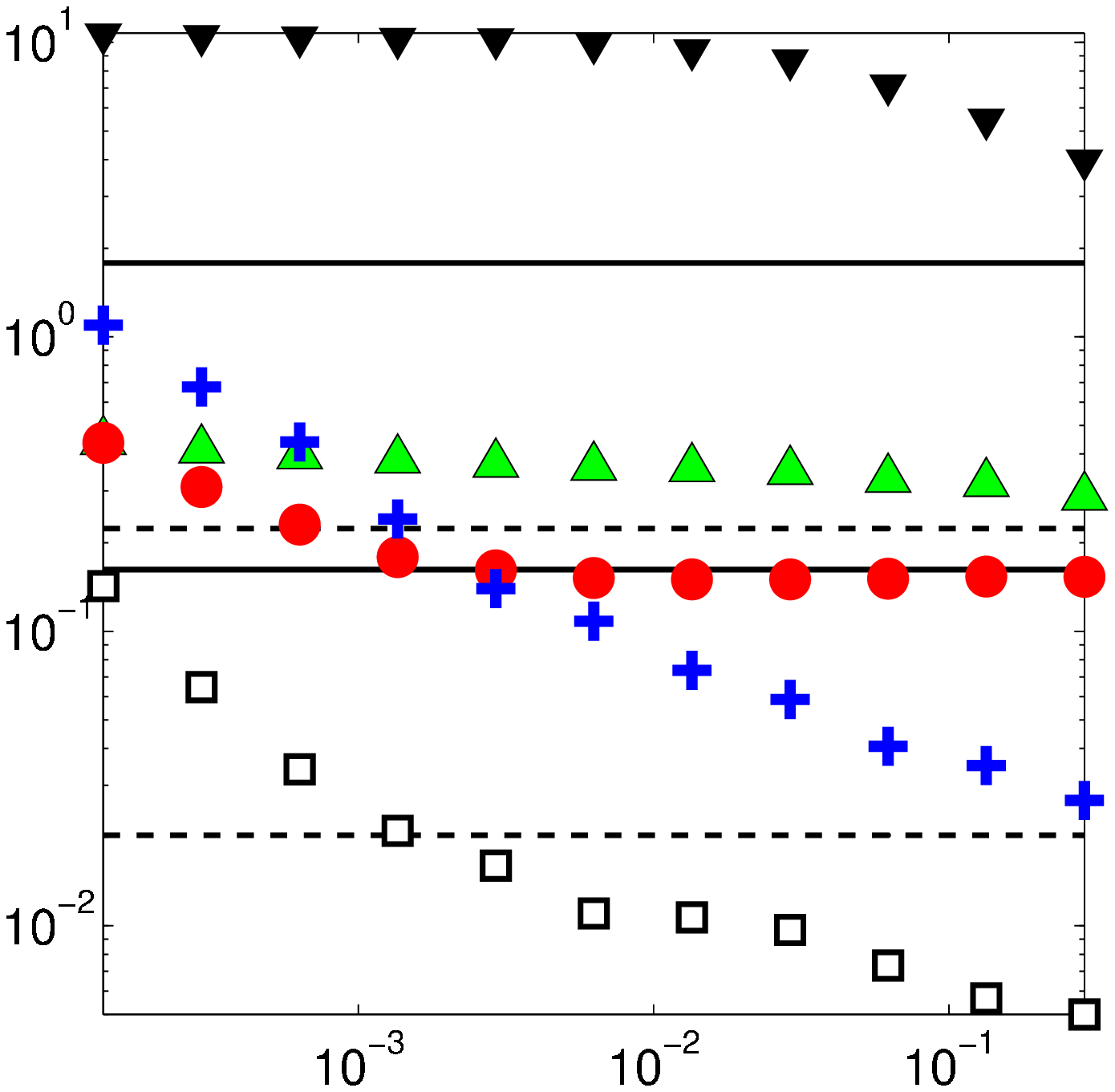}\\
\small{(E)} &\small{(E)} & \small{(F)} & \small{(F)} \\
\figc[width=41mm]{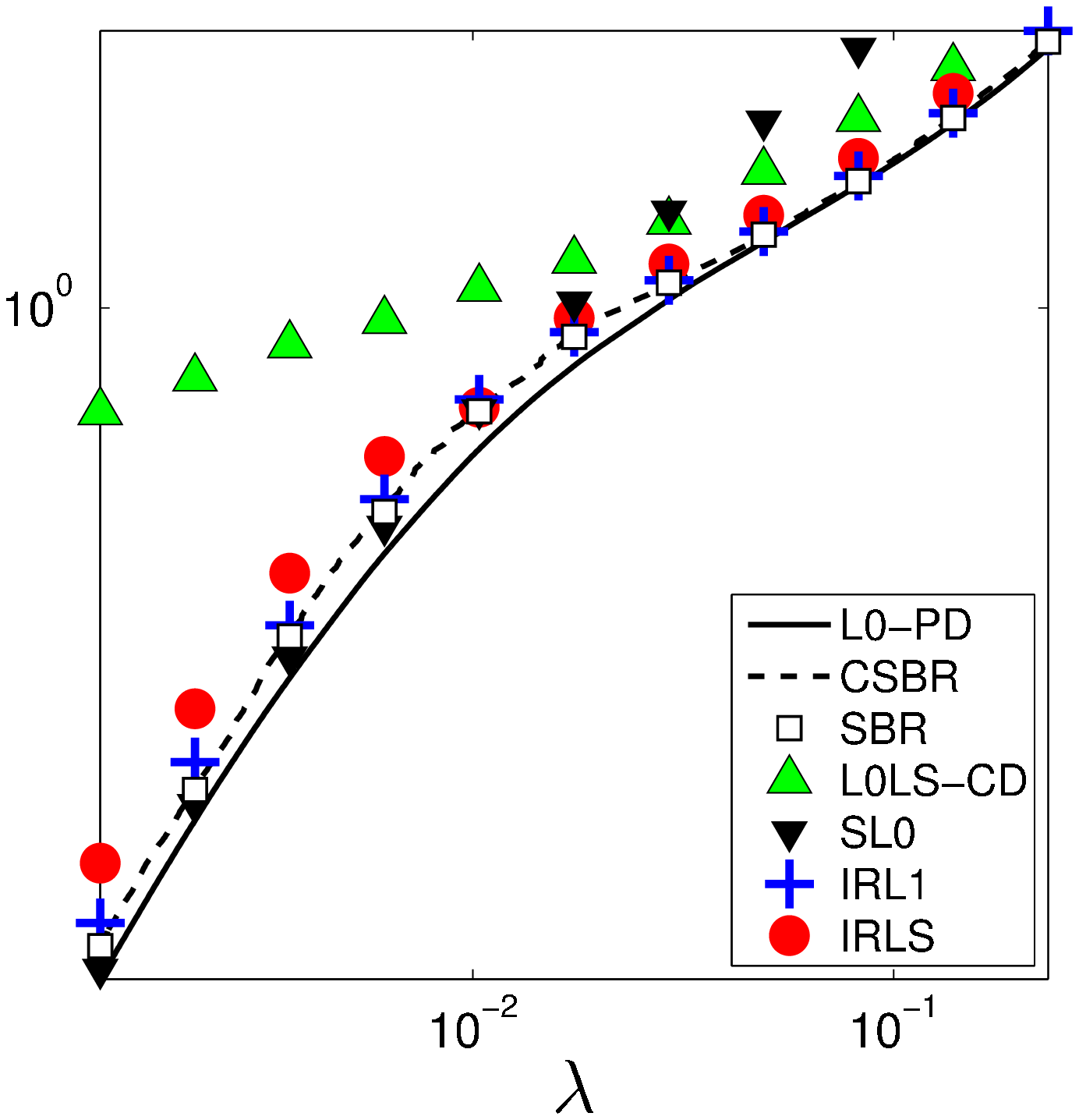}&
\figc[width=41mm]{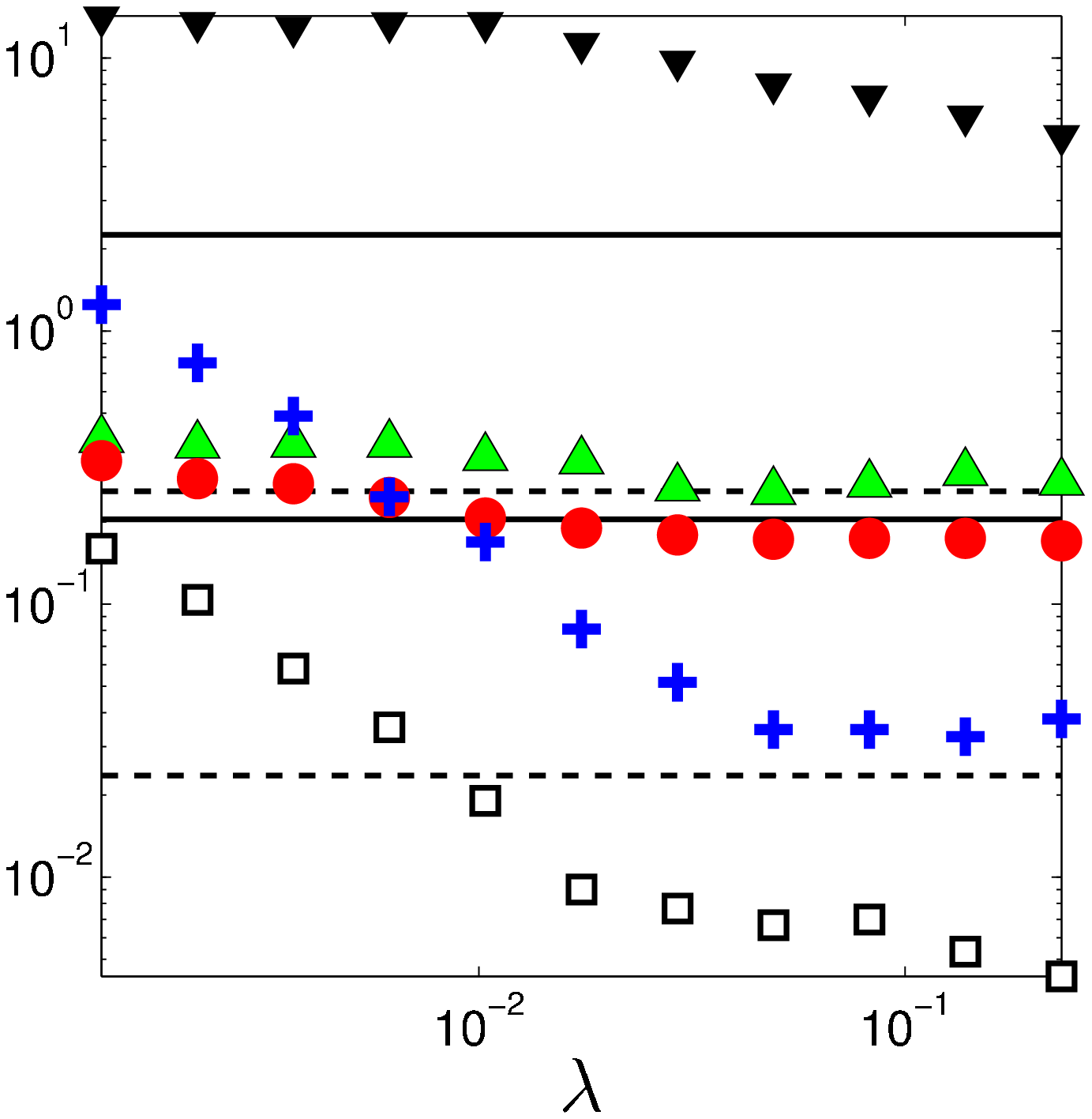}&
&
\\
\small{(G)} &
\small{(G)} & &
\end{tabular}
}
\FINAL{
\setlength{\tabcolsep}{.0cm}
\begin{tabular}{ccc}
& $\min_\xb\Jc(\xb;\lambda)$ & CPU Time (seconds)\\
\small{(E)} &\figc[width=40mm]{S8_Pb3N300fact1SNR25_figJJ_bis}&
\figc[width=40mm]{S8_Pb3N300fact1SNR25_figTT_bis}\\
\small{(F)} & 
\figc[width=40mm]{S8_Pb3N300fact1SNR25K30_figJJ_bis}&
\figc[width=40mm]{S8_Pb3N300fact1SNR25K30_figTT_bis}\\
\small{(G)} &\figc[width=41mm]{S8_Pb3N300fact1SNR10_figJJ_bis}&
\figc[width=41mm]{S8_Pb3N300fact1SNR10_figTT_bis}
\end{tabular}
}
}
\caption[Comparison of algorithms for the jump detection 
problem]{Comparison of algorithms for the jump detection 
problem for the scenarii E, F and G of Table~\ref{tab:scenarios}. }
  \label{fig:bigcompar_jump_noisy}
\end{figure}

\subsubsection{Analysis in the optimization viewpoint}
CSBR and $\ell_0$-PD are always among the most accurate to minimize
the cost function, as illustrated on
Figs.~\ref{fig:bigcompar_deconv_noisy}, \ref{fig:bigcompar_jump_noisy}
and~\ref{fig:bigcompar_noisefree}. We can clearly distinguish two
groups of algorithms on these figures: IRLS, L0LS-CD and SL0 one the
one hand, and the OLS-based algorithms (SBR, CSBR, $\ell_0$-PD) and
IR$\ell_1$ on the other hand, which are the most accurate. We cannot
clearly discriminate the accuracy of SBR and CSBR: one may behave
slightly better than the other depending on the scenarii. On the
contrary, SBR and CSBR are often outperformed by $\ell_0$-PD. The
obvious advantage of CSBR and $\ell_0$-PD over SBR and IR$\ell_1$ is
that they are $\ell_0$-homotopy algorithms, \ie a set of solutions are
delivered for many sparsity levels, and the corresponding
$\lambda$-values are adaptively found. On the contrary, the SBR output
is related to a single $\lambda$ whose tuning may be tricky. Another
advantage over IR$\ell_1$ is that the structure of forward-backward
algorithms is simpler, as no call to any $\ell_1$ solver is
required. Moreover, the number of parameters to tune is lower: there
is a single (early) stopping parameter $\lambda_{\mathrm{stop}}$.

The price to pay for a better performance is an increase of the
computation burden. On Figs.~\ref{fig:bigcompar_deconv_noisy},
\ref{fig:bigcompar_jump_noisy} and~\ref{fig:bigcompar_noisefree}, two
lines are drawn for CSBR (respectively, for $\ell_0$-PD). They are
horizontal because the algorithm is run only once per trial, so there
is a single computation time measurement.  The first line corresponds
to the overall computation time, \ie from the start to the termination
of CSBR / $\ell_0$-PD. This time is often more expensive than for
other algorithms. However, the latter times refer to a single
execution for some $\lambda_i^\mathrm{G}$ value. If one wants to
recover sparse solutions for many $\lambda_i^\mathrm{G}$'s, they must
be cumulated. This is the reason why we have drawn a second line for
CSBR and $\ell_0$-PD corresponding to a normalization (by
$N_\lambda=11$) of the overall computation time. In this viewpoint,
the CPU time of CSBR and $\ell_0$-PD are very reasonable.

\begin{table}[t]
  \caption{Jump detection problem in the noisy setting. The algorithms
    are evaluated in terms of support error (SE) and 
    number of true
    positives (TP). The number of jumps that are found 
    is reported (Order) together with the ``true order''
    corresponding to the ground truth $k$. The scores related
    to the MDLc estimate are indicated similarly.
  }
\label{tab:ERC_jumps}
\setlength{\arraycolsep}{3.5pt}
\centering
{
\setlength{\tabcolsep}{0.12cm}
\begin{tabular}{|l|r|r|r|r|r|r|r|}
\hline
\textbf{Scenario E}& $\ell_0$-PD & CSBR & SBR & $\ell_0$LS-CD & S$\ell_0$ & IR$\ell_1$ & IRLS\\
\hline
SE		 & 1.6 & 1.6 & 1.6 & 5.3 & 4.0 & 1.5 & 1.8\\
TP	 & 8.6 & 8.7 & 8.6 & 5.2 & 7.8 & 8.7 & 8.7\\
Order (true: 10) & 8.8 & 9.0 & 8.9 & 5.7 & 9.6 & 8.9 & 9.1\\
\hline
MDLc-SE	 & 4.7 & 4.3 & 4.1 & 22.7 & 5.6 & 4.1 & 3.5\\
MDLc-TP	 & 8.7 & 8.8 & 8.8 & 6.9 & 8.6 & 8.8 & 8.8\\
MDLc-Order	 & 12.2 & 11.9 & 11.6 & 26.6 & 12.7 & 11.7 & 11.0\\
\hline
\textbf{Scenario F}	& $\ell_0$-PD & CSBR & SBR & $\ell_0$LS-CD & S$\ell_0$ & IR$\ell_1$ & IRLS\\
\hline
SE		 & 11.1 & 11.9 & 11.8 & 22.5 & 11.6 & 10.9 & 11.6\\
TP	 & 21.2 & 20.6 & 20.7 & 9.2 & 20.3 & 20.8 & 20.7\\
Order (true: 30)	 & 23.6 & 23.2 & 23.2 & 10.9 & 22.2 & 22.5 & 23.1\\
\hline
MDLc-SE	 & 13.7 & 13.4 & 13.4 & 39.2 & 14.0 & 13.1 & 13.3\\
MDLc-TP	 & 21.8 & 21.8 & 21.4 & 12.9 & 21.6 & 22.1 & 21.5\\
MDLc-Order	 & 27.3 & 27.0 & 26.3 & 35.0 & 27.2 & 27.2 & 26.4\\
\hline
\textbf{Scenario G}& $\ell_0$-PD & CSBR & SBR & $\ell_0$LS-CD & S$\ell_0$ & IR$\ell_1$ & IRLS\\
\hline
SE		 & 7.3 & 7.5 & 7.5 & 8.9 & 10.3 & 7.2 & 7.5\\
TP	 & 4.0 & 3.6 & 3.6 & 3.1 & 2.9 & 3.9 & 4.0\\
Order (true: 10)	 & 5.23 & 4.73 & 4.73 & 5.17 & 6.07 & 4.97 & 5.57\\
\hline
MDLc-SE	 & 11.4 & 10.7 & 10.9 & 11.7 & 15.1 & 11.2 & 10.7\\
MDLc-TP	 & 4.2 & 4.2 & 4.2 & 3.0 & 3.9 & 4.2 & 4.4\\
MDLc-Order & 9.8 & 9.1 & 9.3 & 7.6 & 12.8 & 9.6 & 9.5\\
\hline
\end{tabular}
}
\end{table}

The computation time depends on many factors among which the
implementation of algorithms (including the memory storage) and the
chosen stopping rules. We have followed an homogeneous implementation
of algorithms to make the CPU time comparisons meaningful. We have
defined two sets of stopping rules depending on the problem
dimension. The default parameters apply to medium size problems
($m=300$). They are relaxed for problems of larger dimension ($m>500$)
to avoid huge computational costs. The stopping rule of CSBR and
$\ell_0$-PD is always
$\lambda\leq\lambda_{\mathrm{stop}}=\alpha\lambda_1^\mathrm{G}$ with
$\alpha=1$ for CSBR and 0.5 (medium size) or 0.8 (large size) for
$\ell_0$-PD. For L0LS-CD, the maximum number of cyclic descents
(update of every amplitude $x_i$) is set to 60 or 10 depending on the
dimension. For SL0, we have followed the default setting
of~\cite{Mohimani09} for the rate of deformation of the nonconvex
penalty. The number of BFGS iterations done in the local minimization
steps for each penalty is set to $L=$ 40 or 5. It is set to $5L$ for
the last penalty which is the most nonconvex. Regarding IRLS and
IRL$\ell_1$, we keep the same settings whatever the dimension since
the computation times remain reasonable for large dimensions. Finally,
SBR does not require any arbitrary stopping rule. The problems of
large dimensions correspond to scenarii C and D. We observe on
Fig.~\ref{fig:bigcompar_deconv_noisy} that the comparison (trade-off
performance \vs computation time) is now clearly in favor of CSBR and
$\ell_0$-PD. IR$\ell_1$ remains very competitive although the average
numerical cost becomes larger.

\begin{figure}[!t]
\DRAFT{\setlength{\tabcolsep}{.7cm}
\begin{center}
\begin{tabular}{ccc}
\figc[width=40mm]{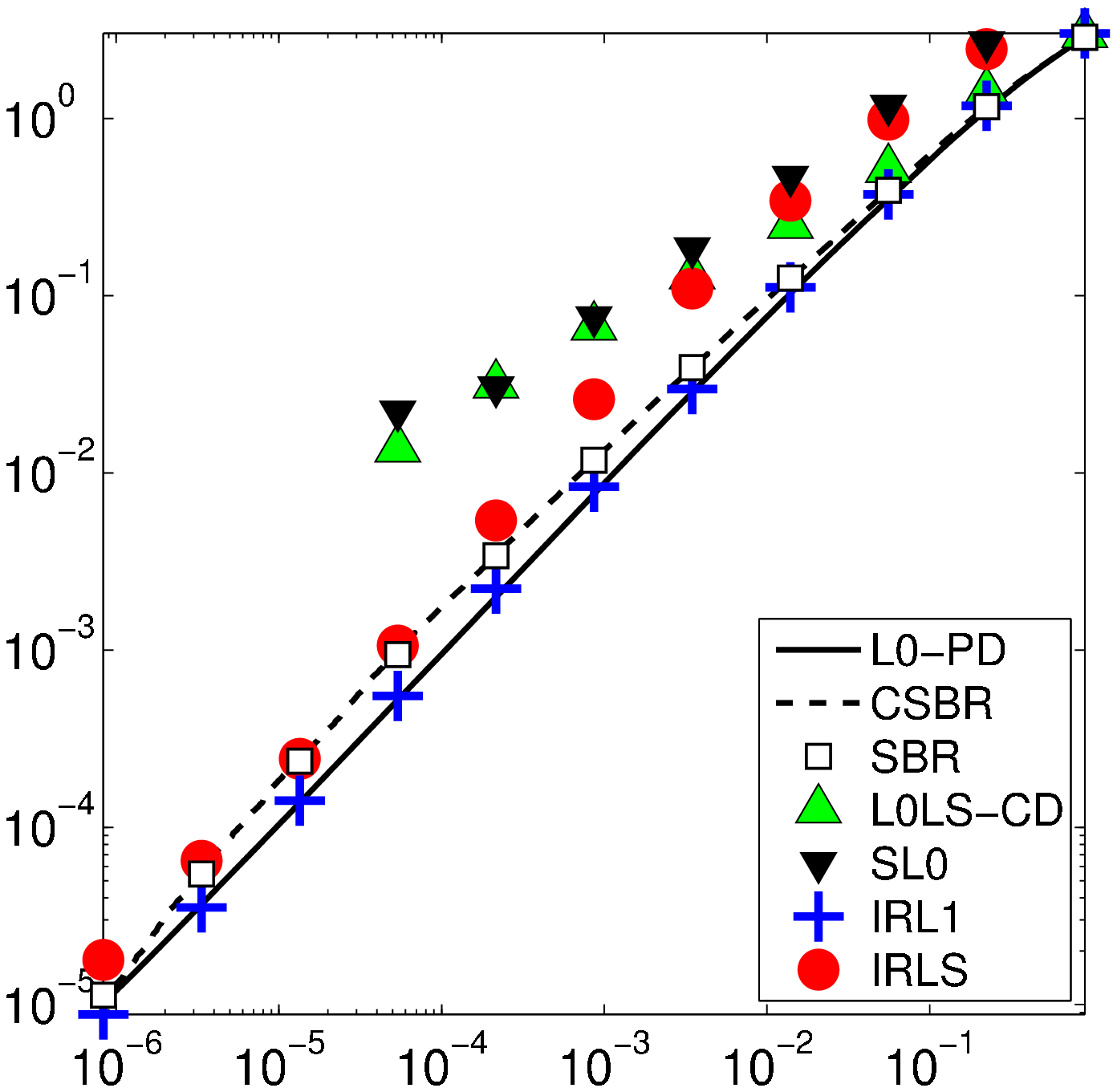}&
\figc[width=40mm]{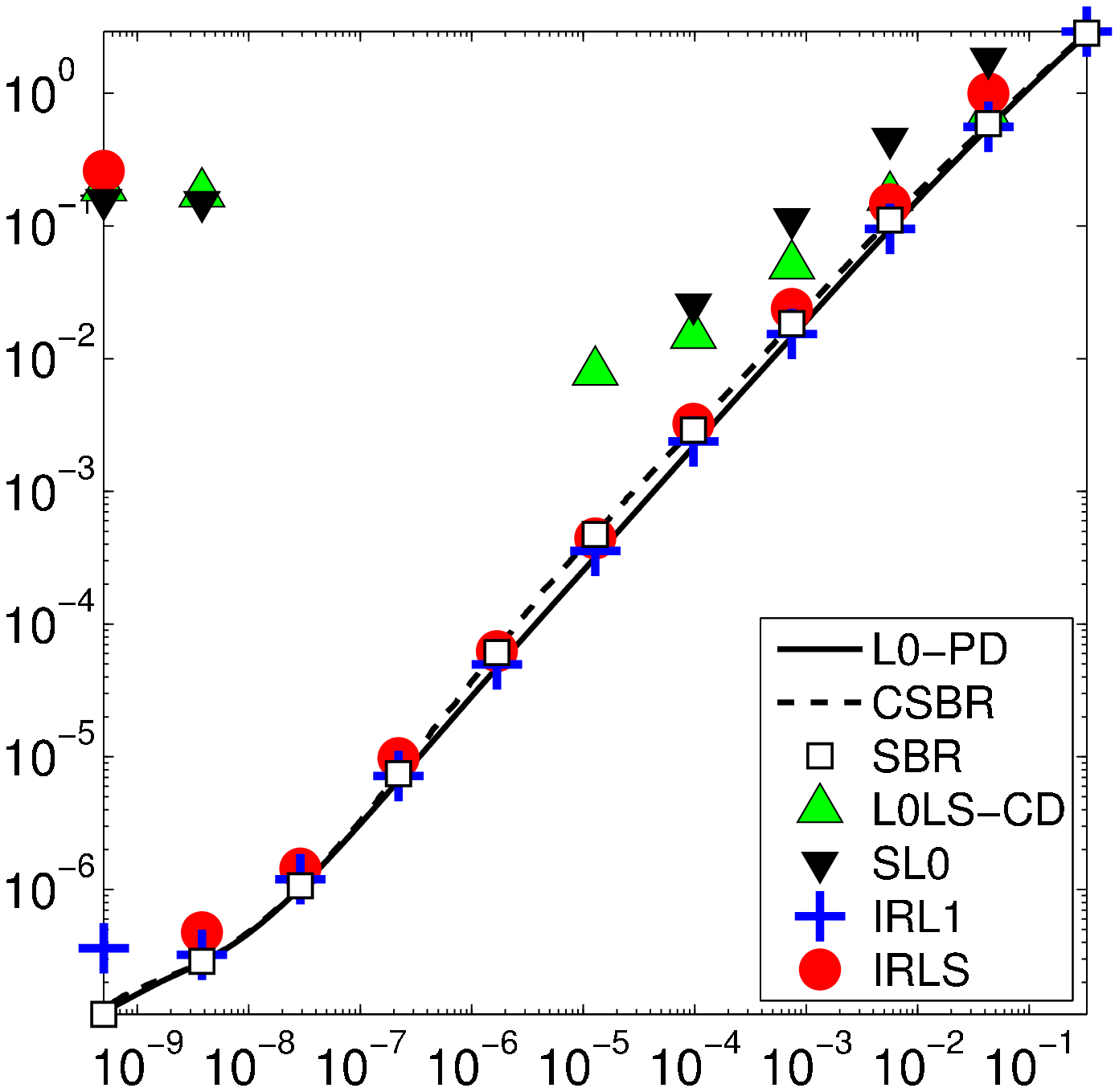}&
\figc[width=41mm]{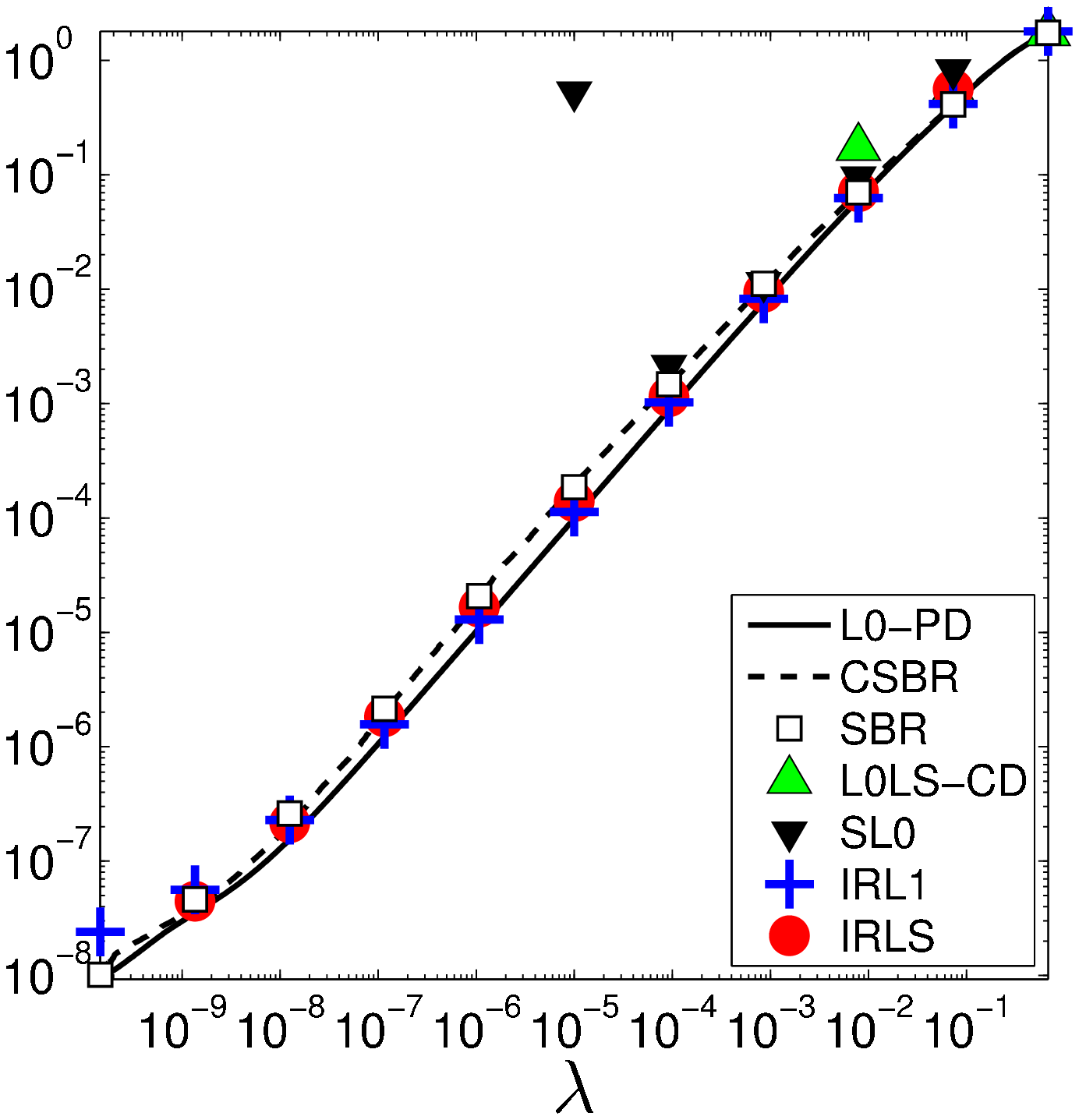}\\
\figc[width=40mm]{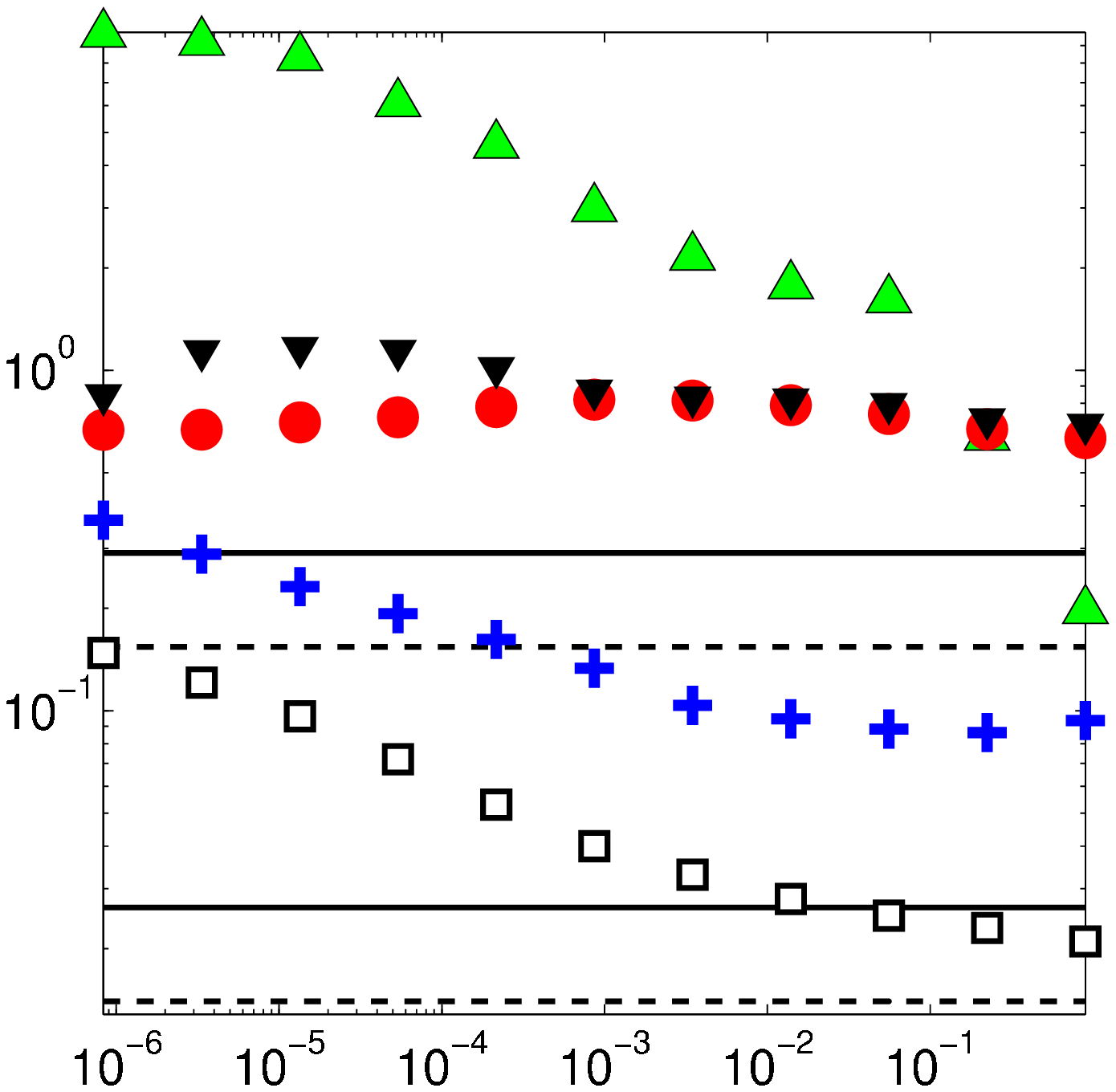}&
\figc[width=40mm]{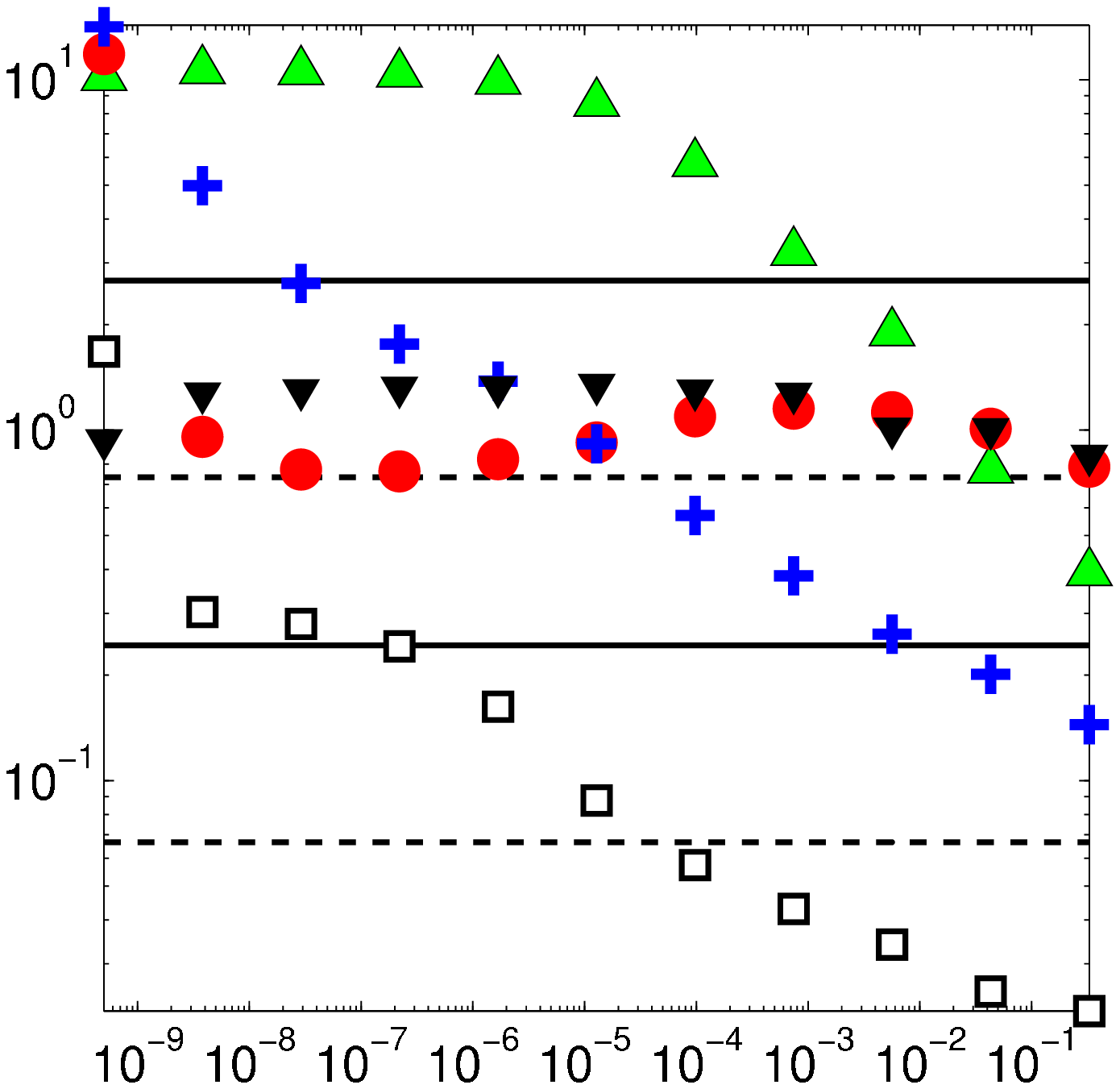}&
\figc[width=41mm]{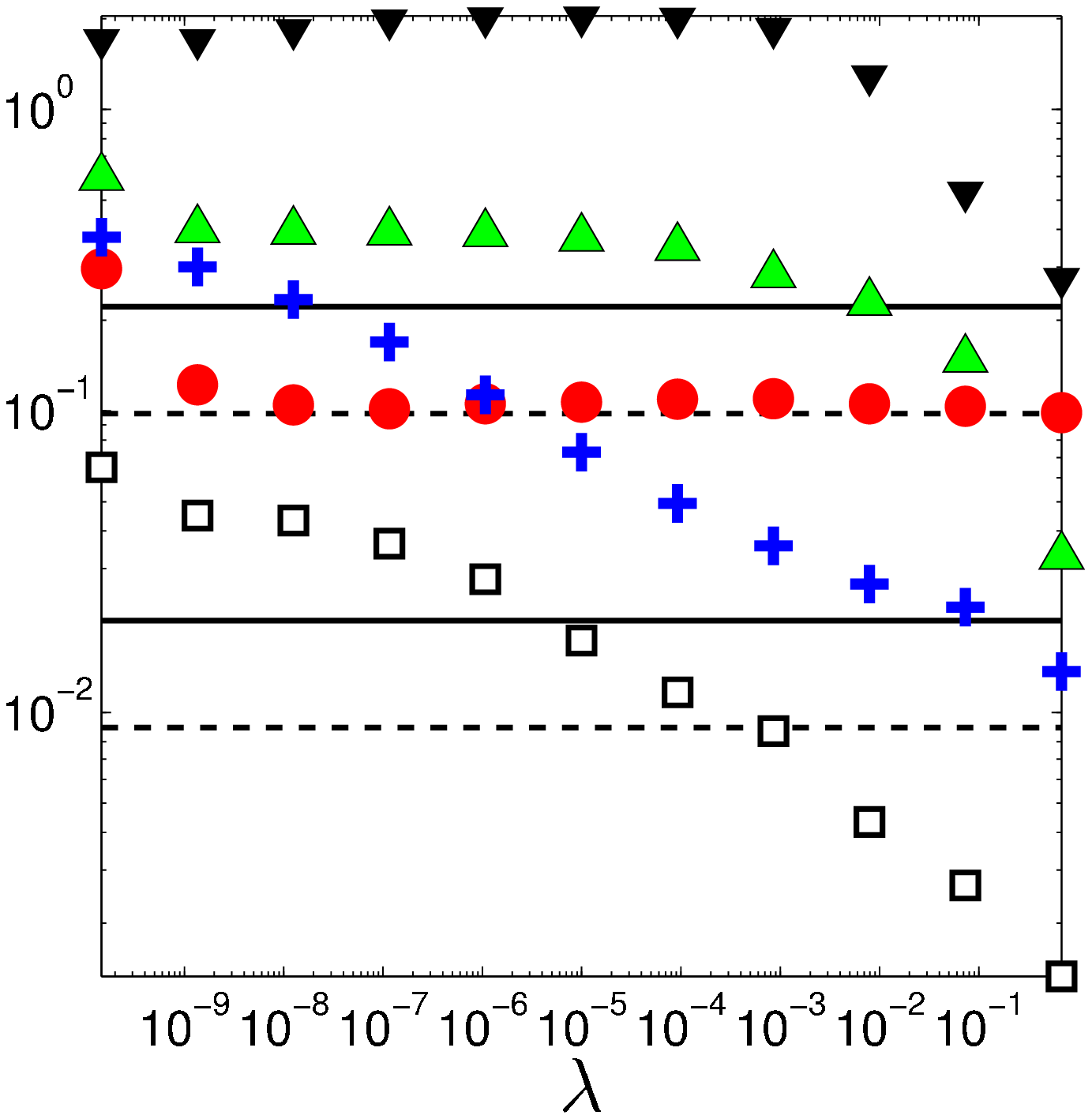}\\
(H)&(I)&(J)
\end{tabular}
\end{center}
}
\FINAL{
\centering
\setlength{\tabcolsep}{.0cm}
\begin{tabular}{ccc}
(H)&\figc[width=40mm]{S8_Pb5N756fact3SNR80_figJJ_bis}&
\figc[width=40mm]{S8_Pb5N756fact3SNR80_figTT_bis}\\
(I)&\figc[width=40mm]{S8_Pb5N756fact3SNR80K30_figJJ_bis}&
\figc[width=40mm]{S8_Pb5N756fact3SNR80K30_figTT_bis}\\
(J)&\figc[width=41mm]{S8_Pb6N252fact1SNR80_figJJ_bis}&
\figc[width=41mm]{S8_Pb6N252fact1SNR80_figTT_bis}
\end{tabular}
}
\caption[Comparison of algorithms for the noise-free deconvolution
problem ]{Comparison of algorithms for the noise-free deconvolution
  problem, \ie for the scenarii H, I and J of
  Table~\ref{tab:scenarios}. Some markers do not appear for low
  $\lambda$'s (L0LS-CD, SL0) in the left figures because they do not
  lay in the zoom-in window (their performance is poor).  }
  \label{fig:bigcompar_noisefree}
\end{figure}

\subsubsection{Analysis in the support recovery viewpoint}
The support recovery performance is only shown for the scenarii E to J
(Tabs.~\ref{tab:ERC_jumps} and~\ref{tab:ERC_noisefree}). For noisy
deconvolution problems, these results are omitted because the support
error is often quite large and the true positive scores are low
whatever the algorithm, especially for scenarii B to D. Specifically,
the least support error always exceeds 20, 10, 10 and 32 for the
scenarii A to D ($k=30$, 10, 10 and 30, respectively). For such
difficult problems, one can hardly discriminate algorithms based on
simple binary tests such as the true positive rate. More sophisticated
localization tests are non binary and would take into account the
distance between the location of the true spikes and their wrong
estimates~\cite{VanRossum01}. It is noticeable, though, that the MDLc
estimator delivers subsets of realistic cardinality for scenarii A to
D (\eg the subsets found with CSBR are of cardinalities 33, 9, 15 and
38, the true cardinalities being 30, 10, 10 and 30). The model orders
are also quite accurate for the noisy jump detection problem
(Tab.~\ref{tab:ERC_jumps}) whereas the true support is often partially
detected by several of the considered algorithms. Here, CSBR and
$\ell_0$-PD are among the best algorithms in terms of support error.

The results of Tab.~\ref{tab:ERC_noisefree} and
Fig.~\ref{fig:bigcompar_noisefree} correspond to the deconvolution
problem in noise-free case. The data \yb are undersampled so that the
dictionary \Ab is overcomplete. The undersampling rate $\Delta\approx
m/n$ is set to 2 in scenarii H and I and 4 in scenario J. Again, CSBR
and $\ell_0$-PD are among the best (SE, TP, MDLc-order)
especially for the most difficult problem J.

\begin{table}[t]
  \caption{Sparse deconvolution problem in the noise-free setting:
    exact support recovery.
  }
\label{tab:ERC_noisefree}
\setlength{\arraycolsep}{3.5pt}
\centering
{
\setlength{\tabcolsep}{0.11cm}
\begin{tabular}{|l|r|r|r|r|r|r|r|}
\hline
\textbf{Scenario H}& $\ell_0$-PD & CSBR & SBR & $\ell_0$LS-CD & S$\ell_0$ & IR$\ell_1$ & IRLS\\
\hline
SE		 & 2.5 & 3.6 & 4.8 & 11.4 & 13.0 & 0.8 & 6.1\\
TP	 & 8.3 & 8.2 & 6.8 & 0.4 & 0.1 & 9.5 & 9.4\\
Order (true: 10)	 & 9.1 & 10.0 & 8.3 & 2.2 & 3.2 & 9.8 & 14.9\\
\hline
MDLc-SE	 & 3.6 & 3.8 & 5.8 & 168.5 & 343.8 & 1.1 & 9.0\\
MDLc-TP	 & 8.6 & 8.6 & 7.9 & 3.3 & 6.6 & 9.5 & 9.6\\
MDLc-Order	 & 10.8 & 11.0 & 11.6 & 153.5 & 347.0 & 10.1 & 18.2\\
\hline
\hline
\textbf{Scenario I}& $\ell_0$-PD & CSBR & SBR & $\ell_0$LS-CD & S$\ell_0$ & IR$\ell_1$ & IRLS\\
\hline
SE		 & 0.9 & 1.3 & 2.1 & 36.7 & 48.5 & 3.8 & 9.4\\
TP	 & 29.4 & 29.3 & 29.1 & 0.7 & 0.8 & 28.0 & 27.7\\
Order (true: 30)	 & 29.7 & 29.8 & 30.2 & 8.2 & 20.1 & 29.8 & 34.8\\
\hline
MDLc-SE	 & 3.8 & 3.5 & 3.7 & 686.0 & 444.9 & 9.5 & 114.3\\
MDLc-TP	 & 29.5 & 29.4 & 29.2 & 28.6 & 17.5 & 28.5 & 26.4\\
MDLc-Order	 & 32.8 & 32.3 & 32.1 & 437.0 & 449.8 & 36.5 & 137.2\\
\hline
\hline
\textbf{Scenario J}& $\ell_0$-PD & CSBR & SBR & $\ell_0$LS-CD & S$\ell_0$ & IR$\ell_1$ & IRLS\\
\hline
SE		 & 0.3 & 3.5 & 5.3 & 10.3 & 10.4 & 2.4 & 4.3\\
TP	 & 9.8 & 7.3 & 5.8 & 0.6 & 2.6 & 8.8 & 9.2\\
Order (true: 10)	 & 9.8 & 8.1 & 6.9 & 1.4 & 5.6 & 10.0 & 12.7\\
\hline
MDLc-SE	 & 2.6 & 7.7 & 12.4 & 176.2 & 78.6 & 7.2 & 69.0\\
MDLc-TP	 & 9.7 & 8.9 & 8.0 & 8.6 & 3.0 & 8.9 & 4.1\\
MDLc-Order	 & 12.0 & 15.5 & 18.5 & 73.0 & 74.6 & 14.9 & 67.2\\
\hline
\end{tabular}
}
\end{table}

\subsubsection{Overcomplete dictionaries with noise}
We now provide arguments indicating that the proposed algorithms are
competitive as well for noisy problems with overcomplete
dictionaries. The detailed experiments commented below are not
reported for space reasons.

We have first considered the noisy deconvolution problem with
$\Delta=2$ or 4 leading to overcomplete dictionaries, the other
parameters being set as in scenarii A to D. Although the data
approximation is qualitatively good for CSBR and $\ell_0$-PD, the SE
and TP scores are very weak. It is hard to discriminate the
performance of algorithms because these measures are very weak for all
considered algorithms. Moreover, the values of $\Jc(\lambda)$ found
for most algorithms are often similar.

We have also considered an adaptive spline approximation problem
generalizing the jump detection problem to the approximation of a
signal using piecewise polynomials of degree $P=1$ or
2~\cite{Soussen11c}. The jump detection problem can indeed be thought
of as the approximation with a piecewise constant signal ($P=0$). The
generalized version~\cite{Soussen11c} is inspired from the regression
spline modeling in~\cite{Friedman91}. Now, the dictionary atoms are
related to the detection of the locations of jumps, changes of slopes
and changes of curvatures in the signal \yb (subdictionaries $\Ab^0$,
$\Ab^1$ and $\Ab^2$). The dictionary then takes the form
$\Ab\leftarrow\stdcro{\Ab^0,\Ab^1}$ or
$\Ab\leftarrow\stdcro{\Ab^0,\Ab^1,\Ab^2}$ where each sub-dictionary
$\Ab^p$ ($p\leq P$) is formed of shifted versions of the one-sided
power function $i\mapsto\stdcro{\max(i,0)}^p$. The size of the full
dictionary $\Ab$ is approximately $m\times (P+1)m$. Hence, it becomes
overcomplete as soon as $P\geq 1$. We have shown~\cite{Soussen11c}
that SBR is competitive when $P=1$ or 2. We have carried out new tests
confirming that CSBR and $\ell_0$-PD are more efficient than their
competitors in terms of values of $\Jc(\lambda)$. However, the rate of
true positives is low for $P\geq 1$ since the location of the change
of slopes and of curvatures can hardly be exactly recovered from noisy
data.

\section{Software}
The Matlab implementation of the proposed CSBR and $\ell_0$-PD
algorithms is available at\\
\url{www.cran.univ-lorraine.fr/perso/charles.soussen/software.html}
including programs showing how to call these functions. 

\section{Conclusion}
The choice of a relevant sparse approximation algorithm relies on a
trade-off between the desired performance and the computation time one
is ready to spend. The proposed algorithms are relatively expensive
but very well suited to inverse problems inducing highly correlated
dictionaries. A reason is that they have the capacity to escape from
local minimizers of
$\Jc(\xb;\lambda)=\|\yb-\Ab\xb\|_2^2+\lambda\|\xb\|_0$~\cite{Soussen11c}.
This behavior is in contrast with other classical sparse algorithms.

We have shown the usefulness and efficiency of the two SBR extensions
when the level of sparsity is moderate to high, \ie $k/\min(m,n)$ is
lower than 0.1. They remain competitive when $k/\min(m,n)$ ranges
between 0.1 and 0.2, and their performance gradually degrade for
weaker levels of sparsity, which is an expected behavior for such
greedy type algorithms.
For a single $\lambda$, CSBR is as efficient as SBR, and $\ell_0$-PD
improves the SBR and CSBR performance within a larger computation
cost. The main benefit over SBR is that 
sparse solutions are provided for a continuum of $\lambda$-values, enabling
the utilization of any classical order selection method. We found that
the MDL criterion yields very accurate estimates of the cardinality
$\|\xb\|_0$.

Our perspectives include the proposal of forward-backward search
algorithms that will be faster than SBR and potentially more
efficient. In the standard version of SBR, CSBR and $\ell_0$-PD, a
single replacement refers to the insertion or removal of a dictionary
element. The cost of an iteration is essentially related to the $n$
linear system resolutions done to test single replacements for all
dictionary atoms. The proposed algorithms obviously remain valid when
working with a larger neighborhood, \eg when testing the replacement
of two atoms simultaneously, but their complexity becomes huge. To
avoid such numerical explosion, one may rather choose not to carry out
all the replacement tests, but only some tests that are likely to be
effective. Extensions of OMP and OLS were recently proposed in this
spirit~\cite{Chatterjee12} and deserve consideration for proposing
efficient forward-backward algorithms.

\appendices

\section{Properties of the $\ell_0$ regularization paths}
\label{app:theory}
In this appendix, we prove that the $\ell_0$-penalized 
path $\Sc^\star_{\mathrm{P}}$ (see Definition~\ref{def:Cpath}) is
piecewise constant (Theorem~\ref{th:1}) and is a subset of the
$\ell_0$-constrained regularization path $\Sc^\star_{\mathrm{C}}$
(Theorem~\ref{th:2complet}).
We will denote the $\ell_0$-curve by
$\lambda\mapsto\Jc^\star(\lambda)=\min_S\{\hat{\Jc}(S;\lambda)\}$.
Let us recall that 
this function is concave and affine on each
interval $(\lambda_{i+1}^\star,\lambda_i^\star)$, with
$i\in\{0,\ldots,I\}$ (Definition~\ref{def:l0curve}). Moreover,
$\lambda_{I+1}^\star=0$ and $\lambda_{0}^\star=+\infty$.

\subsection{Proof of Theorem~\ref{th:1}}
\label{sec:piece_const}
We prove Theorem~\ref{th:1} together with the following lemma, which
is informative about the content of $\Sc^\star_{\mathrm{P}}(\lambda)$
for the breakpoints $\lambda=\lambda_i^\star$.
\begin{lemma}
  \label{lem:2}
  Let $i\in\{1,\ldots,I-1\}$. Then, for all
  $\lambda\in(\lambda_{i+1}^\star,\lambda_{i}^\star)$,
  $\Sc^\star_{\mathrm{P}}(\lambda)\subset\Sc^\star_{\mathrm{P}}(\lambda^\star_{i+1})\cap
  \Sc^\star_{\mathrm{P}}(\lambda^\star_i)$.

For the first and last intervals, we have:
\begin{itemize}
\item For all $\lambda\in(0,\lambda_{I}^\star)$, 
  $\Sc^\star_{\mathrm{P}}(\lambda)\subset\Sc^\star_{\mathrm{P}}(\lambda^\star_{I})$.
\item For all $\lambda\in(\lambda_{1}^\star,+\infty)$, 
  $\Sc^\star_{\mathrm{P}}(\lambda)=\{\emptyset\}\subset\Sc^\star_{\mathrm{P}}(\lambda_1^\star)$.
\end{itemize}
\end{lemma}
\begin{IEEEproof}[Proof of Theorem~\ref{th:1}]
  By definition, the $\ell_0$-curve
  is the concave envelope of the (finite) set of lines $S$ for all
  possible subsets $S$. Because it is affine on the $i$-th interval
  $(\lambda^\star_{i+1},\lambda^\star_{i})$, $\Jc^\star(\lambda)$
  coincides with
  $\hat{\Jc}({S_i};\lambda)=\Ec(S_i)+\lambda\stdbars{S_i}$, where
  $S_i$ is some optimal subset for all
  $\lambda\in(\lambda^\star_{i+1},\lambda^\star_{i})$.

  Let $\lambda\in(\lambda^\star_{i+1},\lambda^\star_{i})$ and
  $S\in\Sc^\star_{\mathrm{P}}(\lambda)$. Then,
  $\hat{\Jc}(S;\lambda)=\hat{\Jc}({S_i};\lambda)$. It follows that
  both lines $S$ and $S_i$ necessarily coincide; otherwise, they would
  intersect at $\lambda$, and line $S$ would lay below $S_i$ on either
  interval $(\lambda^\star_{i+1},\lambda)$ or
  $(\lambda,\lambda^\star_{i})$, which contradicts the definition of
  $S_i$. We conclude that $S\in\Sc^\star_{\mathrm{P}}(\lambda')$ for
  all $\lambda'\in(\lambda^\star_{i+1},\lambda^\star_{i})$.

  We have shown that the content of $\Sc^\star_{\mathrm{P}}(\lambda)$
  does not depend on $\lambda$ when
  $\lambda\in(\lambda^\star_{i+1},\lambda^\star_{i})$.
\end{IEEEproof}
\begin{IEEEproof}[Proof of Lemma~\ref{lem:2}]
  The first result
  $\Sc^\star_{\mathrm{P}}(\lambda)\subset\Sc^\star_{\mathrm{P}}(\lambda^\star_{i+1})\cap
  \Sc^\star_{\mathrm{P}}(\lambda^\star_i)$ is obtained by slightly
  adapting the proof of Theorem~\ref{th:1}: replace
  $(\lambda^\star_{i+1},\lambda^\star_{i})$ by the closed interval
  $[\lambda^\star_{i+1},\lambda^\star_{i}]$, and set $\lambda'$ to
  both endpoints of this interval.

  The second and third results are obtained similarly, by considering
  the intervals $(0,\lambda^\star_{I}]$ and
  $[\lambda^\star_{1},+\infty)$, and setting $\lambda'\leftarrow
  \lambda^\star_{I}$ and $\lambda'\leftarrow \lambda^\star_{1}$,
  respectively. It is obvious that $\Sc^\star_{\mathrm{P}}(\lambda)$
  reduces to the empty support for $\lambda>\lambda_1^\star$ since the
  $\ell_0$-curve is constant for $\lambda>\lambda_1^\star$.
\end{IEEEproof}

\subsection{Proof of Theorem~\ref{th:2complet}}
\label{sec:const_pen}
  The first result is straightforward: for any $\lambda$ and for
  $S\in\Sc^\star_{\mathrm{P}}(\lambda)$, we have $S\in\Sc^\star_{\mathrm{C}}(\stdbars{S})$.
  Otherwise, there would exist $S'$ with
  $\stdbars{S'}\leq\stdbars{S}$ and $\Ec(S')<\Ec(S)$. Then,
  $\hat{\Jc}({S'};\lambda)<\hat{\Jc}({S};\lambda)$ would contradict
  $S\in\Sc^\star_{\mathrm{P}}(\lambda)$.

  To prove the second result, let us first show that for any $i$,
  $\exists
  k_i:\,\forall\lambda\in(\lambda^\star_{i+1},\lambda^\star_{i}),\,
  \Sc^\star_{\mathrm{P}}(\lambda)\subset\Sc^\star_{\mathrm{C}}(k_i)$.

  Let $S\in\Sc^\star_{\mathrm{P}}(\lambda)$ for some
  $\lambda\in(\lambda^\star_{i+1},\lambda^\star_{i})$.
  Theorem~\ref{th:1} implies that
  $S\in\Sc^\star_{\mathrm{P}}(\lambda)$ for \emph{any}
  $\lambda\in(\lambda^\star_{i+1},\lambda^\star_{i})$. Therefore,
  $\Jc^\star(\lambda)=\hat{\Jc}(S;\lambda)$ for
  $\lambda\in(\lambda^\star_{i+1},\lambda^\star_{i})$, and the slope of
  line $S$, \ie $\stdbars{S}$, is constant whatever
  $S\in\Sc^\star_{\mathrm{P}}(\lambda)$ and
  $\lambda\in(\lambda^\star_{i+1},\lambda^\star_{i})$. Let us denote
  this constant by $k_i=\stdbars{S}$. According to the first
  paragraph of the proof, $S\in\Sc^\star_{\mathrm{P}}(\lambda)$ implies that
  $S\in\Sc^\star_{\mathrm{C}}(k_i)$.

  Let us prove the reverse inclusion
  $\Sc^\star_{\mathrm{C}}(k_i)\subset\Sc^\star_{\mathrm{P}}(\lambda)$. Let
  $\lambda\in(\lambda^\star_{i+1},\lambda^\star_{i})$ and
  $S\in\Sc^\star_{\mathrm{C}}(k_i)$. First, we have $\stdbars{S}\leq
  k_i$. Second, for any $S'\in\Sc^\star_{\mathrm{P}}(\lambda)$,
  we have $\stdbars{S'}=k_i$ by definition of $k_i$. We also
  have that $\Ec(S')=\Ec(S)$ because
  $\Sc^\star_{\mathrm{P}}(\lambda)\subset\Sc^\star_{\mathrm{C}}(k_i)$. Finally,
  $\hat{\Jc}(S';\lambda) \geq \hat{\Jc}({S};\lambda)$.
  $S'\in\Sc^\star_{\mathrm{P}}(\lambda)$ implies that
  $S\in\Sc^\star_{\mathrm{P}}(\lambda)$.
  This concludes the proof of the second result.


\end{document}